\documentclass[review]{elsarticle}
\makeatletter
\def\ps@pprintTitle{%
 \let\@oddhead\@empty
 \let\@evenhead\@empty
 \def\@oddfoot{\centerline{\thepage}}%
 \let\@evenfoot\@oddfoot}
\makeatother

\usepackage{graphicx}
\usepackage{amsfonts, amsmath, amssymb}
\usepackage{bm}
\usepackage{cancel}
\usepackage{color}
\usepackage{enumerate}
\usepackage{float}
\usepackage{hyperref}
\usepackage[english]{babel}
\usepackage[margin=1in]{geometry}
\usepackage{mathtools}
\usepackage{setspace}
\usepackage{amsthm}
\usepackage{subfigure}
\usepackage[mathscr]{euscript}
\usepackage{lineno}
\usepackage{array,multirow}
\newcolumntype{K}[1]{>{\centering\arraybackslash}p{#1}}
\newcommand{\RN}[1]{%
	\textup{\uppercase\expandafter{\romannumeral#1}}%
}

\newtheorem{mylem}{Lemma}

\newcommand \D [2]{\frac{\partial #1}{\partial #2}}

\renewcommand{\vec}[1]{\bm{\mathrm{#1}}}
\newcommand{\V}[1]{\bm{\mathrm{#1}}}

\def \div{\nabla \cdot \mbox{}}
\def \grad{\nabla}
\def \lap{\nabla^2}

\def \x{\vec{x}}

\def \u{\vec{u}}

\def \F{\vec{F}}

\def \U{\vec{U}}

\def \cA{\vec{\mathcal{A}}}
\def \cG{\vec{\mathcal{G}}}
\def \cW{\vec{\mathcal{W}}}
\def \C{\vec{C}}

\def \F{\vec{F}}

\def \Comegap{\V{\C}_{\Omega_b^{+}}}

\def \cP{{\mathcal{P}}}

\def \U{\vec{U}}

\def \X{\vec{X}}

\def \vcL{\vec{\mathcal{L}}}
\def \cS{\vec{\mathcal{S}}}
\def \cSIB{ \vec{\mathcal{S}}_{\text{IB}} }
\def \cSMLS{ \vec{\mathcal{S}}_{\text{MLS}} }
\def \cJ{\vec{\mathcal{J}}}
\def \cJIB{ \vec{\mathcal{J}}_{\text{IB}} }
\def \cJMLS{ \vec{\mathcal{J}}_{\text{MLS}} }

\def \g{\vec{g}}
\def \f{\vec{f}}

\def \half{\frac{1}{2}}
\def \3half{\frac{3}{2}}

\def \nref{n_{\text{ref}}}

\def \s{\vec{s}}

\def \u{\vec{u}}

\def \x{\vec{x}}

\def \div{\nabla \cdot \mbox{}}
\def \grad{\nabla}
\def \lap{\nabla^2}

\def \dt{\Delta t}
\def \dx{\Delta x}
\def \dy{\Delta y}
\def \dz{\Delta z}

\def \Ds{{\mathrm d}\s}

\def \Dx{{\mathrm d}\x}

\def \ds{\Delta s}
\def \dt{\Delta t}
\def \dx{\Delta x}

\def \ncycles{n_\text{cycles}}

\def \InterpSumi{{ \sum_{i = 1}^N}}
\def \InterpSumj{{ \sum_{j = 1}^m}}
\makeatletter
\renewcommand*\env@matrix[1][\arraystretch]{%
  \edef\arraystretch{#1}%
  \hskip -\arraycolsep
  \let\@ifnextchar\new@ifnextchar
  \array{*\c@MaxMatrixCols c}}
\makeatother

\newcommand{\REVIEW}[1]{{#1}}

\begin{document}

\let\today\relax

\begin{frontmatter}

\title{A one-sided direct forcing immersed boundary method using moving least squares}

\author[Riken]{Rahul Bale\corref{mycorrespondingauthor}}
\ead{rahul.bale@riken.jp}
\author[SDSU]{Amneet Pal Singh Bhalla\corref{mycorrespondingauthor}}
\ead{asbhalla@sdsu.edu}
\author[UNC]{Boyce E. Griffith}
\author[Riken,KU]{Makoto Tsubokura}

\address[Riken]{Riken Center for Computational Sciences, Japan}
\address[KU]{Graduate School of System Informatics, Kobe University, Kobe, Japan}
\address[UNC]{Departments of Mathematics, Applied Physical Sciences, and Biomedical Engineering, University of North Carolina, Chapel Hill, NC}
\address[SDSU]{Department of Mechanical Engineering, San Diego State University, San Diego, CA}
\cortext[mycorrespondingauthor]{Corresponding author}

\begin{abstract}
This paper presents a one-sided immersed boundary (IB) method using kernel functions constructed via a moving least squares (MLS) method. The resulting kernels effectively couple structural degrees of freedom to fluid variables on only one side of the fluid-structure interface. This reduces spurious feedback forcing and internal flows that are typically observed in IB models that use isotropic kernel functions to couple the structure to fluid degrees of freedom on both sides of the interface. The method developed here extends the original MLS methodology introduced by Vanella and Balaras (J Comput Phys, 2009). Prior IB/MLS methods have used isotropic kernel functions that coupled fluid variables on both sides of the boundary to the interfacial degrees of freedom. The original IB/MLS approach converts the cubic spline weights typically employed in MLS reconstruction into an IB kernel function that satisfies particular discrete moment conditions. This paper shows that the same approach can be used to construct one-sided kernel functions (kernel functions are referred to as generating functions in the MLS literature). We also examine the performance of the new approach for a family of kernel functions introduced by Peskin. It is demonstrated that the one-sided MLS construction tends to generate non-monotone interpolation kernels  with large over- and undershoots. We present two simple weight shifting strategies to construct generating functions that are positive and monotone, which enhances the stability of the resulting IB methodology. Benchmark cases are used to test the order of accuracy and verify the one-sided IB/MLS simulations in both two and three spatial dimensions. This new IB/MLS method is also used to simulate flow over the Ahmed car model, which highlights the applicability of this methodology for modeling complex engineering flows.

\end{abstract}


\end{frontmatter}


\section{Introduction} \label{sec_intro}

Immersed boundary (IB) methods~\cite{Peskin02,Mittal05,Griffith20} are a class of fictitious domain 
methods~\cite{Patankar2000} that can enable the efficient solution to complex moving domain problems. 
The IB methodology is now widely used in modeling large scale engineering~\cite{Mittal05,Borazjani13,Bhalla13,Gazzola2011,Nangia2019,Li16,khedkar2021inertial,bhalla2020simulating} and biological~\cite{Griffith09,Griffith17,Bhalla13FDO,Luo08,Namu2018,Sprikle17} models that have large \REVIEW{boundary} deformations or displacements within the computational domain. IB methods have also seen substantial use for applications that require only fixed geometries because they simplify grid generation~\cite{Thirumalaisamy2021a,Thirumalaisamy2021b,helgadottir2015imposing,chai2021imposing}. Compared to body-fitted and unstructured grid approaches, IB methods have lower memory footprint, are easier to parallelize, and allow for fast linear solvers.

In IB method, the structural displacement, velocity, and forces are described on a Lagrangian mesh, whereas the fluid velocity and pressure are described on a background Eulerian grid.  The Lagrangian mesh is allowed to cut the background Cartesian grid arbitrarily, which allows the IB methodology to be flexible and efficient. Interactions between the Lagrangian and Eulerian variables are typically mediated by integral transforms with a smooth delta function kernel. The regularized delta function kernels effectively smear the interface over a few Eulerian grid cells. The smearing process reduces the order of accuracy of the methodology, which can necessiate using high grid resolution near the interface, particularly at higher Reynolds numbers. It, however, permits a continuous solution of velocity and pressure fields across the interface. A limitation of the standard diffuse-interface IB method is that for immersed bodies, the continuous solution can lead to spurious flow inside the structure~\cite{Jiang19}. The internal flow may not influence the external flow at moderate Reynolds numbers, but it can interfere with the external boundary layer at high Reynolds number (see Sec.~\ref{sec_spherecase}).

Sharp-interface versions of the IB method have also been developed that reconstruct the velocity along the interface using procedures that are not naturally expressed as integral transforms~\cite{RMittal08,Uday01,Borazjani2008}. In some of these approaches, fluid boundary conditions are directly imposed along the immersed boundary, and the momentum and continuity equations are solved only on one side of the boundary. When the structure moves, fluid cells are ``covered" and ``uncovered", and this can lead to spurious force oscillations in the time history of the integrated drag and lift quantities. Although these approaches can improve the order of accuracy as compared to formulations that use regularized forcing along the boundary, they require rich geometric information about the interface and its location relative to the background grid. Sharp IB methods also require the grid nodes to be classified as “IB”, “fluid”, and “solid” nodes, which increases the coding complexity, especially for large scale distributed computing models. In particular, for moving body problems, the node classification has to be reevaluated at every time step.

This work uses a direct forcing immersed boundary method~\cite{Uhlmann05,Bhalla13,Vanella09} to model flows around complex bodies with prescribed kinematics. This is motivated by vehicular aerodynamics applications, for the which the “dirty” CAD geometry (which may include features such as sharp or thin edges, or artifacts such as holes) of the vehicle is relatively easier to handle via diffuse-interface IB methods than their sharp-interface counterparts~\cite{Jansson18,Bale20b}.

Various approaches have been employed in the literature to construct regularized delta functions for diffuse IB methods.
A common approach is to construct them by imposing certain discrete moment conditions. This approach was pioneered
by Peskin~\cite{Peskin02} and has been extended in subsequent work, including by Yang et al.~\cite{Yang09} and Stein et al.~\cite{Stein2016}. Another approach to constructing kernels functions is to employ the moving least squares (MLS) methodology~\cite{Wendland01}. The MLS technique can be used to construct interpolation kernels (which are referred to as generating functions in the MLS literature) for a particular Lagrangian point dynamically by solving a weighted least squares problem defined over a set of positions (interpolation points) where data are sampled for interpolation (to locations on the Lagrangian mesh). The generating functions produced by the MLS procedure are constructed to satisfy prescribed conditions, such as the polynomial reproduction constraints, which are equivalent to discrete moment conditions like those suggested by Peskin. One of the strengths of the MLS method is that it does not require the interpolation points to be arranged  in a particular way around a Lagrangian point. This makes the MLS technique quite general and even applicable for unstructured grid IB methods~\cite{Krishnan17,Saadat18}. As typically done in regular IB methods,  prior IB/MLS works have also used a full support of interpolation points around a Lagrangian marker. In contrast, we employ the MLS technique to construct \emph{one-sided} regularized delta functions for the direct forcing IB method in this work.

\begin{figure}
	\centering
	\subfigure[A closed interface]{
	   	\includegraphics[scale=0.28]{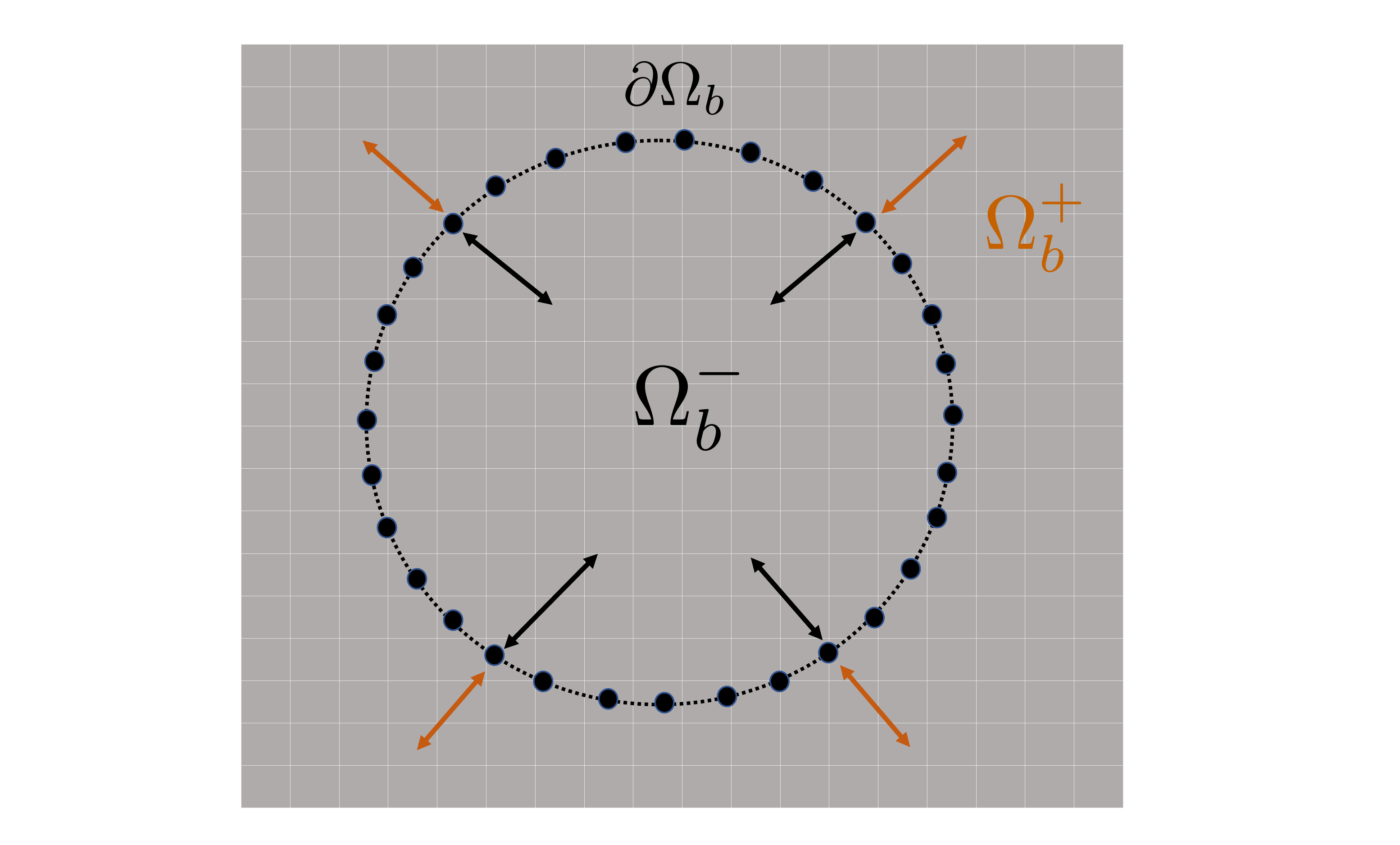}
	   	\label{subfig:closed}
	   }
	   \subfigure[An open interface]{
	   	\includegraphics[scale=0.32]{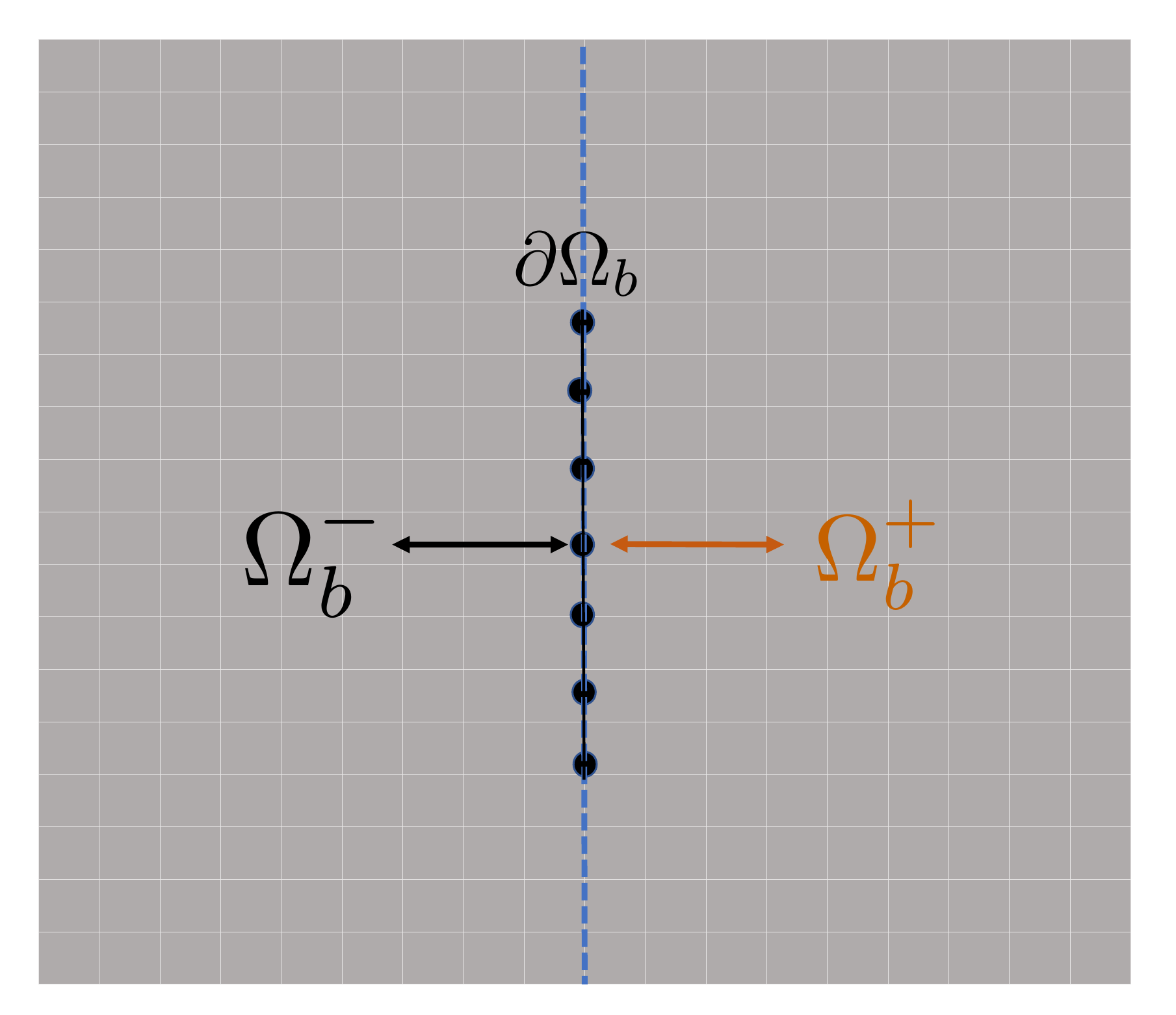}
		\label{subfig:open}
	   }
	\caption{Demarcation of Eulerian domain $\Omega = \Omega_b^{-} \cup \Omega_b^{+}$ based on the location of the interface $\partial \Omega_b$. Lagrangian marker points depicted as $[\bullet]$ on the interface interpolate velocity from and spread force to  $\Omega_b^{-}$ and $\Omega_b^{+}$ regions separately.~\subref{subfig:closed} A self-closing interface partitioning the domain into interior $\Omega_b^{-}$  and exterior $\Omega_b^{+}$ regions.~\subref{subfig:open} A hypothetical plane (shown with a dashed line) passing through an open interface and partitioning the domain into interior and exterior regions. The choice of the hypothetical plane for an open interface is arbitrary.}
	\label{fig:MLSSchematic}
\end{figure}

The MLS methodology has been extensively used in the mesh-free continuum mechanics literature~\cite{Liu95,Liu97,Mirzaei12,Sukumar07}. Its use with the direct forcing IB method was introduced by Vanella and Balaras~\cite{Vanella09}. Recently, Li et al.~\cite{Li15} and Tullio and Pascazio~\cite{Tullio16}  have also used the MLS technique of Vanella and Balaras to construct regularized IB kernels on structured Eulerian meshes. In prior studies~\cite{Vanella09,Li15,Tullio16}, the MLS technique was used to transform cubic spline weight functions into regularized kernel functions that satisfy zeroth- and first-order discrete moment conditions.  Because the fluid is described using a regular Cartesian grid, a Cartesian arrangement of interpolation points across the interface was used in defining the least squares problem. It is important to notice, however,  that if the weight function already satisfies the conditions imposed by the MLS construction, then the generating function produced by the MLS procedure is the same as the weight function. For example, Peskin's kernel functions and certain spline functions, already satisfy the reproducing properties, and hence are not modified by the MLS construction. We prove this property of MLS weight generation in Sec.~\ref{sec_mls_theory}. Consequently, there is no advantage to using a MLS construction with weighting functions that already satisfy the desired properties unless the construction aims to reduce the support of the weighting function. In contrast, if the MLS kernel is to be supported on a region that does not contain the full support of the basic kernel function, the MLS kernel will generally be different from the basic kernel. This approach allows us to interpolate velocity from and spread force to only one side of the interface (see Fig.~\ref{fig:MLSSchematic}). 

We remark that although the one-sided generating functions will satisfy the desired moment conditions, they will generally include large over- and undershoots, even if the original weights are non-negative and monotonically decreasing away from the Lagrangian point. For diffuse-interface IB methods, the positive and negative values in the generating functions can produce oscillatory feedback through force spreading operation that can destabilize the FSI simulation. By construction, the velocity interpolation operation is not affected by (possibly larger-in-magnitude) negative weights. The source of the aforementioned instability can be understood as follows. For the direct forcing IB method, the feedback force is defined in terms of a slip velocity at a Lagrangian point. If the weights used in the force spreading operator are of similar magnitude, but with opposite sign, the Eulerian representation of the Lagrangian feedback force will be spatially oscillatory. This causes the numerical error to grow over time. We remark that for exactly-constrained IB methods, wherein the Lagrangian forces are computed by solving an extended saddle-point system~\cite{Kallemov16,Usabiaga17}, the (larger magnitude of) positive and negative weights do not affect the stability of the system --- the IB force is dynamically adjusted to satisfy the imposed kinematic constraints. To overcome this issue, this paper introduces two shifting techniques to make weights non-negative and to reduce the overshoots for the one-sided kernels. A side-effect of the shifting procedure is that the first and other higher-order moment conditions are not satisfied in general, and as a result, the shifted one-sided kernels are limited to first-order accurate reconstructions.

One-sided MLS kernels have also been used by Le and Khoo~\cite{Le2017}  and Mohammadi et al.~\cite{Haji19} in the 
context of sharp-interface IB methods. These developments are relatively recent. In their approach, velocity reconstruction at ``IB forcing nodes" is done via an MLS method instead of grid based interpolations. For a given IB forcing node (which is defined to be near the solid interface and towards the fluid side), nearby Eulerian fluid nodes and interfacial Lagrangian nodes are selected as a set of interpolation points to carry out the velocity reconstruction process. Internal solid nodes and neighboring IB forcing nodes are excluded from this interpolation set~\footnote{These nodes are also excluded from grid based interpolation schemes that are typically used in sharp IB methods.}. Le and Khoo~\cite{Le2017} and Mohammadi et al.~\cite{Haji19} found MLS based reconstructions to be more robust than interpolating the velocity along a specific direction, such as grid-aligned directions or local normal directions and as commonly used with sharp-interface IB schemes. Their MLS based sharp-interface IB method also retains second-order accuracy. We remark that there is a fundamental difference in the MLS interpolation procedure for sharp-interface IB methods and our current diffuse-interface IB method. For sharp-interface IB methods, the interfacial Lagrangian nodes participate in the interpolation process because they are not considered to be forcing points, and IB nodes near the solid interface and towards the fluid side act as forcing points. For diffuse-interface IB methods, the interfacial Lagrangian nodes are the forcing points, and therefore, only neighboring Eulerian fluid nodes are taken in the support of the interpolating kernel. The reason for exclusion stems from the fact that the range of spreading operator, which is typically constructed to be the adjoint of interpolation operator, is restricted to the fluid nodes. Conversely, the domain of the interpolation operator is also restricted to the same set of Eulerian fluid nodes.

\REVIEW{Specialized forms of one-sided IB kernels dealing with interactions between immersed structures and the boundaries of the computational domain have also been proposed in the literature. For example, Griffith and co-workers algebraically derived one-sided regularized delta functions that continuously transition to Peskin’s four-point kernel as the structure position moves away from the boundary of the computational domain~\cite{Griffith2009simulating}. As in the simplest version of the MLS weighting functions constructed herein, this one-sided kernel included both positive and negative weights. An alternative approach was proposed by Yeo and Maxey to deal with boundaries in the computational domain~\cite{Yeo2010dynamics}; see also Delong et al.~\cite{Delong2014brownian}. However, these constructions do not deal with one-sided constructions along general (curved) internal boundaries, which is the focus of this paper. Moreover, the current MLS approach can also naturally handle immersed structures near the computational domain boundaries.}


\section{Equations of motion}

\subsection{Immersed boundary method} \label{sec_cont_eqs}

We begin by describing the equations of motion of a rigid body moving in an incompressible viscous fluid 
of constant density $\rho$ and viscosity $\mu$ in a region $\Omega \subset \mathbb{R}^d$, $d = $ 2 or 3.  
In the immersed boundary formulation, a Lagrangian description is employed for the structural location 
and forces, and an Eulerian description is employed for the fluid velocity and pressure.
We denote the position of the boundary of the Lagrangian structure $\partial \Omega_{b}(t) \subset \Omega$ at time 
$t$ by $\X (\s,t)$, in which $\s = (s_1, \ldots s_{d-1}) \in U$ denote a fixed material coordinate system 
attached to the boundary, and $U \subset \mathbb{R}^{d-1}$ is the Lagrangian curvilinear coordinate
domain. The fixed physical domain is fully described by the  Cartesian coordinates $\x = (x_1, \ldots, x_d) \in \Omega$.
The open or closed boundary of the immersed structure partitions the Cartesian domain into two regions 
$\Omega = \Omega_{b}^{+}(t) \cup \Omega_{b}^{-}(t)$ (see Fig.~\ref{fig:MLSSchematic}).

The combined equations of motion for 
fluid-structure system are~\cite{Peskin02,Bhalla13}

\begin{align}
\rho\left(\D{\u}{t}(\x,t) + \u(\x,t) \cdot \grad \u(\x,t) \right) &= -\grad p(\x,t) + \mu \lap \u(\x,t) + \f(\x,t), \label{eqn_momentum}\\
  \div \u(\x,t) &= 0, \label{eqn_continuity} \\
\f(\x,t)  &= \int_{U} \F(\s,t) \, \delta(\x - \X(\s,t)) \, \Ds, \label{eqn_F_f} \\
  \U(\s,t) &= \int_{\Omega} \u(\x,t) \, \delta(\x - \X(\s,t)) \, \Dx, \label{eqn_u_interpolation} \\
   \D{\X}{t} (\s,t) &= \U(\s,t). \label{eqn_body_motion} 
\end{align}
Eqs.~\eqref{eqn_momentum} and \eqref{eqn_continuity} are the incompressible Navier-Stokes equations written 
in Eulerian form, in which $\u(\x,t)$  is the velocity and $p(\x,t)$ is the pressure. Eq.~\eqref{eqn_F_f}  converts 
the Lagrangian force density $\F(\s,t)$ into an equivalent Eulerian density $\f(\x,t)$.
Conversely, Eq.~\eqref{eqn_u_interpolation} determines the physical velocity of each Lagrangian material point 
from the Eulerian velocity field, so that the immersed structure moves according to the local value of the velocity
field  $\u(\x,t)$  (Eq.~\eqref{eqn_body_motion}). The Lagrangian and Eulerian quantities of Eqs.~\eqref{eqn_F_f} and~\eqref{eqn_u_interpolation} are mediated by integral equations using Dirac delta function, which is taken to be $d$-dimensional tensor 
product of one-dimensional delta functions: $\delta(\x) = \Pi_{i=1}^{d}\delta(x_i)$. In this formulation, the Lagrangian force density $\F(\s,t)$ is a Lagrange multiplier that ensures that the body moves according to a prescribed velocity.

\subsection{Discrete equations of motion} \label{sec_discrete_eqs}

\subsubsection{Lagrangian discretization}
In conventional IB methods, Eqs.~\eqref{eqn_F_f} and~\eqref{eqn_u_interpolation}  are approximated by replacing the singular delta function kernel by a regularized delta function $\delta_h(\x)$ that is supported on both sides of the interface $\partial\Omega_b(t)$. Using a short-hand notation, we denote the force spreading operation by $\f = \cSIB[\X] \, \F$, in which $\cSIB[\X]$ is the \emph{force-spreading operator}  associated with the boundary configuration and regularized delta function. The velocity interpolation operation as defined in Eq.~\eqref{eqn_u_interpolation} is expressed as $ \D{\X}{t}= \U = \cJIB[\X]\, \u$, in which $\cJIB[\X]$ is the \emph{velocity-interpolation operator}. It can be shown that if $\cSIB$ and $\cJIB$ are adjoint operators, i.e. $\cSIB = \cJIB^{*}$,
then Lagrangian-Eulerian coupling conserves energy~\cite{Peskin02}. Later in Sec.~\ref{sec_mls_ops} we will obtain spreading and interpolation operators using generating functions produced by the MLS technique, which we denote by $\cSMLS$ and $\cJMLS$, respectively. Generic force-spreading and velocity-interpolation operators are denoted $\cS$ and $\cJ$, respectively.


We impose velocity boundary conditions only along the fluid-structure interface, and the volume enclosed by the closed interface is not constrained.  Discretely, Lagrangian markers with 
curvilinear mesh spacing $(\Delta s_1, \Delta s_2)$ represent the codimension-1 interface $\partial \Omega_b$; see Fig.~\ref{fig:MLSSchematic}. 
The Lagrangian markers are indexed by the tuple $(l,m)$. Based on our prior experience with 
direct forcing IB method, we typically take the Lagrangian marker spacing 
approximately equal to the Eulerian grid spacing on the finest grid level, i.e. $\Delta s_1  \approx \Delta s_2 \approx h$, 
in which $h$ is the uniform grid cell size on the finest level of the locally refined grid.  A discrete approximation to any general quantity defined on marker points
is described by $\Phi^n_{l,m} \approx \Phi(\s_{l,m}, t^n) = \Phi(l\ds_1, m\ds_2,t^n)$ at time $t^n$.
More specifically, the position, velocity, and force of a marker point are denoted as 
$\X_{l,m}$, $\U_{l,m}$, and $\F_{l,m}$, respectively. Using the regularized IB kernel $\delta_h$, 
the discrete velocity interpolation of the collocated grid fluid velocity $\u \equiv (u,v,w)$ onto a specific configuration 
of Lagrangian markers (i.e. $\U \equiv (U,V,W) = \cJIB[\X] \, \u)$ reads 
\begin{align}
U_{l,m} & = \sum_{\x_{i,j,k} \in \Omega} u_{i,j,k} \delta_h\left(\x_{i,j,k} - \X_{l,m}\right)  \dx\dy\dz, \label{eqn_U_interp}\\  
V_{l,m} & = \sum_{\x_{i,j,k} \in \Omega} v_{i,j,k} \delta_h\left(\x_{i,j,k} - \X_{l,m}\right)  \dx\dy\dz, \\
W_{l,m} & = \sum_{\x_{i,j,k} \in \Omega} w_{i,j,k} \delta_h\left(\x_{i,j,k} - \X_{l,m}\right)  \dx\dy\dz, 
\end{align}
in which the tuple $(i,j,k)$ is the Cartesian grid cell index and $\x_{i,j,k}$ is the spatial location of the cell centroid.
Conversely, the discrete spreading of Lagrangian force density $\F \equiv (F_1,F_2,F_3)$  
onto cell centers of the collocated grid
(i.e. $\f \equiv (f_1, f_2, f_3) = \cSIB[\X] \, \F$) reads
\begin{align}
(f_{1})_{i,j,k} & = \sum_{\X_{l,m} \in \Omega_b} (F_{1})_{l,m} \delta_h\left(\x_{i,j,k} - \X_{l,m}\right) \ds_1 \ds_2, \\
(f_{2})_{i,j,k} & = \sum_{\X_{l,m} \in \Omega_b} (F_{2})_{l,m} \delta_h\left(\x_{i,j,k} - \X_{l,m}\right) \ds_1 \ds_2, \\
(f_{3})_{i,j,k} & = \sum_{\X_{l,m} \in \Omega_b} (F_{3})_{l,m} \delta_h\left(\x_{i,j,k} - \X_{l,m}\right) \ds_1 \ds_2. \label{eqn_f3_spread}
\end{align}

The regularized IB kernel $\delta_h$ interpolates from and spreads to both sides of the interface $\partial \Omega_b$. 
Note that the arguments of an IB kernel $\delta_h$ depend only on the difference between Eulerian and 
Lagrangian locations, i.e. $\delta_h = \delta_h(\x - \X)$, and not the locations themselves; 
see Eqs.~\eqref{eqn_U_interp}-\eqref{eqn_f3_spread}. This property implies that the weights of an IB kernel are 
independent of the Lagrangian marker identity and the same kernel is used for each marker.  In Sec.~\ref{sec_mls_theory}, 
we describe a procedure to obtain the one-sided version of IB kernel $\psi_h(\x,\X)$ using the moving least squares method whose weights depend upon the Lagrangian marker position explicitly. In particular, a different kernel function is generated for each marker.

\subsubsection{Eulerian discretization and time-stepping scheme}
In this work, a collocated grid discretization for the momentum and continuity equations is used,
in which the Eulerian velocity, pressure, and force variables are defined at the centers of Cartesian grid cells of grid spacing $h$. 
Second-order finite differences are used to approximate the Eulerian equations  
on locally refined grids~\cite{Jansson18,Bhalla13}.  The spatially discretized cell-centered operators
are denoted with a `cc' subscript. To avoid the velocity and pressure decoupling on 
collocated grid~\cite{Patankar2018}, auxiliary face-centered variables and operators are introduced, 
which are distinguished using a `fc' subscript.  

A version of the second-order accurate pressure projection algorithm of Brown et al.~\cite{Brown2001} is used 
to solve the incompressible Navier-Stokes system, and the direct forcing approach of Bhalla et al.~\cite{Bhalla13}  is used
to approximately impose the rigidity constraint of the immersed body. A fixed-point iteration time stepping scheme using $\ncycles = 2$ cycles 
per time step is used to evolve quantities from time level $t^n$ to time level $t^{n+1} = t^n + \Delta t$. A superscript ``$k$" is used to denote the
cycle number of the fixed-point iteration. At the beginning of each time step, the solutions
from the previous time step are used to initialize cycle $k = 0$:
$\u^{n+1,0} = \u^{n}$ and $p^{n+\half,0} = p^{n-\half}$.  At the initial time $n = 0$, the Eulerian velocity is prescribed via an initial condition, and the Eulerian pressure is taken to be zero.

With $\V{N} =   (\u \cdot \grad_{h} \u)$ denoting the non-linear convective term, the equations of motion are
\begin{align}
\rho \left(\frac{\u^{n+ 1, k+1} - \u^n}{\dt} +  \V{N}^{n+\half, k}\right) &=  - \grad_{h} p^{n + \half,k+1} + \mu \V{\grad}^2_{h} \left( \frac{\u^{n+ 1,k} + \u^{n}}{2}\right) + \f^{n+\half, k+1}, \label{eqn_lm_momentum} \\
\V{\grad}_h \cdot \u^{n+ 1, k+1} & = 0, 
\end{align}
and are integrated from time step $n$ to $n+1$ in an operator-splitting manner as follows: 

%
%
%

\begin{enumerate}

\item An intermediate velocity $\widetilde{\u}_{*}^{n+ 1,k+1}$ is obtained by integrating the advection-diffusion 
momentum equation, which ignores the pressure and the constraint forces

\begin{equation}
\rho \left(\frac{\widetilde{\u}_{*}^{n+ 1,k+1} - \u^n}{\dt} +  \V{N}^{n+\half,k}\right) =  \mu \V{\grad}^2_h \left( \frac{\widetilde{\u}_{*}^{n+ 1,k+1} + \u^{n}}{2}\right).\label{eqn_lm_momentum}
\end{equation}
A geometric multigrid solver using a Gauss-Seidel smoother is employed to solve for  $\widetilde{\u}_{*}^{n+ 1,k+1}$
on a locally refined grid. For the first ($k = 0$) cycle, the nonlinear convective term is discretized using the explicit Adams-Bashforth scheme, so that $\V{N}^{n+\half,0} =  \frac{3}{2} (\u \cdot \grad_h \u)^{n} -  \frac{1}{2} (\u \cdot \grad_h \u)^{n-1}$, while for the remaining cycles, a midpoint approximation is used, so that $\V{N}^{n+\half,k} =   (\u \cdot \grad_h \u)^{n+\half,k}$, in which  $\u^{n+\half,k} = \frac{1}{2} (\u^{n+ 1,k} + \u^{n})$.

\item Impose the rigidity constraint of the immersed body by computing the constraint force $\F$ 

\begin{align}
	\Delta \U^{n+\half, k+1} &= \U_b^{n+\half} -  \cJ[\X] \left(\frac{\widetilde{\u}_{*}^{n+1,k+1} + \u^n}{2} \right),  \label{eqn_vel_incr} \\
	\widetilde{\u}^{n+ 1,k+1} &=  \widetilde{\u}_{*}^{n+1, k+1} + \cS[\X] \, \Delta \U^{n+\half,k+1},  \label{eqn_vel_correction} \\
	\F^{n+\half,k+1} &= \frac{\rho}{\Delta t} \Delta \U^{n+\half,k+1}, \label{eqn_lag_force}      \\
	\f^{n+\half,k+1} &=  \cS[\X] \,  \F^{n+\half,k+1}, \label{eqn_eul_lag_force}
\end{align}
in which $\U_b$ is the desired rigid body velocity of the structure, and $\Delta \U$ is the slip velocity  
computed from the interpolated (unconstrained) fluid velocity on the Lagrangian markers. In our implementation, we directly 
update the velocity field $\widetilde{\u}^{n+1,k+1}$ using Eq.~\ref{eqn_vel_correction} on the Eulerian grid, without 
explicitly calculating $\f^{n+\half,k+1}$ as given in Eq.~\ref{eqn_eul_lag_force}. This equivalence can be readily verified by
considering the operator-split between momentum and rigidity constraint equations~\cite{Bhalla13}. However, Eq.~\ref{eqn_eul_lag_force} is
useful for computing the net hydrodynamic forces on the immersed structure~\cite{Nangia17}.


At this stage, the fluid velocity $\widetilde{\u}^{n+ 1,k+1}$ accounts for the rigidity of the structure, but not the incompressibility 
of the system. This is corrected in the next step. 

\item Impose the incompressibility constraint by solving the pressure Poisson equation for the auxiliary variable $\phi$ and 
estimate the fluid pressure

\begin{align}
        \widetilde{\u}^{n+ 1,k+1}_{\textrm{fc}} &= \mathbb{I}(\widetilde{\u}^{n+ 1,k+1}), \label{eqn_adv_vel} \\	
	-\grad^2  \phi^{n+1,k+1} &= -\frac{\rho}{\Delta t} \grad \cdot \widetilde{\u}^{n+ 1,k+1}_{\textrm{fc}}, \label{eqn_projection} \\
	\u^{n+ 1,k+1} & =  \widetilde{\u}^{n+ 1,k+1}    -\frac{\Delta t}{\rho} \grad \phi^{n+1,k+1}_{\text{cc}}, \label{eqn_vel_projection} \\
	\u^{n+ 1,k+1}_{\textrm{fc}} & =  \widetilde{\u}^{n+ 1,k+1}_{\textrm{fc}}    -\frac{\Delta t}{\rho} \grad \phi^{n+1,k+1}_{\textrm{fc}}, \label{eqn_vel_projection} \\	
	p^{n + \half,k+1} & =   \phi^{n+1,k+1} - \frac{\Delta t}{\rho} \frac{\mu}{2} \grad^2 \phi^{n+1,k+1}_{\textrm{cc}},  \label{eqn_pressure_incr} 
\end{align}
in which $\mathbb{I}$ is the interpolation operator that determines the face-centered velocity field $\u^{n+ 1,k+1}_{\textrm{fc}}$ by averaging the adjacent cell-centered velocity components.

\end{enumerate}

We use adaptive mesh refinement (AMR) framework for some of the cases presented in 
Sec.~\ref{sec_results}. A grid with $\ell_{\text{max}}$ refinement levels has grid spacings $\dx_0$, $\dy_0$, and $\dz_0$ on the coarsest grid level and grid spacings $\dx = \dx_0/\nref^{\ell-1}$, 
$\dy = \dy_0/\nref^{\ell-1}$, and $\dz = \dz_0/\nref^{\ell-1}$ on a grid level $\ell$,
in which, $\nref$ is the integer refinement ratio. The Lagrangian mesh is embedded on the 
finest grid level to adequately resolve the thin boundary layers. For all of the cases considered 
in this work a constant time step size $\dt = \min(\dt^{\ell})$ is chosen, in which the time step size $\dt^{\ell}$ on grid level $\ell$ satisfies 
the convective CFL condition $\dt^{\ell} \le C \min \left(\frac{\dx}{\|u_x\|_{\infty}}, \frac{\dy}{\|u_y\|_{\infty}}, \frac{\dz}{\|u_z\|_{\infty}}\right)^{\ell}$. 
In this work, the convective CFL number is set to $C = 0.1$ unless otherwise stated.


\section{Moving least squares method} \label{sec_mls_theory}

\subsection{Review}
\REVIEW{The Backus-Gilbert formulation of the moving least squares method seeks} the quasi-interpolant 
\begin{equation}
\cP g(\X) = \sum_i^N g(\x_i) \psi_i(\X),      
\end{equation}
in which $\g = [g(\x_1), \ldots, g(\x_N)]^T$ are given data at $N$ \emph{interpolation} points, $\X$ is the \emph{evaluation} point, and $\cP$ is the interpolation operator. 
The moving least squares method computes \emph{generating functions} $\psi_i(\X) = \{\psi(\x_i,\X)\}$ subject 
to the polynomial reproduction constraints

\begin{equation}
\InterpSumi p(\x_i) \psi_i(\X) = p(\X), \quad \text{for all } p \in \Pi_{d}^s,  \label{eqn_poly_reproduce}
\end{equation}
in which $\Pi_{d}^s$ is the space of s-variate polynomials (centered around any arbitrary 
point) of total degree at most $d$. The polynomial reproduction constraints 
correspond to \emph{discrete moment conditions} for the function $\psi_i(\X)$. Eq.~\eqref{eqn_poly_reproduce} can be written in a matrix form as
\begin{equation}
\cA \V{\Psi} (\X) = \V P(\X),  \label{eqn_linsys}
\end{equation}
in which the entries of the polynomial matrix $\cA \in \mathbb{R}^{m \times N}$ are the values of the basis functions at the data point 
locations, $\cA_{ij} = p_i(\x_j), i = 1, \ldots,m, j = 1, \ldots, N$, and the right-hand side vector $\V P = [p_1, \ldots, p_m]^T$
contains the values of the polynomials at the evaluation point $\X$. The unknown generating function 
vector $\V{\Psi} = [\psi_1, \ldots, \psi_N]^T$ is obtained by solving a \emph{least squares problem}.
Since the set of generating functions changes by considering a different evaluation point $\V Y \ne \X$, $\{\psi_i(\V Y) \} \ne \{\psi_i(\V X)\}$,
this least squares procedure is called the \emph{moving least squares problem}.

Generally $N \gg m$, implying that Eq.~\ref{eqn_linsys}  is an underdetermined system which can be solved in a weighted least 
squares sense with the help of Lagrange multipliers $\V{\lambda} (\X)$ to enforce the reproducing conditions. \REVIEW{Specifically, the Backus-Gilbert MLS problem is solved by posing it as a constrained quadratic minimization problem, which reads as

\begin{equation}
\label{eqn_probBG}
\text{Backus-Gilbert MLS :}  \begin{cases}
  \min\limits_{\V{\Psi}(\X) \in  \mathbb{R}^{N} }  & J  \; = \; \half\, \V{\Psi}^\intercal(\X) \cW^{-1}(\X) \V{\Psi}(\X) \\
  \text{subject to:} & \cA \V{\Psi}(\X) =  \V{P}(\X).
  \end{cases}
\end{equation} 
Here, $\cW(\X) =  \text{diag}\left( \V W \right)$ is the diagonal matrix containing weights of the interpolation points with respect to the evaluation point, and 
$\V W = \left(W(\x_1,\X), \dots, W(\x_N,\X) \right)$ is the main diagonal of $\cW(\X)$. The weights are defined by a nonnegative weight function $W(\x_i,\X)$ that decreases in magnitude with distance away from the evaluation point $\X$. Eq.~\eqref{eqn_probBG} is solved by minimizing the Lagrangian $\mathfrak{L}$ of the problem, defined using the Lagrange multiplier  $\V{\lambda} (\X) \in \mathbb{R}^m$ as
\begin{equation}
\label{eqn_LB}
\mathfrak{L} = \half\, \V{\Psi}^\intercal(\X) \cW^{-1}(\X) \V{\Psi}(\X) - \V{\lambda}^\intercal(\X)[\cA \V{\Psi}(\X) -  \V{P}(\X)]. 
\end{equation} 
The Lagrangian is minimized by using the essential conditions of extrema,   $\D{\mathfrak{L}}{\V{\Psi}} = \V{0}$ and $\D{\mathfrak{L}}{\V{\lambda}} = \V{0}$, which yields 
\begin{align}
\label{eqn_dLB_1}
\cW^{-1}  \V{\Psi} - \cA^\intercal \V{\lambda} & = \V{0}, \\
\label{eqn_dLB_2}
\cA \V{\Psi} - \V{P} &= \V{0}.
\end{align} 
Solving the above two equations, we obtain 
\begin{align}
\label{eqn_BG_lambda}
\V{\lambda} & =   \left( \cA \cW \cA^\intercal \right)^{-1} \V{P} = \cG^{-1} \V{P}, \\
\label{eqn_BG_psi}
 \V{\Psi}  &= \cW \cA^\intercal \V{\lambda} = \cW \cA^\intercal \left( \cA \cW \cA^\intercal \right)^{-1} \V{P} = \cW \cA^\intercal \cG^{-1} \V{P}.
\end{align} 

In Eqs.~\eqref{eqn_BG_lambda} and~\eqref{eqn_BG_psi},  $\cG(\X)  = \cA \cW \cA^\intercal \in \mathbb{R}^{m \times m}$ is the symmetric positive-definite Gram matrix, whose entries are the weighted $L^2$ inner products of the polynomials
\begin{equation}
\cG_{jk}(\X) = \left< p_j, p_k \right>_{W(\X)} = \InterpSumi p_j(\x_i) p_k(\x_i) W(\x_i, \X), \quad j,k = 1,\ldots,m. \label{eqn_gram_mat}
\end{equation}  
It is also instructive to write the generating functions (obtained in Eq.~\eqref{eqn_BG_psi}) in component form as
\begin{equation}
 \psi(\x_i, \X) = \psi_{i}(\X)  = W(\x_i,\X) \InterpSumj \lambda_j(\X) p_j(\x_i), \quad i = 1, \ldots, N. \label{eqn_psi}
\end{equation}
In matrix notation, Eq.~\eqref{eqn_psi} can be written as
\begin{equation}
 \V{\Psi} (\X)  =  \cW (\X) \odot \vcL(\X), \label{eqn_matrix_psi}
\end{equation}
in which $\vcL(\X) = \cA^\intercal \V{\lambda}(\X) $ and $\odot$ indicates the Hadamard (component-wise) product of two matrices.

In the IB literature, the standard formulation of the MLS problem is more popular~\cite{Vanella09,Tullio16,Krishnan17,Saadat18,Le2017,Haji19} compared to the Backus-Gilbert theory~\cite{Backus1968resolving}. In the standard formulation of MLS, the quasi-interpolant  $\cP g$ to $g$ is expressed as
\begin{equation}
\label{eqn_stdmls}
\cP g(\X) = \sum_{j= 1}^m c_j(\X) p_j(\X)  = \V{P}^\intercal (\X)\V{C}(\X),      
\end{equation}
and the unknown coefficient vector $\V{C}(\X)$ appearing in Eq.~\eqref{eqn_stdmls} above is found by minimizing the weighted $L^2$-norm of the error function $J(\V{C})$: 
\begin{equation}
 \label{eqn_stderror}
J(\V{C}(\X)) = \sum_{i= 1}^N  W(\x_i,\X) [\V{P}^\intercal (\x_i)\V{C}(\X) - g(\x_i)]^2,      
\end{equation}   
with respect to $\V{C}$, i.e by setting $\D{J}{\V{C}} = 0$. However,  both formulations ultimately produce the same generating function  $\psi(\x_i, \X)$. There are two main advantages of the Backus-Gilbert MLS formulation: (i) the polynomial reproducing conditions are directly included in the problem formulation; and (ii) the unknown vector $\V{\Psi} $ is independent of the data values $\V{g}$, unlike the coefficient vector \V{C}, which implicitly depends upon the sampled data values, i.e. $\V{C} = \V{C}(\X;\V{g})$.  
}  

As an example, consider the univariate polynomials in two spatial dimensions, i.e. $d = 2$ and $s = 1$, so that  
 $m = 3$. The univariate basis functions defined relative to an evaluation point are   
$p(\x) = (p_1(\x - \X), p_2(\x - \X), p_3(\x - \X)) = (1, x- X, y - Y)$, with $\X = (X, Y)$ denoting
the coordinates of the evaluation point $\X$.  In the context of the immersed boundary method, $p_1$ polynomial (constant) reproducing condition implies that Lagrangian and Eulerian forms of force and power are equivalent. Furthermore, $p_2$ and $p_3$ polynomial (linear) reproducing conditions make Lagrangian and Eulerian representations of torque equivalent~\cite{Peskin02}. Notice that defining the polynomial basis function relative to the evaluation point $\X$ makes the right-hand of Eq.~\eqref{eqn_linsys} a constant vector, i.e.,  
\begin{equation}
\V P (\X) = p(\x)|_{\x = \X} =  (1, 0, 0)^T,
\end{equation}
independent of $\X$.

\begin{mylem}
\label{lem:trivial_psi}
If the weight functions already satisfy the reproducing conditions then $\vcL$ is a matrix of ones.  
\begin{proof}
In this case the weight functions are the generating functions $\V{\Psi} (\X)   = \cW(\X)$. Therefore, 
from Eq.~\ref{eqn_matrix_psi} $\vcL_{ij} = 1$ follows.
\end{proof}
\end{mylem}

\begin{figure}
	\centering
	\subfigure[]{
		\includegraphics[scale=0.5]{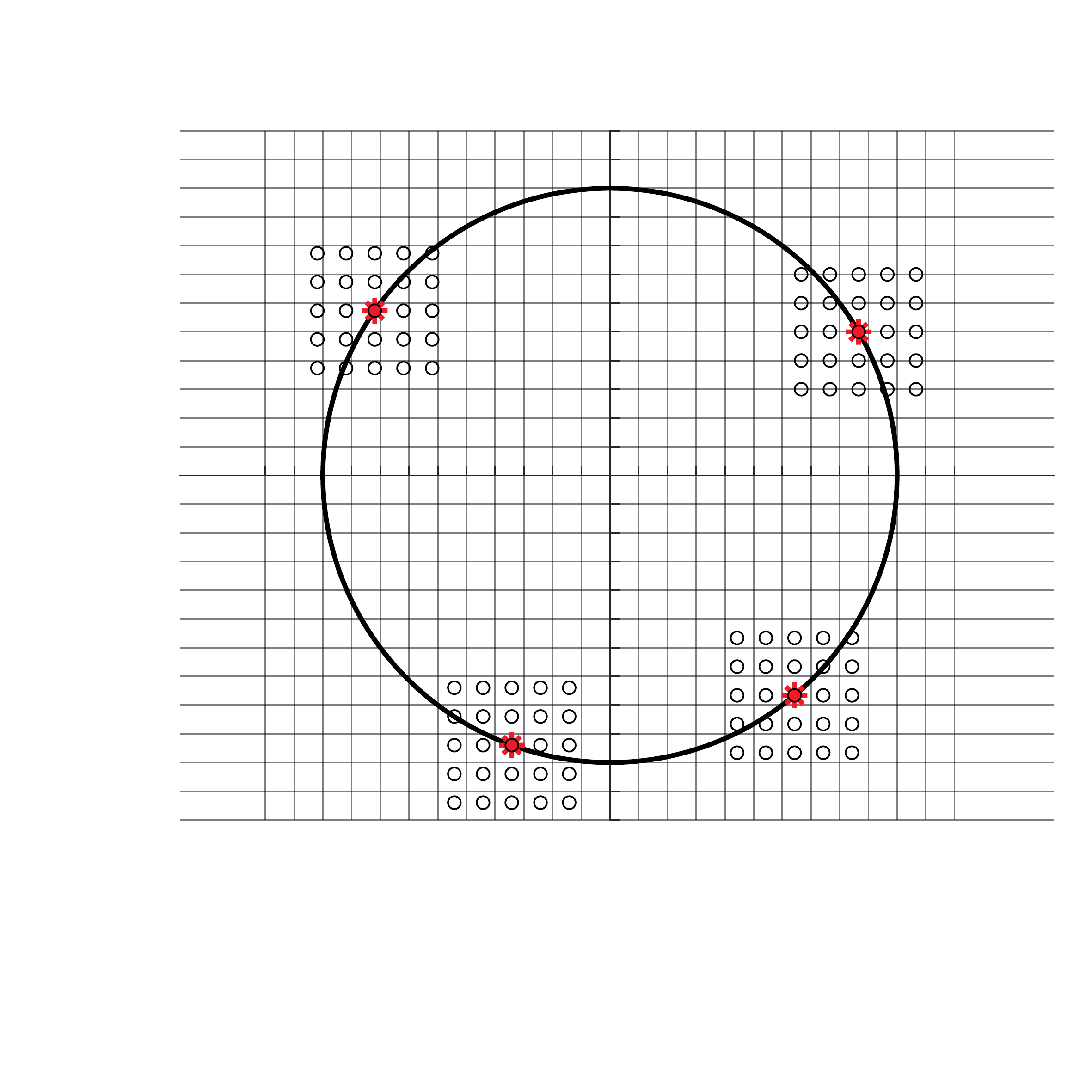} 
		\label{fig:error_psi_w_setup}
	}
	\subfigure[]{
		\includegraphics[scale=0.3]{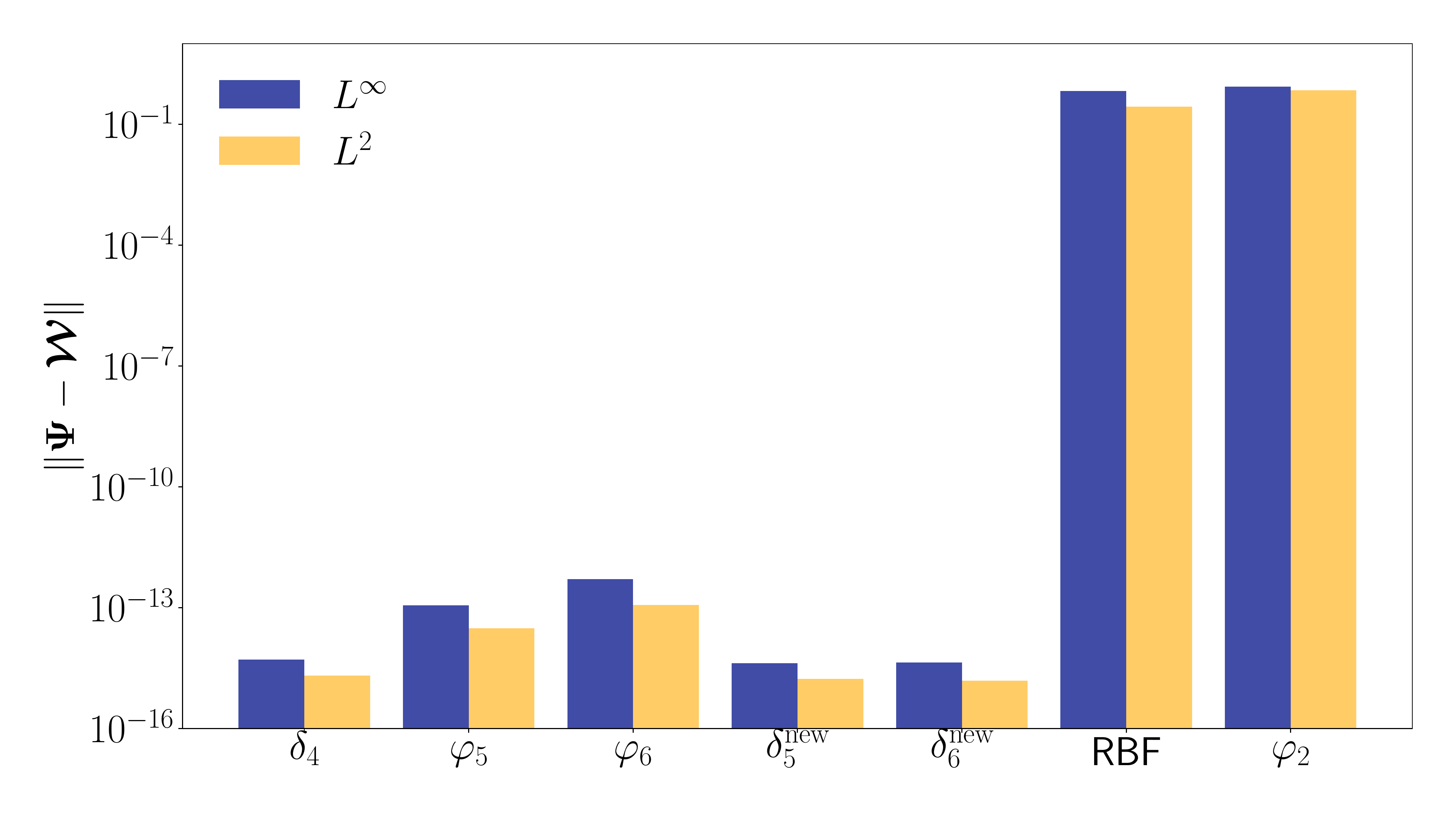} 
		\label{fig:error_psi_w}
	}	
	\caption{\subref{fig:error_psi_w_setup} Representative Lagrangian markers (${\color{red}*}$) on a circular interface embedded in a Cartesian grid, and full support of the kernel ($\circ$) around marker points.~\subref{fig:error_psi_w} Error norms of $|| \V{\Psi} - \cW ||$ for various kernels.}
	\label{fig:diff_psi_w}
\end{figure}

Next, we consider some commonly used weight functions $\cW$ in the IB literature to generate MLS weights $\V{\Psi}$ for four representative points on a circular interface (see Fig.~\ref{fig:error_psi_w_setup}). The chosen weight functions are: (i) Peskin's four-point delta kernel $\delta_4$; (ii) new five-point $\delta_5^{\text{new}}$ and (iii) six-point $\delta_6^{\text{new}}$ delta kernels; (iv) five-point  $\varphi_5$ and (v)  six-point  $\varphi_6$ spline functions; (vi)  two-point cubic spline function $\varphi_2$ as used in Vanella and Balaras~\cite{Vanella09}; and (vii) four-point radial basis function (RBF).  The one-dimensional functional form of these kernels are provided in Appendix~\ref{sec_kernel_appendix}. Out of these kernels, the first five kernels satisfy the zeroth- and first-order moment conditions, when they have a full support around the Lagrangian marker --- as in this case. Fig.~\ref{fig:error_psi_w} plots $L^2$ and $L^\infty$ norm of difference in $\cW$ and $\V{\Psi}$ for the four Lagrangian points. As observed in the figure, the MLS procedure does not transform the weights of those kernels that already satisfy the moment conditions it tries to impose. The numerical results shown in Fig.~\ref{fig:error_psi_w} are consistent with Lemma~\ref{lem:trivial_psi}.  
 
\subsection{One-sided IB kernels} \label{sec_mls_ops}

To obtain one-sided IB kernel $\psi_h$ using the MLS procedure, we use the Heaviside function $H(\x)$
\begin{align}
H_{i,j,k} &= 
\begin{cases} 
       0,  & \Omega_b^{-},\\
       1  &  \Omega_b^{+} \cup \partial \Omega_b,  \label{eq_heaviside_body}
\end{cases}
\end{align}
to select the appropriate side of the interface $\partial \Omega_b$ in order to define the domain of influence of the evaluation/Lagrangian point $\X$. 
This is achieved by multiplying the standard/unrestricted weights by the Heaviside function to obtain restricted weights $W_{\text{MLS}}$
\begin{equation}
W_{\text{MLS}}(\x_i,\X) = W(\x_i,\X)H(\x_i) \quad \text{s.t} \quad  W_{\text{MLS}}(\x_i,\X) = 0 \quad  \forall \quad \x_i \in \Omega_b^{-}.  \label{eqn_wmls}
\end{equation}
It can be easily verified that by using $W_{\text{MLS}}$ in Eq.~\ref{eqn_psi}, $\psi_h(\x_i,\X) = 0\; \forall\;  \x_i \in \Omega_b^{-}$.  
Therefore, by using a  simple weight manipulation process we can easily control the domain of influence or the interaction region of point $\X$.
    
Many choices for the weight function have been introduced in the MLS literature, including  the Gaussian function, B-splines, Wendland function, inverse distance square function, and other. In this work we use four-point, and the new/modified five-point and the six-point Peskin's delta functions; the new kernels remove the negative tails of the original Peskin kernels by imposing a weaker second moment condition. Because Peskin's delta functions are specifically constructed to satisfy the zeroth- and first-order moment conditions, they automatically qualify as the generating function when the domain of influence of the evaluation point $\X$ is 
not restricted. Furthermore, Peskin's delta functions show good grid translation invariance property in spite of a compact support~\cite{Bao2016}. In Sec.~\ref{sec_results},  we also compare the performance of Peskin's delta function weights to B-splines weights.      

Finally, using the one-sided IB kernel $\psi_h$,  the discrete velocity interpolation of the fluid velocity onto Lagrangian markers 
(i.e., $\U = \cJMLS[\X] \, \u)$ reads 
\begin{align}
U_{l,m} & = \sum_{\x_{i,j,k} \in \Omega} u_{i,j,k} \psi_h\left(\x_{i,j,k} , \X_{l,m}\right)  \dx\dy\dz, \label{eqn_U_interp_one_sided}\\  
V_{l,m} & = \sum_{\x_{i,j,k} \in \Omega} v_{i,j,k} \psi_h\left(\x_{i,j,k} , \X_{l,m}\right)  \dx\dy\dz, \\
W_{l,m} & = \sum_{\x_{i,j,k} \in \Omega} w_{i,j,k} \psi_h\left(\x_{i,j,k} , \X_{l,m}\right)  \dx\dy\dz, 
\end{align}
Conversely, the discrete spreading of Lagrangian force density $\F$  
onto cell centers using one-sided IB kernel
(i.e. $\f = \cSMLS[\X] \, \F$) reads
\begin{align}
(f_{1})_{i,j,k} & = \sum_{\X_{l,m} \in \Omega_b} (F_{1})_{l,m} \psi_h\left(\x_{i,j,k} , \X_{l,m}\right) \ds_1 \ds_2, \\
(f_{2})_{i,j,k} & = \sum_{\X_{l,m} \in \Omega_b} (F_{2})_{l,m} \psi_h\left(\x_{i,j,k} , \X_{l,m}\right) \ds_1 \ds_2, \\
(f_{3})_{i,j,k} & = \sum_{\X_{l,m} \in \Omega_b} (F_{3})_{l,m} \psi_h\left(\x_{i,j,k} , \X_{l,m}\right) \ds_1 \ds_2. \label{eqn_f3_spread_one_sided}
\end{align}
Notice that the above formulas are same as those employed for regular IB interpolation and spreading operations, except that here we use MLS weighting functions.

\subsection{Nullspace of the polynomial matrix} \label{sec_nullspace}

An interesting property of the MLS generating functions is that the odd-degree polynomial reproduction conditions 
are not affected by a constant shift, given a symmetric arrangement of the interpolation points $\x$ 
around the evaluation point $\X$. The following Lemma makes it precise by considering two spatial dimensions and 
univariate polynomials.

\begin{mylem}
\label{lem:constant_shift}
With $N$ denoting the number of Eulerian points in the support of the weight function in both $\Omega_b^{+}$ and 
$\Omega_b^{-}$ regions, 
and $\V C = c \V 1$ denoting a constant vector with entries $c$, the modified weights 
$\V{\Psi}^m = (\V{\Psi}+ \V C)/\InterpSumi (\V \Psi + \V C) = (\V{\Psi}+ \V C)/(1+ Nc)$ satisfy the 
first moment condition if the evaluation point is coincident with one of the interpolation points.
\begin{proof}
Since $\V{\Psi}$ satisfies the zeroth-order moment condition by construction, $\sum \V{\Psi} = 1$ and since $\V C $ is a 
constant vector $\sum \V C = Nc$.  Therefore, the denominator of $\V{\Psi}^m$ is $\sum (\V \Psi + \V C) = 1 + Nc$.  
Next, considering Eq.~\eqref{eqn_linsys}, the modified weights satisfy the relation 
\begin{equation}
\cA \V{\Psi}^m  =  \V{P}^m =  (\V P + \cA \V C)/(1+ Nc).  \label{eqn_psim}
\end{equation}     
In two spatial dimensions with univariate polynomials, the $\cA \V C$ term of $\V{P}^m$ is
\begin{equation}
  \cA \V C = c (\cA \; \V 1) =  \begin{bmatrix}                             
           Nc  \\
           c \sum_{i=1}^N (x_i - X)  \\
           c \sum_{i=1}^N (y_i - Y) 
         \end{bmatrix}  
\end{equation}
If the evaluation point $\X = (X, Y)$ is coincident with one of the interpolation points, odd-degree terms like  $\sum_{i=1}^N (x_i - X)$ or 
$\sum_{i=1}^N (y_i - Y)$ evaluate to zero because of the symmetric arrangement. With $\V P = (1, 0, 0)^T$ 

\[	\V{P}^m =  (\V P + \cA \V C)/(1+ Nc) = \begin{bmatrix}                             
           									1  \\
          									0  \\
           									0 
         									\end{bmatrix}  = \V P.
\].   
\end{proof}
\end{mylem}
Consequently, $\V{\Psi}^m$ satisfies the same discrete moment conditions as  $\V{\Psi}$ for a symmetric arrangement of interpolation 
points centered around a moving point of interest $\X$. With a cell-centered discretization, this assumption holds for an IB point that is coincident with cell centroid. An analogous argument holds in three spatial dimensions with univariate 
polynomials.  Notice that the zeroth-order moment condition of $\V{\Psi}^m$ is satisfied due to the normalization factor $(1 + Nc)$, irrespective
of the relative location of the evaluation point with respect to the background grid. We remark that higher-order moments, 
e.g., second-order moments with bivariate polynomials, won't be satisfied exactly even when 
$\X$ is collocated on cell centroid for a constant shift. This is because even-degree entries of $\V P^m$, such as   
$\sum_{i=1}^N (x_i - X)^2$, do not provide cancellations. However, 
odd-degree terms like $\sum_{i=1}^N (x_i - X)(y_i - Y)$ would still equal zero.

\subsection {Mollifying one-sided MLS kernels}
\label{subsec_positivity}
In an essence, the moving least squares solution transforms the input weights $W_{\text{MLS}}(\x_i,\X)$ 
to $\psi_h(\x_i,\X)$ that satisfy the polynomial reproduction conditions. For IB spreading operation, it is desirable
that generating functions $\psi_h(\x_i,\X)$ are monotonically decreasing and remain positive (or have negligible negative tails) for all evaluation points $\X$. In fact, the new/modified five- and six-point Peskin's 
delta functions were designed to eliminate the negative tails of the standard five- and six-point IB kernels, respectively. 

\begin{figure}
	\centering
	\subfigure[ $\cW$]{
		\includegraphics[scale=0.23]{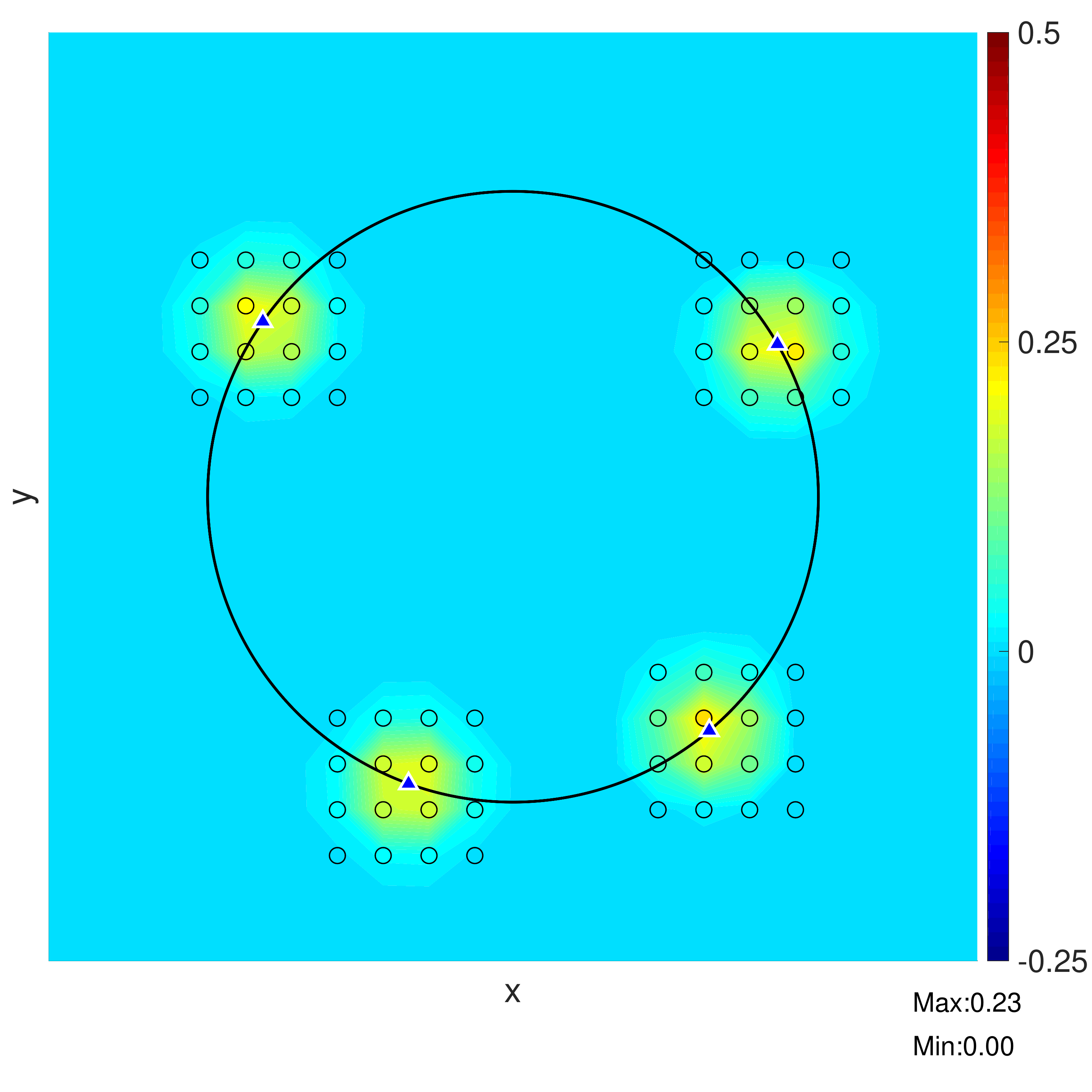} 
		\label{fig:MLSProc:W}
	}
	\subfigure[$ \cW_{\text{MLS}}$]{
		\includegraphics[scale=0.23]{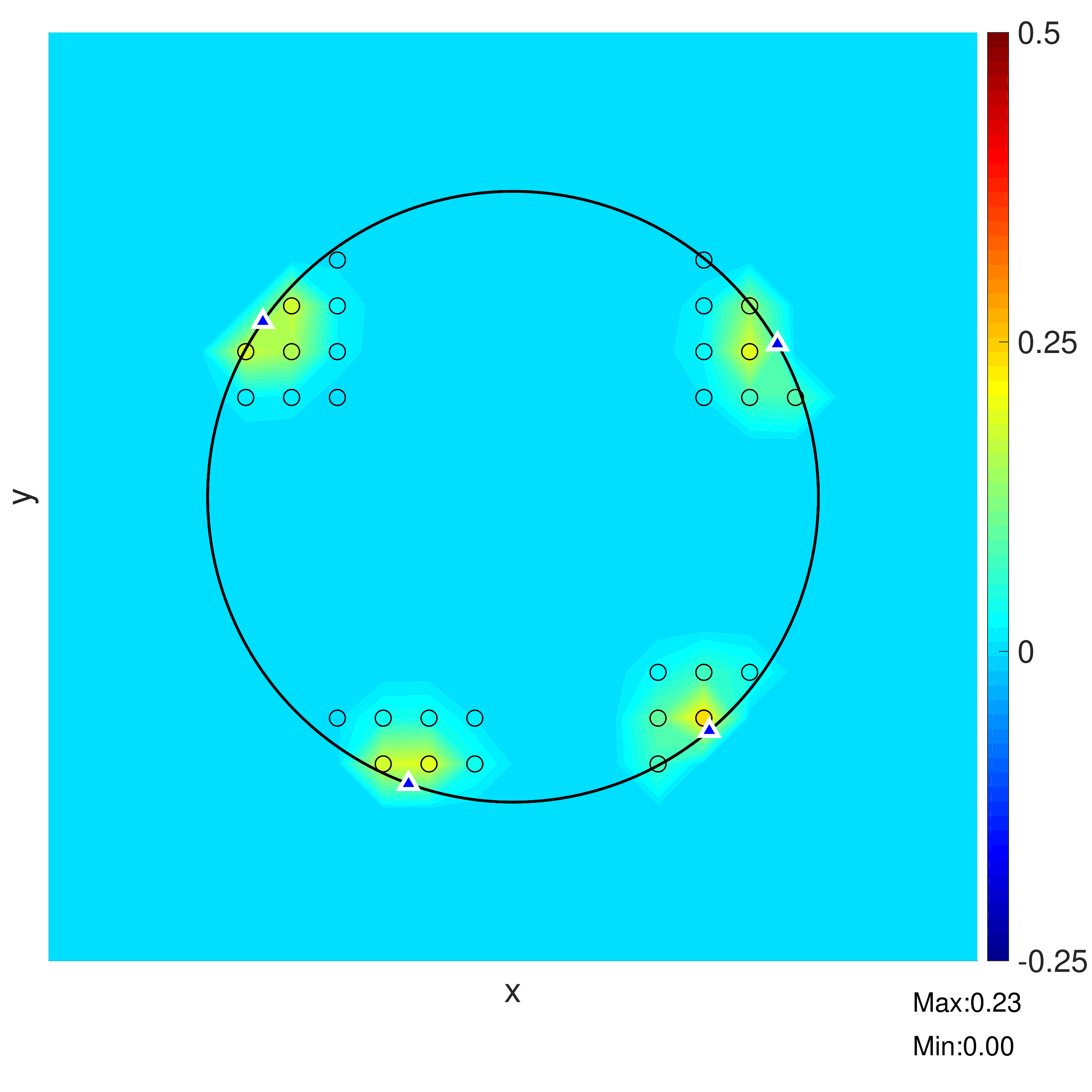} 
		\label{fig:MLSProc:WMLS}
	}
	\subfigure[$\V{\Psi}$]{
		\includegraphics[scale=0.23]{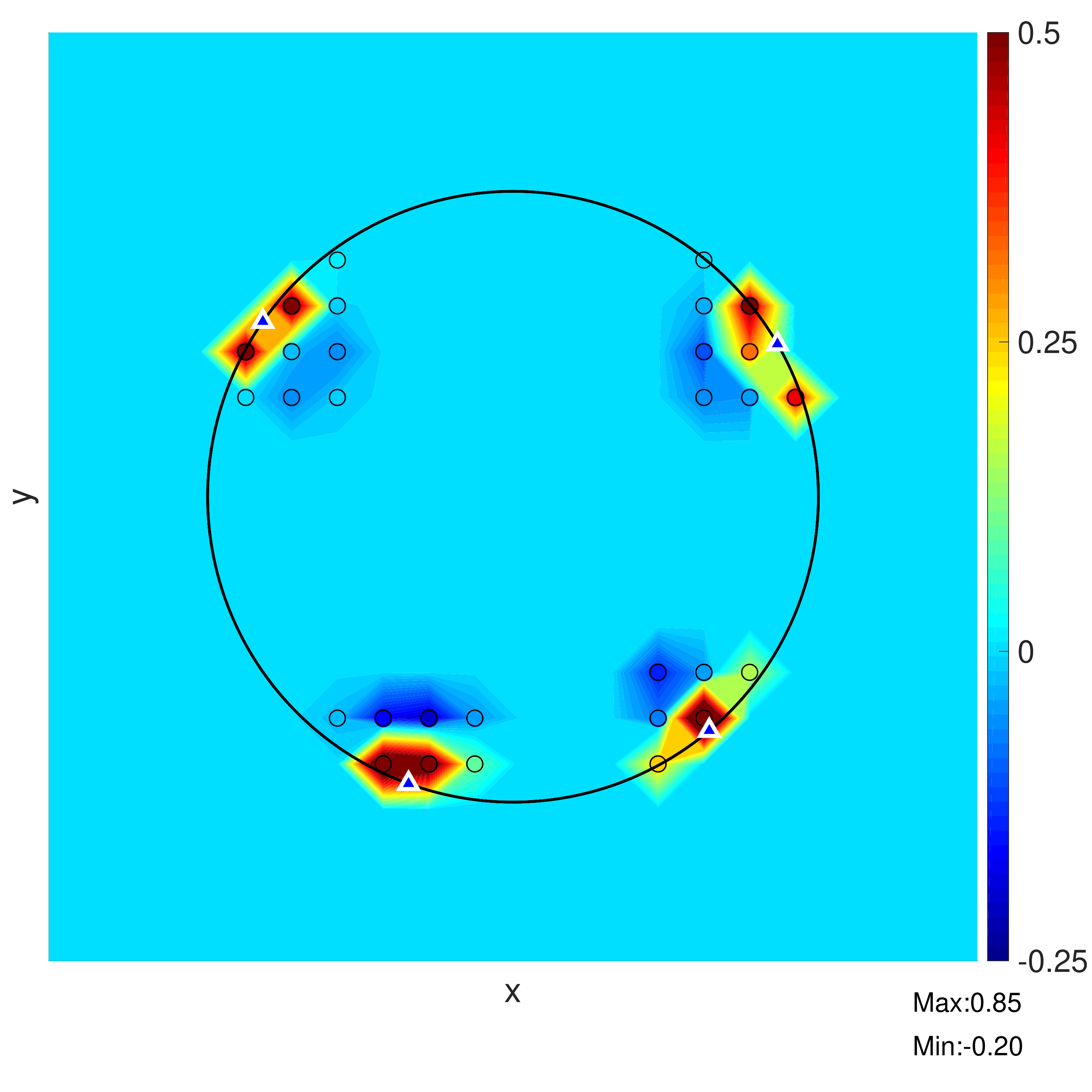} 
		\label{fig:MLSProc:Psi}
	}
	\subfigure[$\V{\Psi}^{m}_{\text{CVS}}$]{
		\includegraphics[scale=0.23]{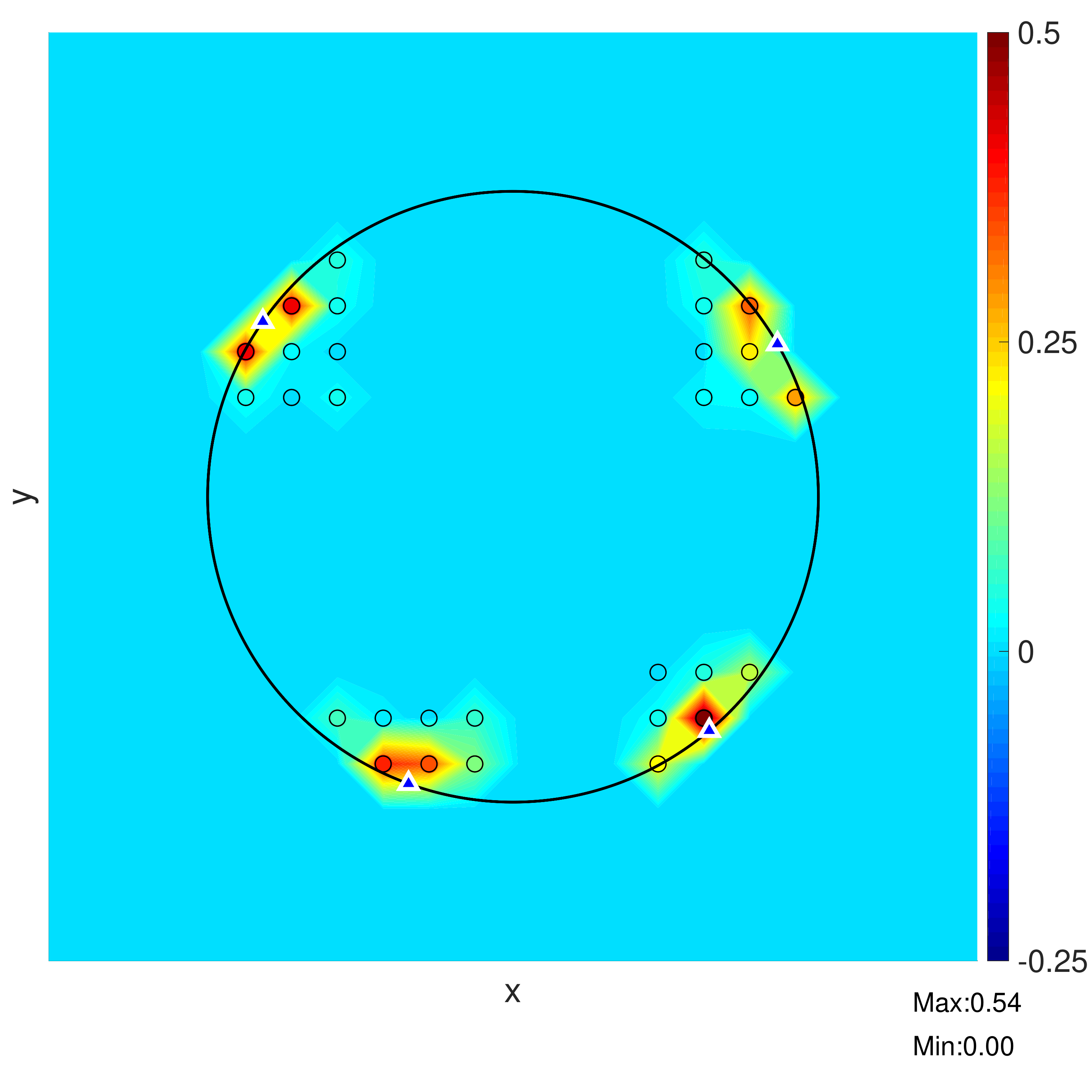} 
		\label{fig:MLSProc:PsiCVS}
	}
	\subfigure[$\V{\Psi}^{m}_{\text{NCVS}}$]{
		\includegraphics[scale=0.23]{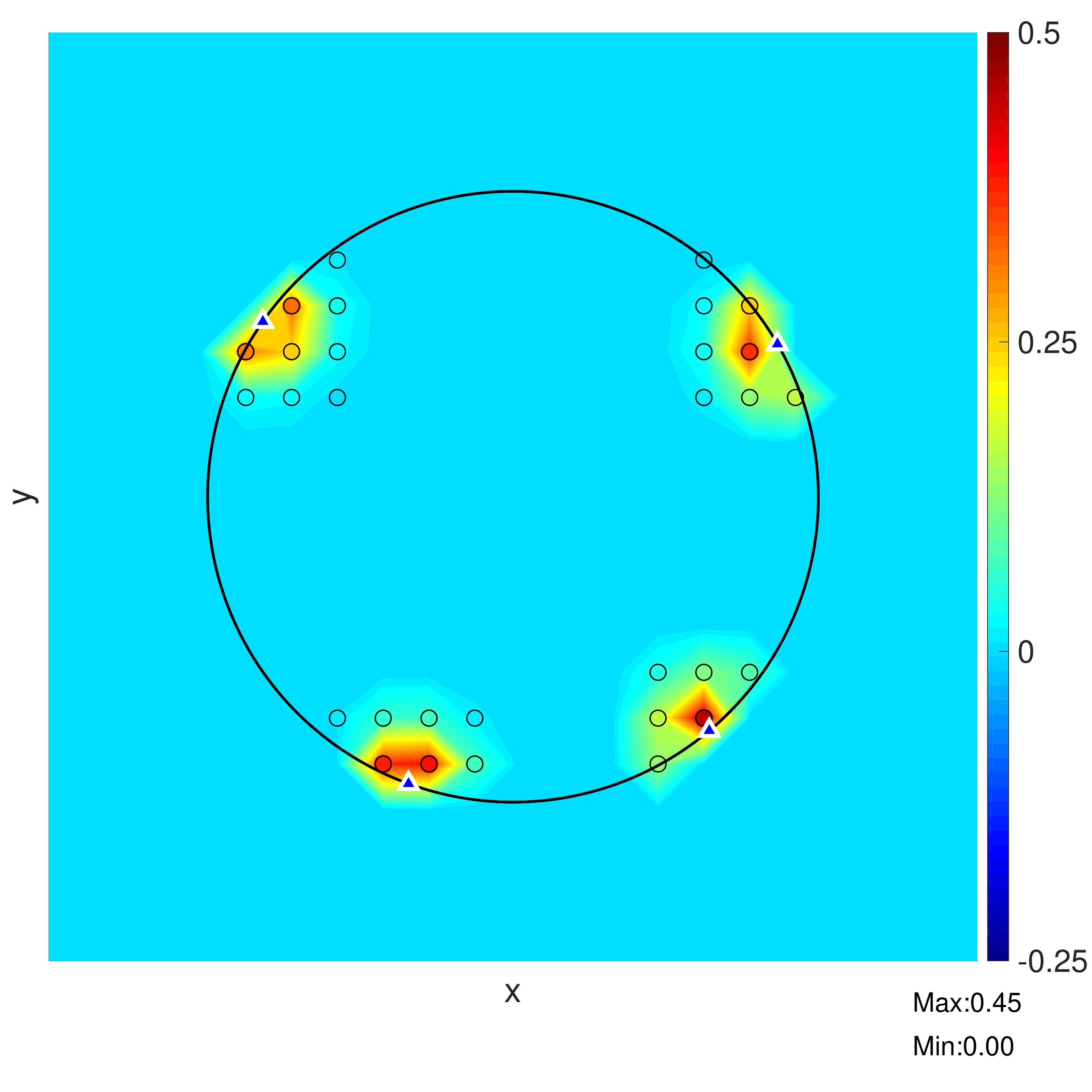} 
		\label{fig:MLSProc:PsiNCVS}
	}
	
	\caption{Color plot of~\subref{fig:MLSProc:W} four-point delta function weights;~\subref{fig:MLSProc:WMLS} restricted four-point delta function weights;~\subref{fig:MLSProc:Psi} MLS generated weights;~\subref{fig:MLSProc:PsiCVS} CVS shifted MLS weights; and~\subref{fig:MLSProc:PsiNCVS} 
		NCVS shifted MLS weights. Here, Eulerian grid points outside of the cylinder are masked.}
	\label{fig:MLSProc}
\end{figure}

The increased positivity of kernels also lead to better grid translation invariance property as shown 
empirically by Bao et al.~\cite{Bao2016}. In the context of one-sided IB kernels, the MLS procedure can generate larger (compared to weighting function) weights for nearby Eulerian grid nodes and negative weights for far-away grid nodes to satisfy the linear conditions. Fig.~\ref{fig:MLSProc} in the next section (Sec.~\ref{sec_mlsproc}) describes this situation. We remark that large positive and negative weights in the kernel function do not pose stability concerns for certain IB methods that compute constraint forces exactly. It however, can induce flow instabilities in a direct forcing IB simulation through an oscillatory feedback. Indeed this was observed for test cases simulated using direct forcing IB method in this work. We also tested the one-sided MLS kernels  with exactly-constrained IB methods~\cite{Kallemov16,Usabiaga17}. Our preliminary tests suggest that the direct use of one-sided MLS kernels with exactly-constrained IB methods do not appear to produce flow instabilities; such specialized IB methods are not the focus of the current work and shall be explored more in the future. 

We explore two \emph{weight-shifting} strategies to mollify the generating functions $\psi_h(\x_i,\X)$. 
The shifting procedure is motivated by the constant nullspace property defined in the previous section.  
We shall compare the two shifting strategies for flow past a sphere  
at $Re =  10,000$ in Sec.~\ref{sec_cvs_ncvs_perf} and make recommendations based on physical observations for this test case.

In Sec.~\ref{sec_nullspace}, we examined shifting the generating functions $\V \Psi$ by a constant vector. It was noted that the constant shift preserves 
the linear polynomial reproduction constraint under certain conditions, although doing so eliminates the one-sided 
property of  $\psi_h(\x_i,\X)$; $\V{\Psi}^m(\x_i,\X) = c \ne 0$ for $\x_i \in \Omega_b^{-}$.  To remedy this situation we propose 
two modifications: 
\begin{enumerate}
\item a constant vector shift (CVS) based on a restricted version of constant vector $\Comegap = c \V 1_{\Omega_b^{+}}$;
\item a more general non-constant vector shift (NCVS) based on MLS weights $\Comegap = c \cW_{\text{MLS}}$.
\end{enumerate} 
The form of the modified generating functions remain the same, which read as 
\begin{align}
\V{\Psi}^m_{\text{CVS}} &= \dfrac{ \V{\Psi}+  c \V 1_{\Omega_b^{+}} }{ \sum (\V \Psi + c \V 1_{\Omega_b^{+}}) }, \label{eqn_psim_cvs} \\
\V{\Psi}^m_{\text{NCVS}} &= \dfrac{ \V{\Psi}+ c \cW_{\text{MLS}} } {\sum (\V \Psi + c \cW_{\text{MLS}} )} = \dfrac{ \cW_{\text{MLS}} \odot (\vcL + c \V 1) }{ \sum (\V \Psi + c \cW_{\text{MLS}}) }. \label{eqn_psim_ncvs}
\end{align}
Notice that both shifts are zero in the $\Omega_b^{-}$ region which preserves the one-sidedness of $\V{\Psi}^m$.  
  
We now consider the effect of two restricted shifts on the polynomial reproduction constraints.
Let $N = N^{+} \cup N^{-}$ denote the number of Eulerian points in the support of the weight function, 
with $N^{+} \in \Omega_b^{+}$ and  $N^{-} \in \Omega_b^{-}$.  Using the zeroth moment condition of original $\V{\Psi}$, 
and the more general $\Comegap = c \cW_{\text{MLS}}$, we have
\begin{equation}  
\V{\Psi}^m = \dfrac{\V{\Psi}+ \Comegap}{1 + c \sum_{i \in N^{+}}  W_i},  \label{eqn_psim_restricted}
\end{equation}
in which $W_i = W_{\text{MLS}}(\x_i,\X)$. Using Eq.~\eqref{eqn_linsys} for the modified weights $\V{\Psi}^m$, 
\begin{equation}
\cA \V{\Psi}^m  =  \V{P}^m = \dfrac{ \V P + \cA \Comegap}{1+ c \sum_{i \in N^{+}} W_i}.     \label{eqn_psim_ws}
\end{equation}     
Again considering two spatial dimensions and univariate polynomials, the $\cA \Comegap$ term of $\V{P}^m$ is
\begin{align}
\cA \Comegap & = c  \begin{bmatrix}                             
\sum_{i \in N^{+}} W_i  \\
\sum_{i \in N^{+}}   (x_i - X) W_i  \\
\sum_{i \in N^{+}}   (y_i - Y) W_i  
\end{bmatrix}.
\end{align}
With $\V P = (1, 0, 0)^T$,
\begin{align}
\V{P}^m & =   \begin{bmatrix} [2.1]                            
1  \\
\dfrac{c \sum_{i \in N^{+}}  (x_i - X) W_i }{1+ c \sum_{i \in N^{+}} W_i}  \ne 0  \\
\dfrac{c \sum_{i \in N^{+}}  (y_i - Y) W_i}{1+ c \sum_{i \in N^{+}} W_i}  \ne 0
\end{bmatrix}.   \label{eqn_pm_restricted}
\end{align}
Similarly, for $\Comegap = c \V 1_{\Omega_b^{+}}$,  $\V{P}^m$ can be obtained by simply setting $W_i = 1$ in Eq.~\ref{eqn_pm_restricted}. 
We make several remarks about these two weight shifting strategies:

\begin{itemize}

\item For the Heaviside function, $W_i = 0$ for $\x_i \in \Omega_b^{-}$, so only data in $\Omega_b^{+}$ contribute to the summation. Hence,
odd-degree polynomial terms may not cancel fully. Therefore, both strategies do not satisfy the linear polynomial reproducing condition in general. 
In practice, the IB marker points are not collocated on cell centroids, especially for moving bodies. 
Therefore, the modified weights do not satisfy the first moment condition even for the unrestricted 
constant shift $\V C = c \V 1$ for practical FSI applications. 

\item For the NCVS strategy, the shifts are inversely proportional to the distance between the evaluation and interpolation points. The CVS strategy is a special case of the NCVS strategy, with an  
equal shift applied to each interpolation point as $W_i = 1$.  We will see the physical manifestation of these two shifting 
strategies in the context of flow past sphere in Sec.~\ref{sec_cvs_ncvs_perf}.

\item The scalar $c$ is a free-parameter, which can be chosen  to enhance the positivity of $\V\Psi$.  In this work we choose 
$c = |\text{min}(0, \text{min}(\V{\Psi}))|$ and $c = |\text{min}(0, \text{min}(\vcL))|$ for the CVS and NCVS strategy, respectively. From Eqs.~\eqref{eqn_psim_cvs} and~\eqref{eqn_psim_ncvs},  it can be easily verified that these are the minimum values of $c$ that ensure $\V{\Psi}^m \ge 0$.

\item The $\V{\Psi}^m_{\text{NCVS}}$ weights in Eq.~\eqref{eqn_psim_ncvs} can be viewed as a specific version of a more general 
combination of moment-satisfying one-sided MLS weights $\V{\Psi}$ and a restricted version of two-sided 
weights $\cW_{\text{MLS}}$ (that may originally also satisfy the discrete moments),

\begin{equation}
\V{\Psi}^m_{\text{NCVS}} = \dfrac{ \alpha \V{\Psi}+ \beta \cW_{\text{MLS}} } {\sum ( \alpha \V \Psi + \beta \cW_{\text{MLS}} )}.
\end{equation}
Here $\alpha$ and $\beta$ can be selected to satisfy the desired properties for the weights $\V{\Psi}^m$. The specific values of $\alpha = 1$ and 
$\beta = |\text{min}(0, \text{min}(\vcL))|$ considered here ensure that the weights are positive and revert to the MLS weights $\V{\Psi}$ to achieve reproducing conditions, whenever possible.

\end{itemize}

\subsection{MLS procedure} \label{sec_mlsproc}

 In this section we take a simple case of a two-dimensional closed circular interface embedded in a uniform Cartesian grid to demonstrate the procedure 
 of generating original and modified MLS weights, $\V{\Psi}$ and $\V{\Psi}^m$, respectively. This example also serves to demonstrate that the MLS method 
 can lead to over- and undershoots in $\V{\Psi}$ for the restricted domain problem. It also serves to illustrate the process of remedying $\V{\Psi}$ through  CVS and NCVS shifting. For generality, the procedure is demonstrated for four Lagrangian points placed non-symmetrically along the interface:

\begin{figure}
	\centering
	\subfigure[CVS]{
		\includegraphics[width=0.3\textwidth]{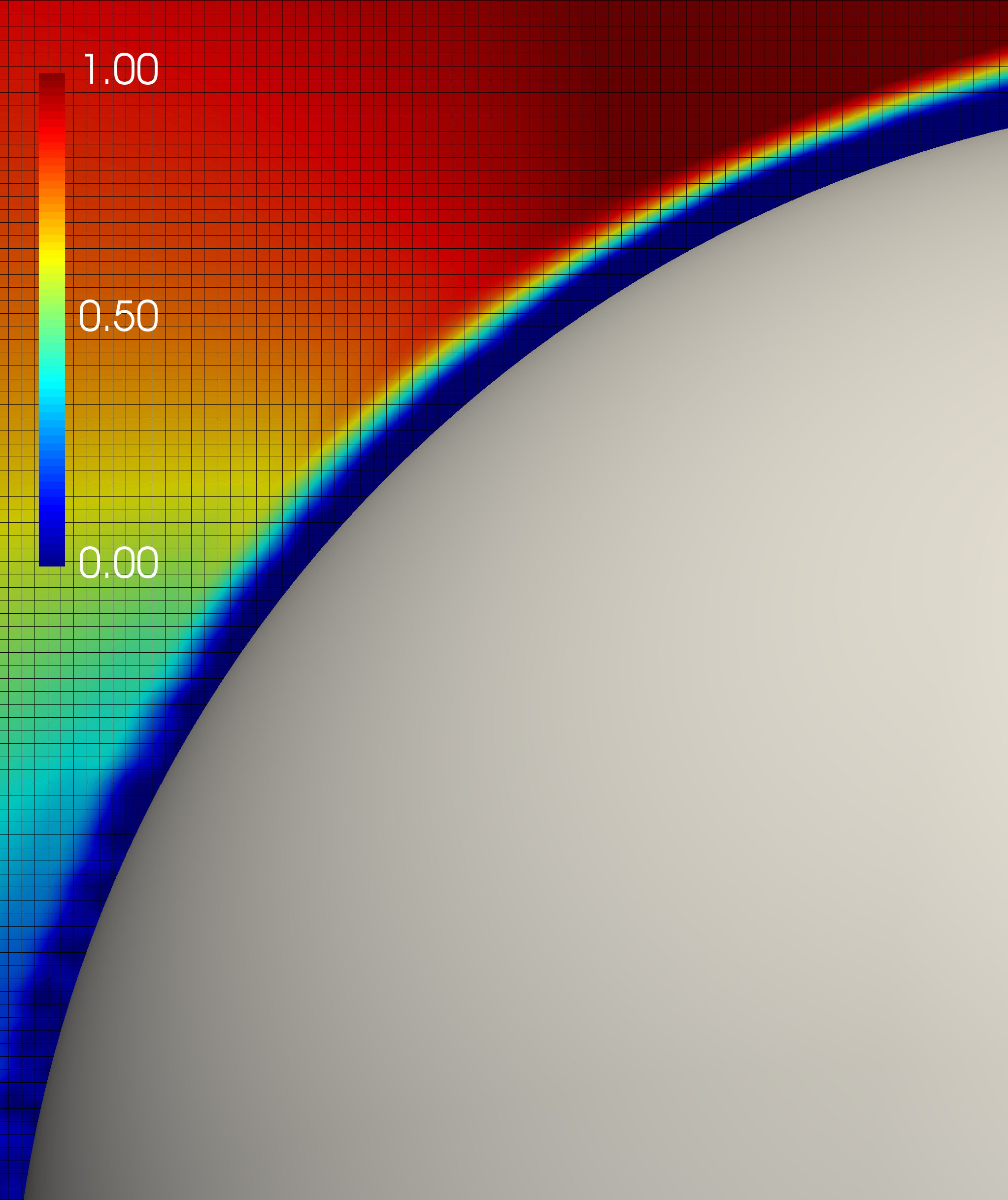} 
		\label{fig:sphere:cvs}
	}
	\subfigure[NCVS]{
		\includegraphics[width=0.3\textwidth]{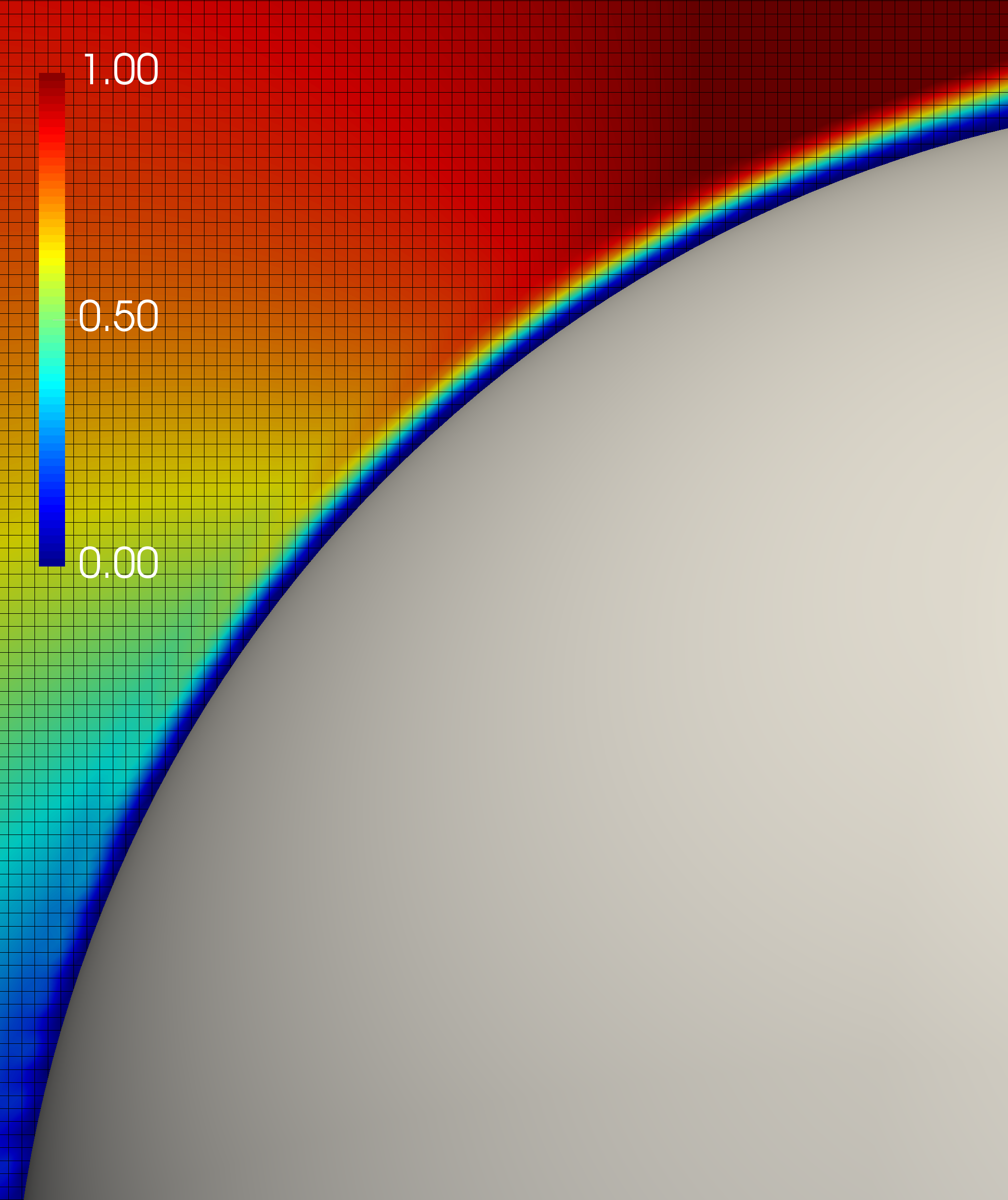} 
		\label{fig:sphere:ncvs}
	}
	\caption{Flow past a sphere at $Re = 10,000$ using~\subref{fig:sphere:cvs} $\V{\Psi}^{m}_{\text{CVS}}$ and~\subref{fig:sphere:ncvs}  $\V{\Psi}^{m}_{\text{NCVS}}$ weights. The former approach leads to a thick numerical boundary layer, which is unphysical. }
	\label{fig:sphere:cvs-ncvs}
\end{figure}

 \begin{enumerate}
 	\item  The first step involves generating the $\cW$ weights at the Eulerian grid points that are in the support of  the Lagrangian material points. 
Fig.~\ref{fig:MLSProc:W} shows the color plot of Peskin's four-point delta function weights. The Lagrangian points are highlighted with a filled triangle and the neighboring Eulerian grid points in the compact support are highlighted through black circles.
 	
 	\item Next, the restricted weights $\cW_{\text{MLS}}$ are evaluated using a masking criterion. In this demonstration we mask the Eulerian grid points outside of the cylinder. The restricted weights are shown in Fig.~\ref{fig:MLSProc:WMLS}. 
 	
 	\item The restricted weights $\cW_{\text{MLS}}$ are then used to produce the generating weights $\V{\Psi}$, which are plotted in Fig.~\ref{fig:MLSProc:Psi}.  As mentioned in Sec.~\ref{subsec_positivity}, and also for this example, the MLS construction produces large positive 
weights near the immediate vicinity of the interface and negative values near the tail-end region of $\V{\Psi}$; compare the color scale of Figs.~\ref{fig:MLSProc:W} and~\ref{fig:MLSProc:Psi}. 	
 	
 	\item To mollify $\V{\Psi}$,  the modified MLS weights $\V{\Psi}^m$ are evaluated using CVS and NCVS  approaches described in the prior Sec.~\ref{subsec_positivity}. The modified weights $\V{\Psi}^{m}_{\text{CVS}}$ and $\V{\Psi}^{m}_{\text{NCVS}}$ are plotted in Fig.~\ref{fig:MLSProc:PsiCVS} and Fig.~\ref{fig:MLSProc:PsiNCVS}, respectively. It can be seen that the modified weights $\V{\Psi}^m$ are positive towards the tail-end region and have mollified near the immediate vicinity of the interface. 
 	
\end{enumerate}

\subsection{CVS and NCVS comparison} \label{sec_cvs_ncvs_perf}

While both approaches mollify $\V{\Psi}$, there is an important difference in their weight distribution. $\V{\Psi}^{m}_{\text{NCVS}}$ has monotonically decreasing weights away from the evaluation point, which is also consistent with the distribution of the original weight $\cW$. In contrast, $\V{\Psi}^{m}_{\text{CVS}}$ has a weight distribution that does not decrease monotonically with distance. For IB simulations, it is desirable that weights have an inverse relationship with distance. If this condition is not met, the interpolation weights can artificially increase the  boundary layer thickness at the fluid-IB interface. 

To elucidate this effect, we consider flow past a sphere at a Reynolds number of $10,000$.  At $ Re = 10,000 $, 
the boundary layer is expected to be very thin and an IB simulation is expected to reflect this behavior. Peskin's four-point regularized delta kernel is used for evaluating $\cW_{\text{MLS}}$. Velocity magnitude of the flow near the stagnation region of the sphere is shown in Fig.~\ref{fig:sphere:cvs-ncvs} for the two mollifying schemes.  Because of the non-monotonically decreasing $\V{\Psi}^{m}_{\text{CVS}}$ weights, grid nodes lying towards the tail-end region possess relatively large weights (see Fig.~\ref{fig:MLSProc:PsiCVS}). The resulting magnitude of the IB force at these grid nodes is large (compared to standard IB kernels), and consequently it leads to a thick numerical boundary layer as illustrated in Fig.~\ref{fig:sphere:cvs}. In contrast, flow at the fluid-IB interface corresponding to NCVS scheme is typical of an external flow around a sphere using a regular diffuse-interface IB method. \REVIEW{The boundary layer over the sphere using the standard four-point
kernel is quite similar to Fig.~\ref{fig:sphere:ncvs}, and the comparison is omitted for brevity. Although this comparison is qualitative, and more systematic investigation is needed for high $Re$ flows than what is presented in this section, quantitative comparison of CVS and NCVS approaches at moderate $Re$ flows is demonstrated for few cases in Sec.~\ref{sec_results}. The results of Sec.~\ref{sec_results} also favor the NCVS approach.}  Therefore, NCVS shifting scheme is preferred over CVS scheme, and all the results in what follows will be based on $\V{\Psi}^{m}_{\text{NCVS}}$, unless stated otherwise.  As discussed in the introduction, the unmodified $\V{\Psi}$ weights produce unphysical flow oscillations (due to force spreading) at the fluid-IB interface which destabilizes the direct forcing IB simulation; flow instability data not presented for this case.

\section{Software Implementation}
The immersed boundary methods presented in this work are implemented in the solver \textsc{Cube}~\cite{Jansson18}, which is a multi-physics flow solver for massively parallel simulations. It is built on a hierarchical meshing technique known as the Building Cube Method (BCM)~\cite{naka03}. \textsc{Cube} supports hybrid parallelism through a combination of  Message Passing Interface (MPI) and shared-memory parallelism using OpenMP. 


\section{Results} \label{sec_results}

In this section, we first benchmark the one-sided IB kernels on problems from the IB literature to check the accuracy of FSI solutions. 
We also highlight the ability of one-sided kernels to avoid spurious flows inside the solid domain that are typically produced 
by standard diffuse IB kernels. Next, the MLS kernels are used to simulate flow past the Ahmed body~\cite{ahme84}, which demonstrates the 
potential of the one-sided IB approach to treat complex engineering geometries.

\begin{figure}
	\centering
	\subfigure[Interior and exterior IB forcing with Peskin's four-point delta function]{
	\includegraphics[scale=0.32]{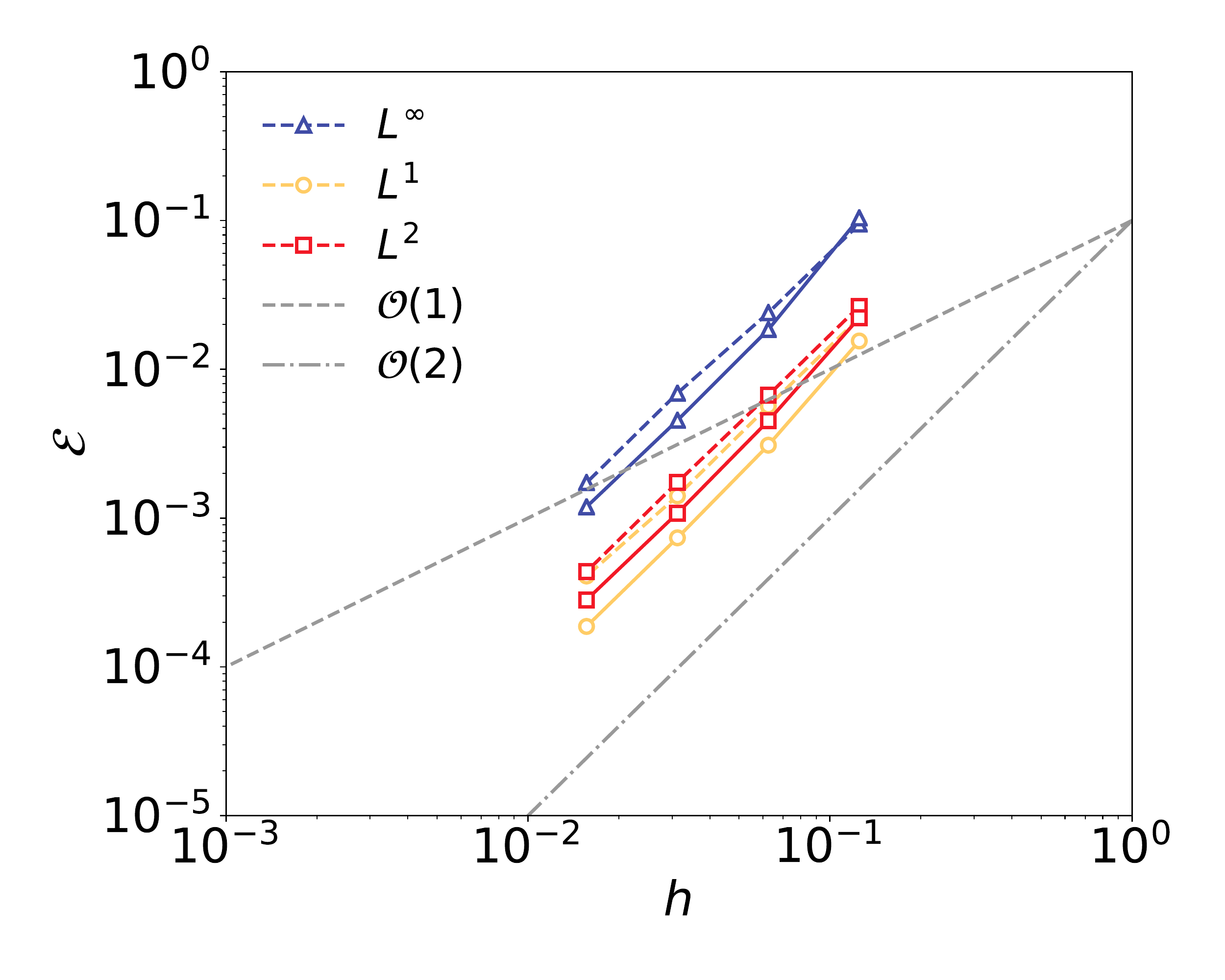} 
	\label{fig:TG:error1cIB}
    }
	\subfigure[Interior and exterior IB forcing with $\V{\Psi}^{m}_{\text{NCVS}}$ in $\cJ$ and $\cS$]{
		\includegraphics[scale=0.32]{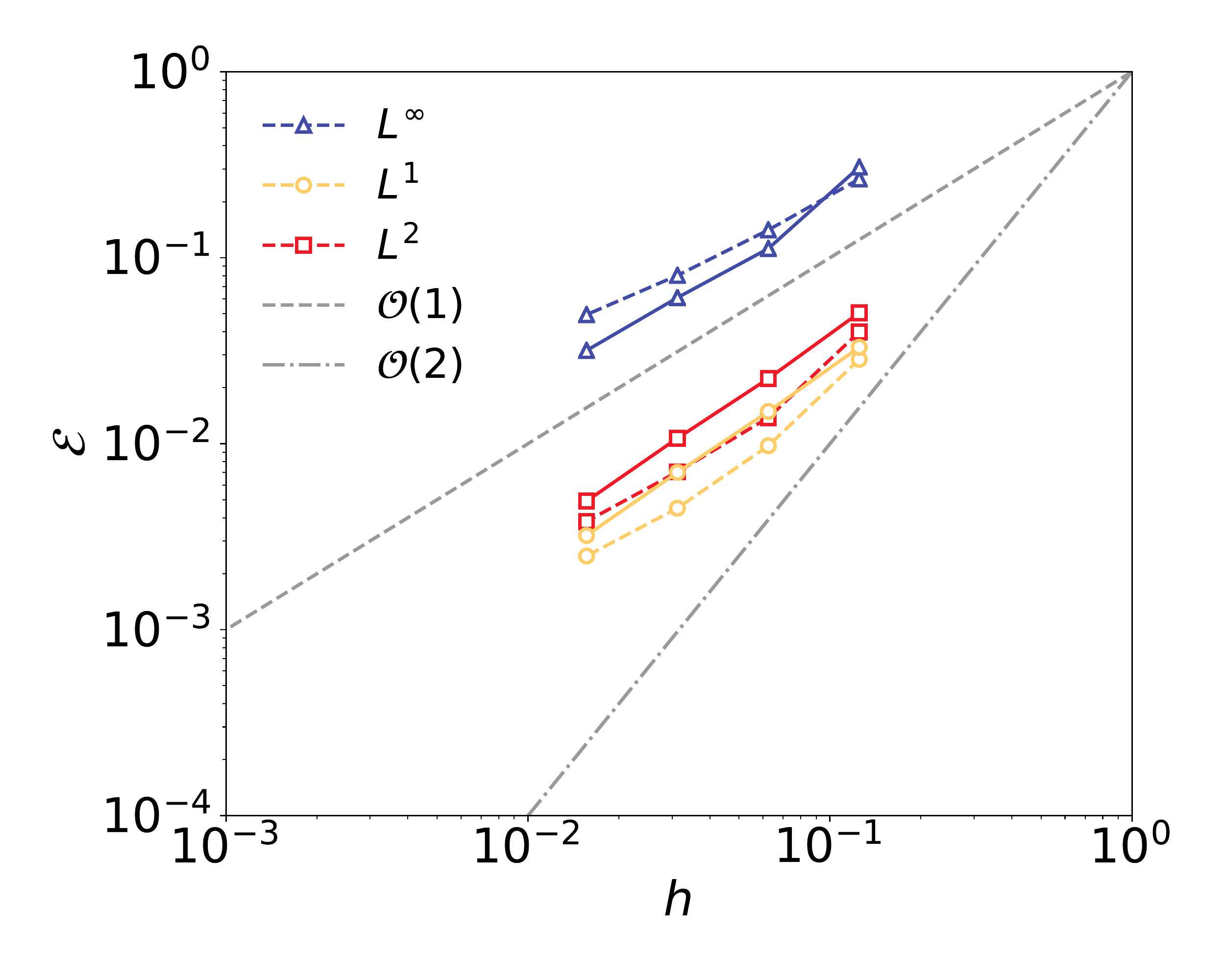} 
		\label{fig:TG:error1}
	}
	\subfigure[Exterior IB forcing with $\V{\Psi}^{m}_{\text{NCVS}}$ in $\cJ$ and $\cS$]{
		\includegraphics[scale=0.32]{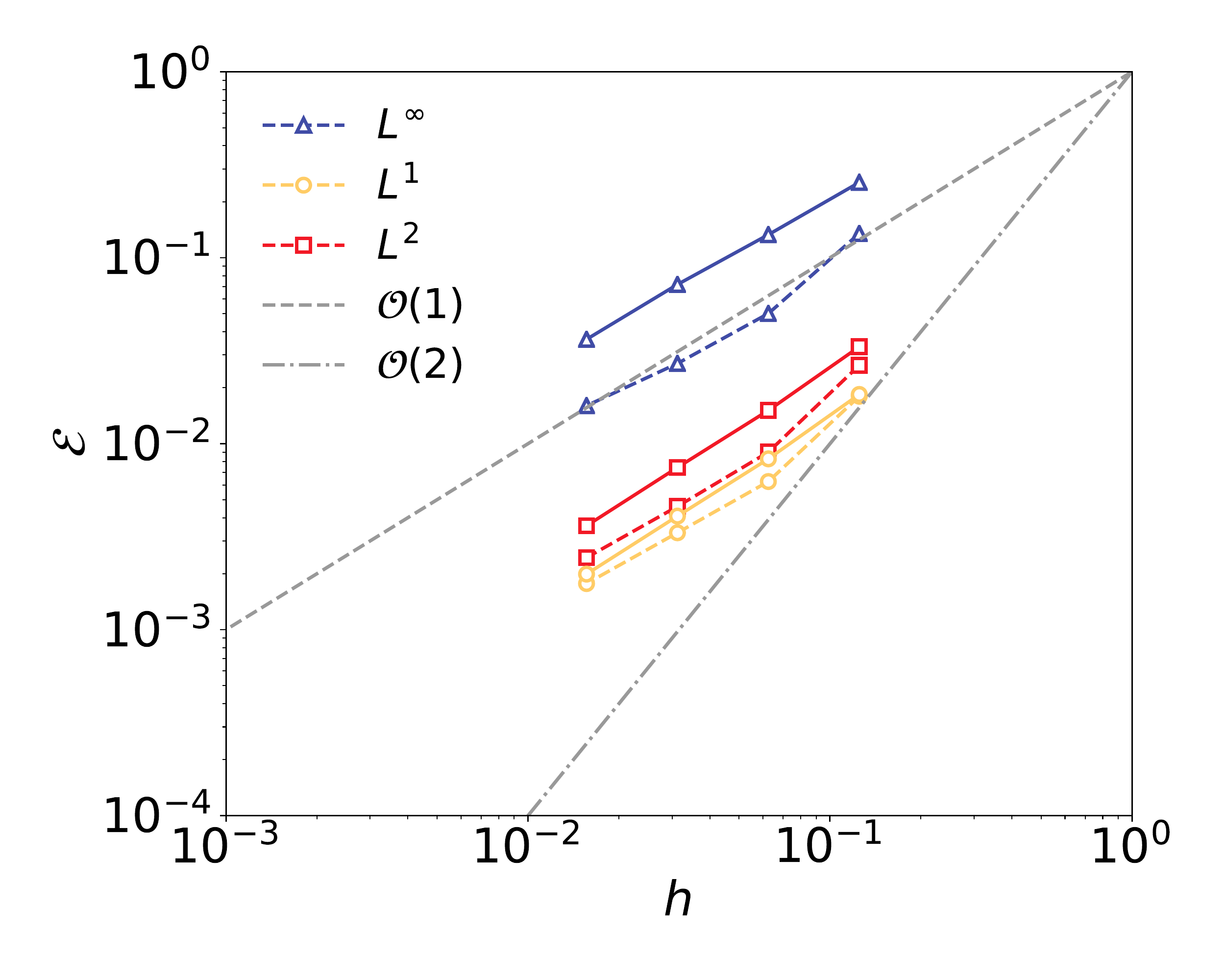} 
		\label{fig:TG:error2}
	}
        \subfigure[Interior and exterior IB forcing with $\V{\Psi}^{m}_{\text{CVS}}$  in $\cJ$ and $\cS$]{
		\includegraphics[scale=0.32]{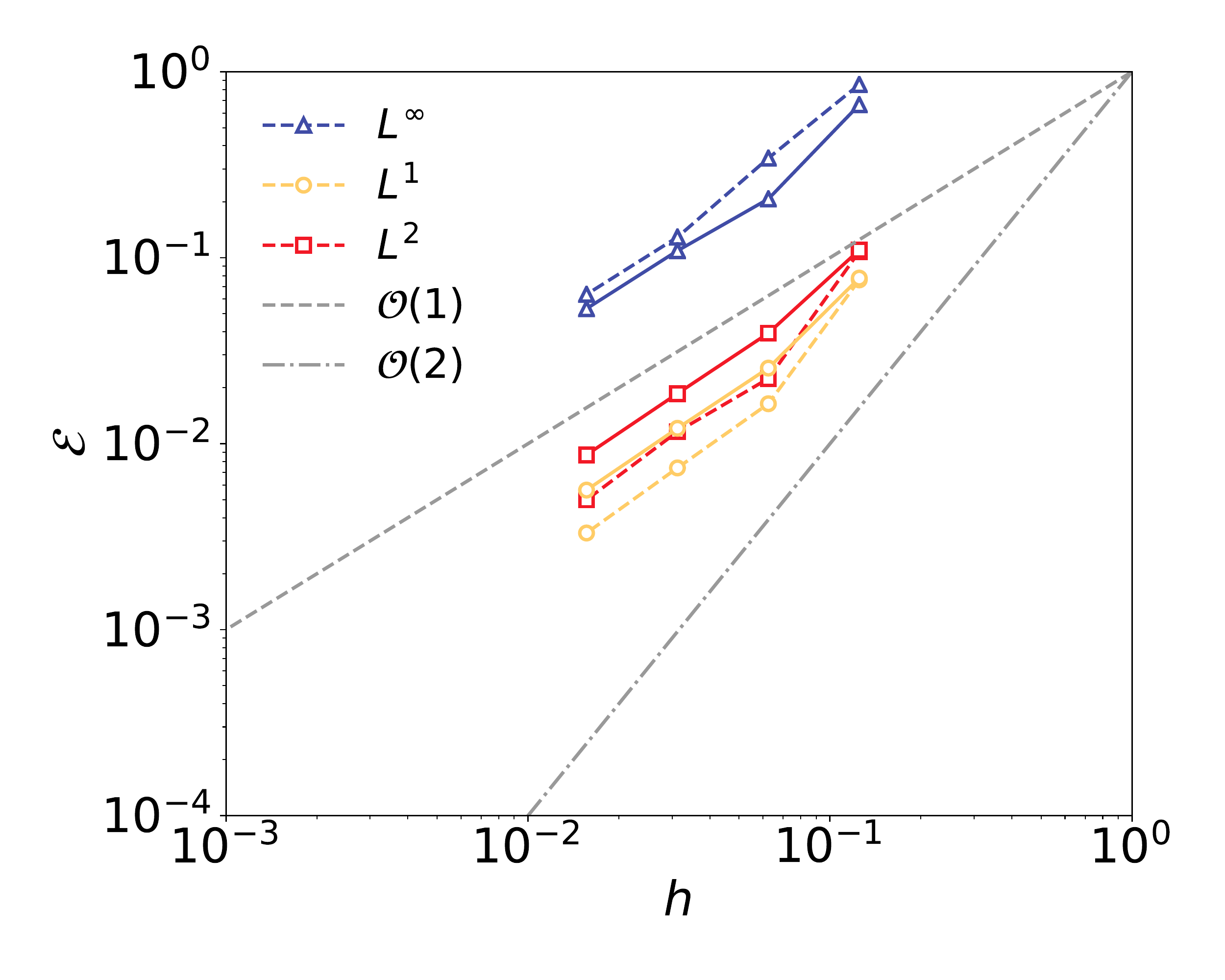} 
		\label{fig:TG:error1CVS}
	}
	\subfigure[Interior and exterior IB forcing with $\V{\Psi}$ in $\cJ$ and $\V{\Psi}^{m}_{\text{NCVS}}$ in $\cS$]{
	\includegraphics[scale=0.32]{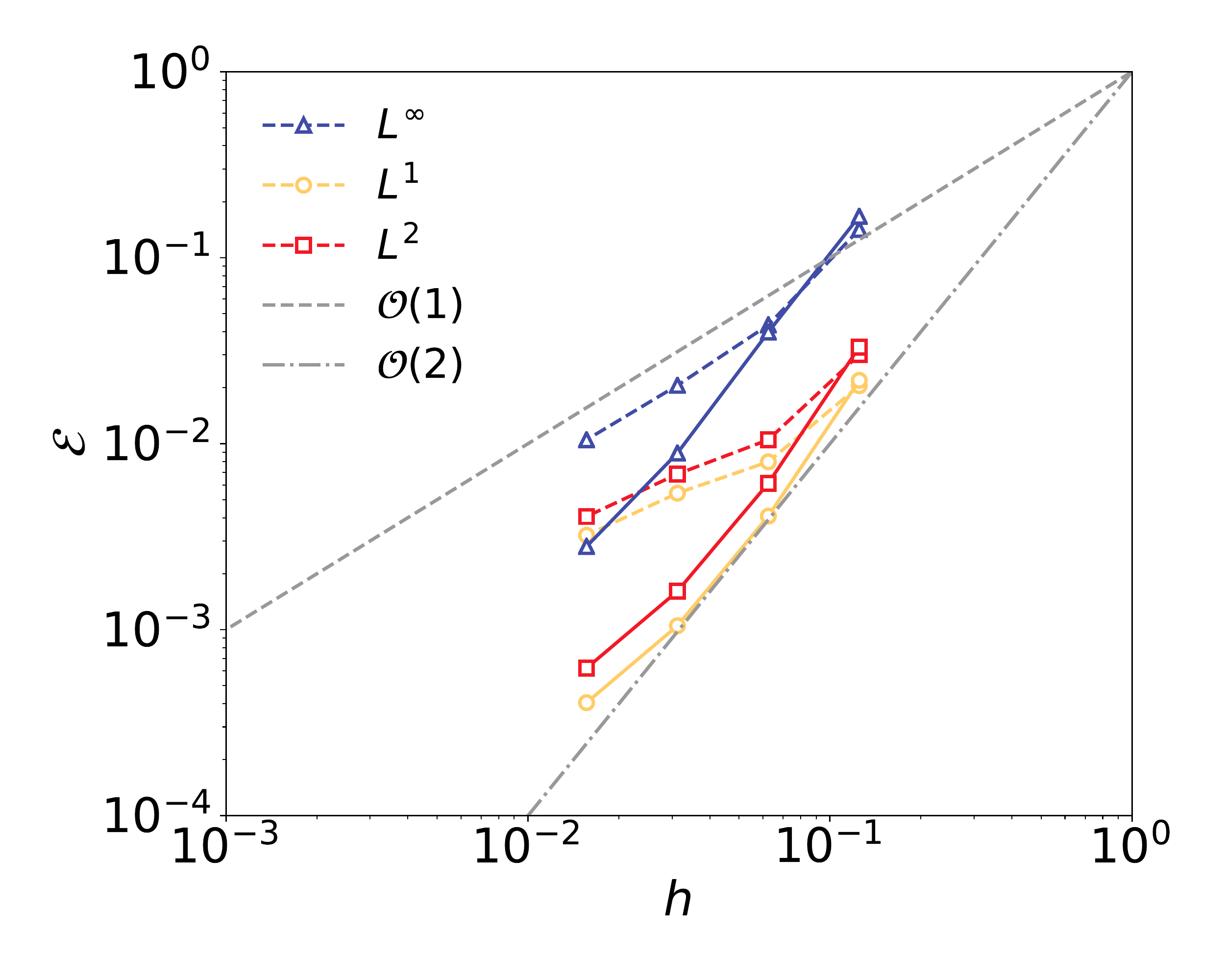} 
	\label{fig:TG:errorIpsiSncvs}
       }
	\caption{$L^1$, $L^2$, and $L^{\infty}$ error norms of $u$ velocity component (solid lines) and pressure (dashed lines) for the Taylor-Green vortex problem plotted against the grid cell size.}
	\label{fig:TGError}
\end{figure}

\begin{figure}
	\centering
	\subfigure[]{
		\includegraphics[width=0.24\textheight]{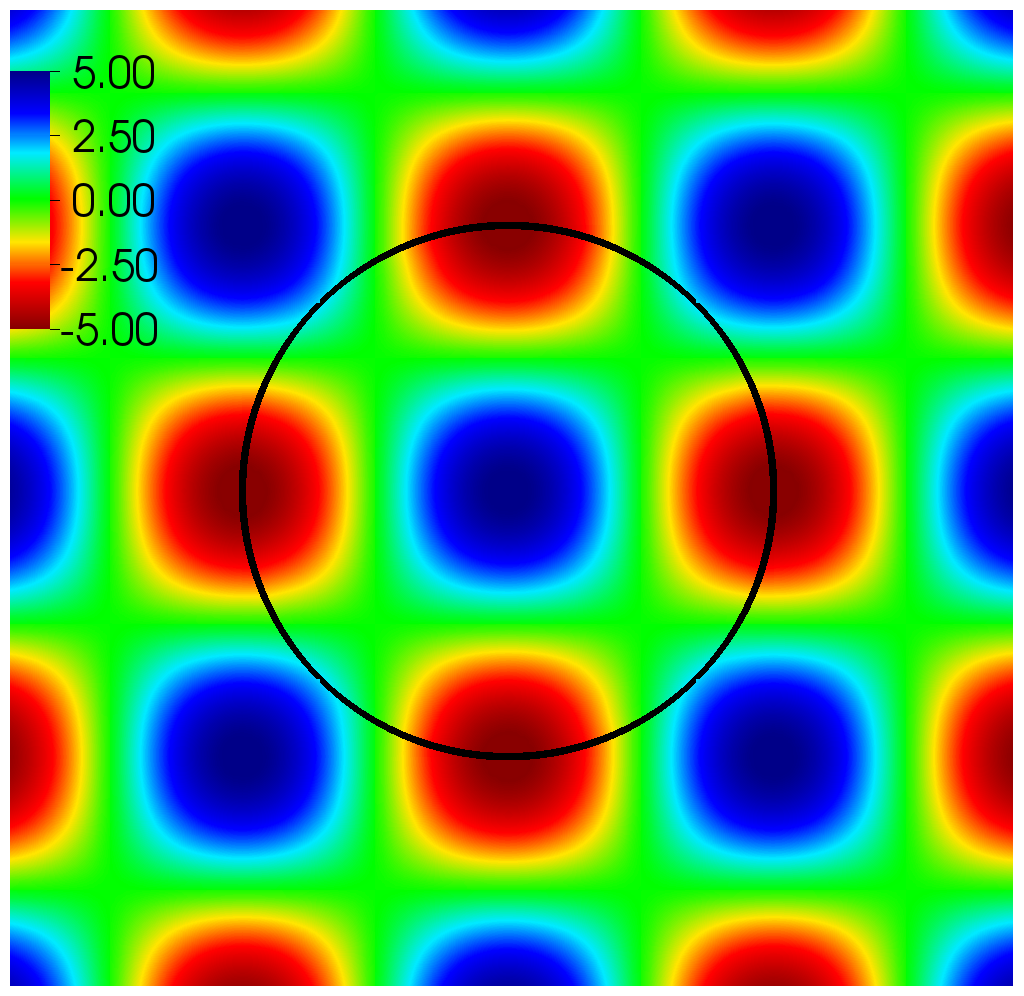} 
		\label{fig:TG:vort_peskin}
	}
	\subfigure[]{
		\includegraphics[width=0.24\textheight]{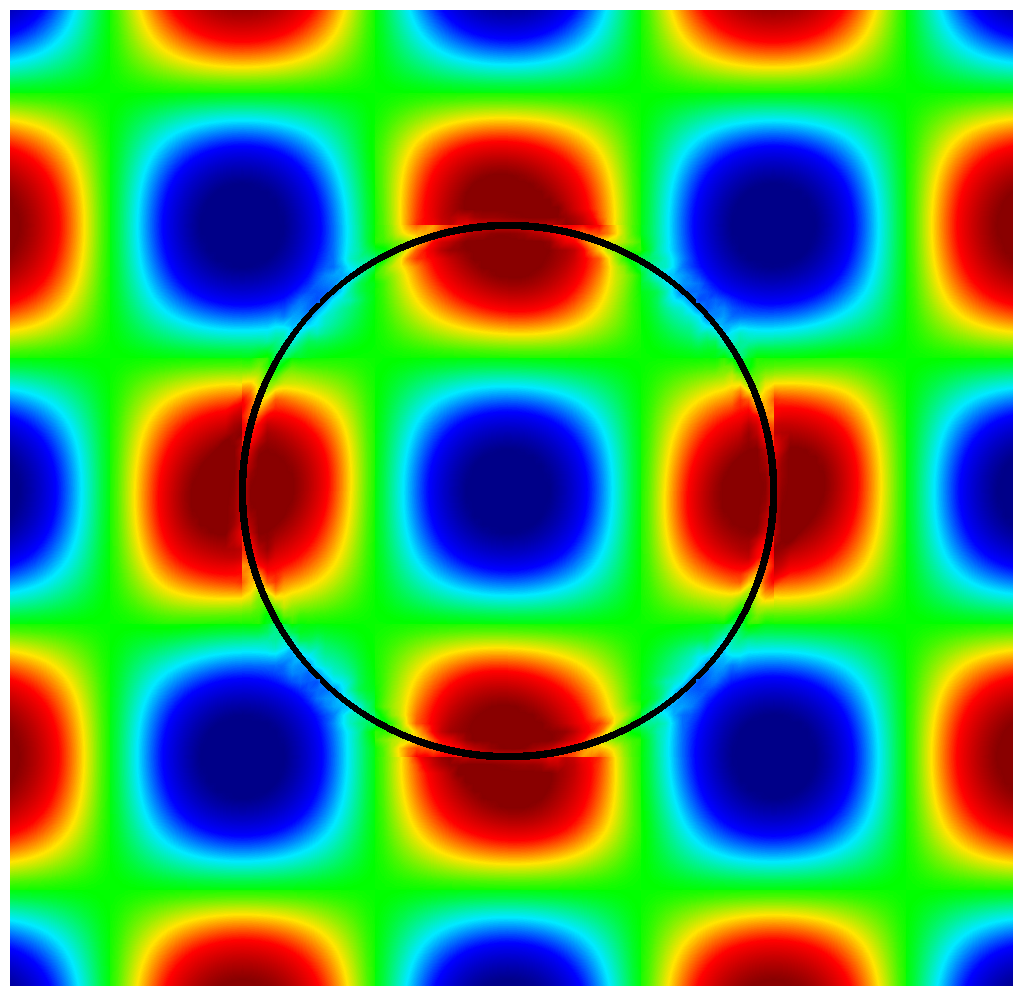} 
		\label{fig:TG:vort_ncvs_both}
	}
	\subfigure[]{
    \includegraphics[width=0.24\textheight]{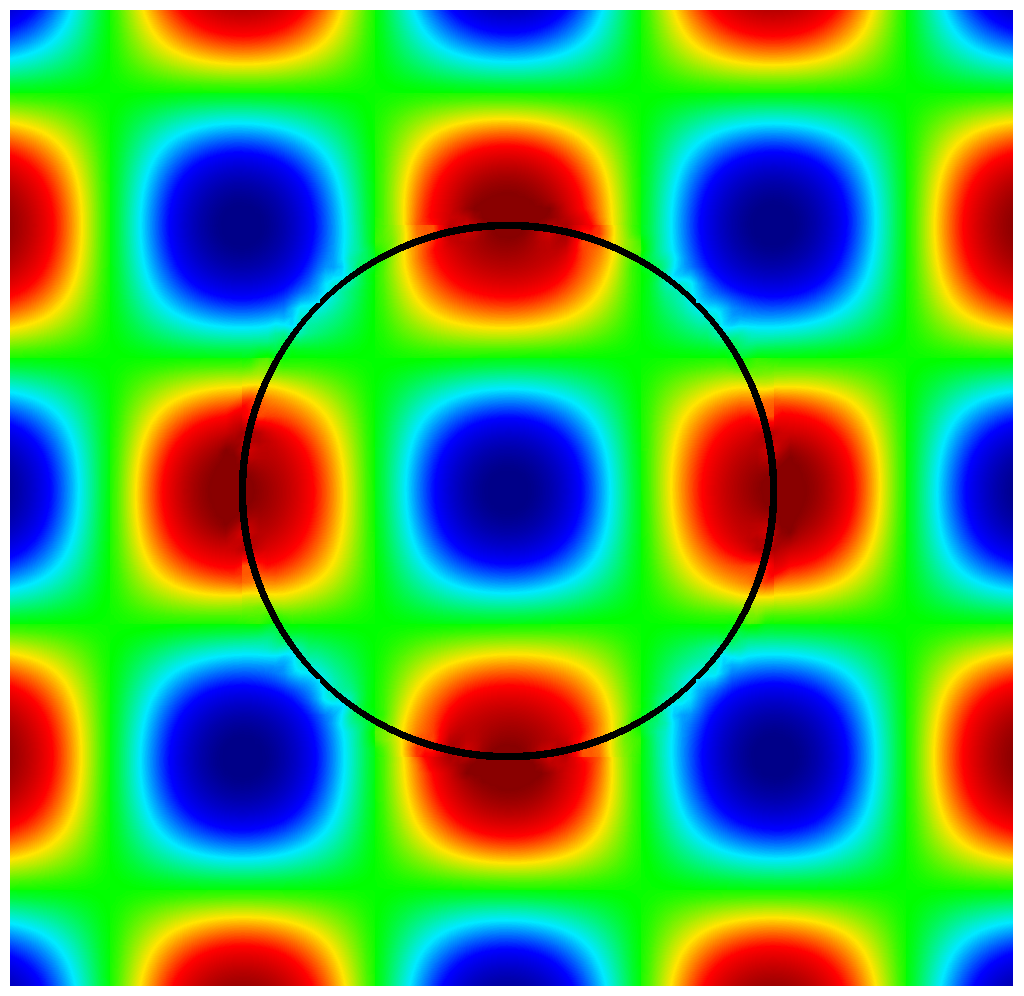} 
    \label{fig:TG:vort_ncvs_one}
}
	\subfigure[]{
    \includegraphics[width=0.24\textheight]{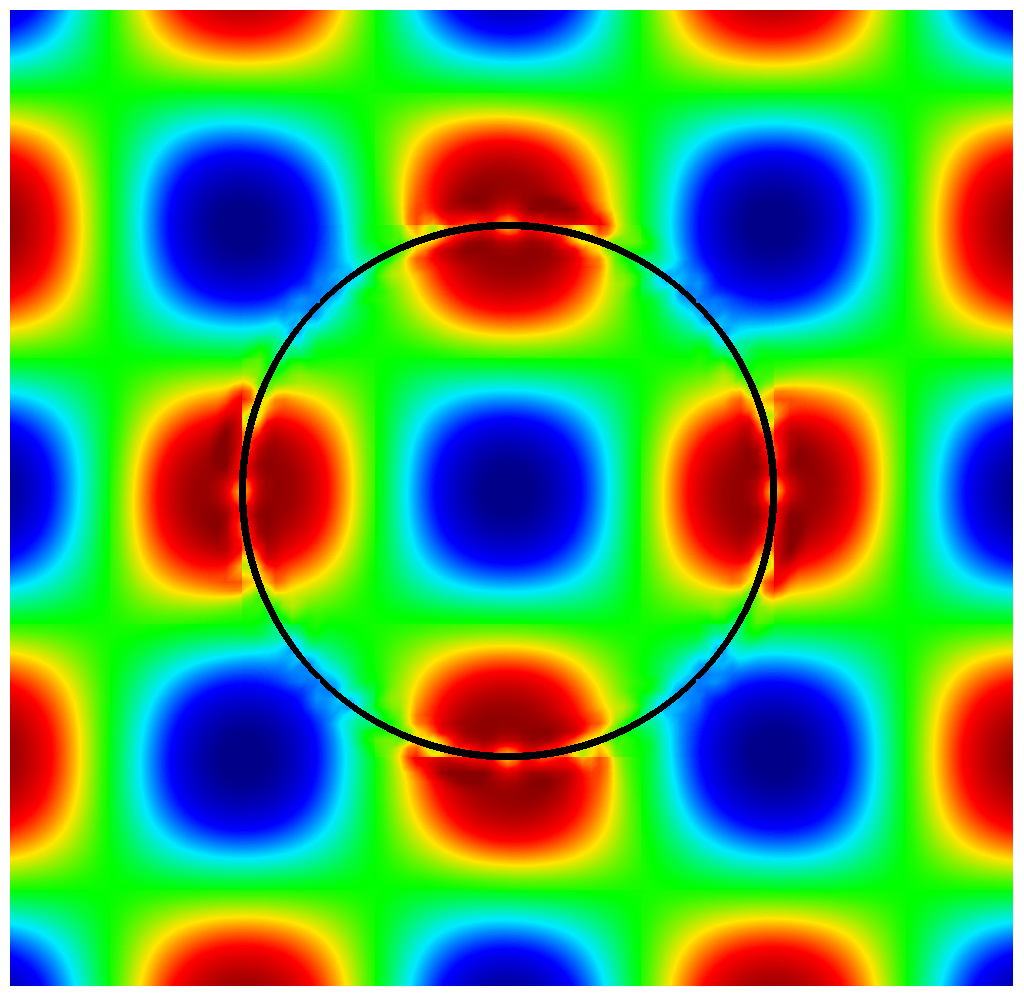} 
    \label{fig:TG:vort_cvs_both}
}
	\subfigure[]{
	\includegraphics[width=0.24\textheight]{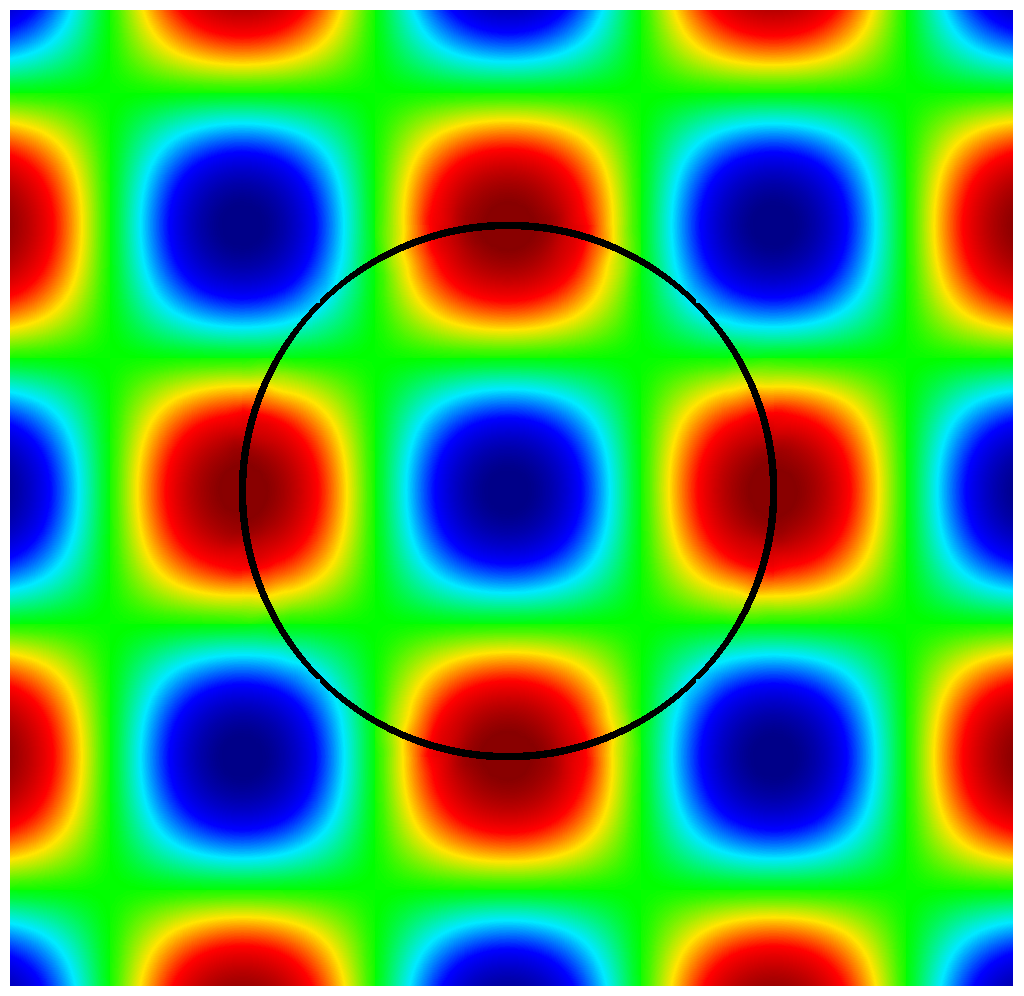} 
	\label{fig:TG:vort_Ipsi_Sncvs}
}

	\caption{Vorticity plots of the Taylor-Green vortex problem for~\subref{fig:TG:vort_peskin} interior and exterior IB forcing with Peskin's four-point delta function;~\subref{fig:TG:vort_ncvs_both}  interior and exterior IB forcing with $\V{\Psi}^{m}_{\text{NCVS}}$ weights in $\cJ$ and $\cS$;~\subref{fig:TG:vort_ncvs_one} exterior IB forcing with $\V{\Psi}^{m}_{\text{NCVS}}$ weights in $\cJ$ and $\cS$;~\subref{fig:TG:vort_cvs_both} interior and exterior IB forcing with $\V{\Psi}^{m}_{\text{CVS}}$ weights in $\cJ$ and $\cS$; and~\subref{fig:TG:vort_Ipsi_Sncvs} interior and exterior IB forcing with $\V{\Psi}$ in $\cJ$ and $\V{\Psi}^{m}_{\text{NCVS}}$ in $\cS$.}
	\label{fig:TGvort}
\end{figure}

\subsection{Taylor-Green vortex flow}
In this test we follow Li et al.~\cite{Li15} and simulate the problem of Taylor-Green vortex flow with an embedded interface upon which exact velocity boundary conditions are imposed. Li et al. tested the accuracy of two-sided IB kernels using the MLS technique of Vanella and Balaras with direct forcing IB method. Here, we test the accuracy of one-sided MLS kernels for this problem.  The analytical solution to the Taylor-Green vortex flow is  

\begin{align}
u(x, y, t) &= -\cos\left(\pi x\right) \sin(\pi y) e^{-2\pi^2t/Re}, \label{eqn:decay_vortex_u} \\
v(x, y, t) &= \sin\left(\pi x\right) \cos(\pi y) e^{-2\pi^2t/Re}, \label{eqn:decay_vortex_v}\\
p(x, y, t) &= -\frac{1}{4}\left(\cos(2\pi x) + \cos(2\pi y)\right) e^{-4\pi^2t/Re}.
\end{align}

The numerical simulation of the decaying vortex is performed with $\Omega = [-2,2]\times[-2, 2]$, and a CFL number of $C= 0.05 $ is chosen. The exact solution at $t=0$ acts as the initial condition. A circular cylinder is chosen to represent the IB surface where a velocity based on the exact solution is imposed. The cylinder is located with its center at the origin of the domain, and its radius is taken to be 1. The Reynolds number of the flow is taken to be $Re = 100$, and the fluid density and viscosity are taken as,  $\rho = 1$ and $\mu = 10^{-2}$, respectively. 

Because the Taylor-Green vortex problem is smooth, formulations using Peskin's second-order accurate delta functions  provide high accuracy. This test also serves as a reference to compare the maximum order of accuracy for the proposed MLS kernels. In particular, more general models will typically include stress discontinuities along the fluid-structure interface, and in such conditions, all of the diffuse-interface methods, including the direct forcing IB method,  will be only first-order accurate. Fig.~\ref{fig:TGError} compares the numerical simulation at $t=1$ with the corresponding exact solution using $L^1$, $L^2$, and $L^{\infty}$ error norms. The $u$ velocity error norms are shown with solid lines, and the pressure error norms are shown with dashed lines.  As observed in Fig.~\ref{fig:TG:error1cIB},  formulations using Peskin's four-point delta function in both $\cJ$ and $\cS$ yields second-order accuracy for the velocity and pressure error norms. Similar results were obtained by Li et al.~\cite{Li15}, who employed two-sided MLS technique of Vanella and Balaras with the five-point cubic spline function weight.

To evaluate the order of accuracy of the proposed MLS kernels, we consider two methods:  
\begin{enumerate}
\item applying interior and exterior IB forcing to interior $\Omega_{b}^{-}$ and exterior $\Omega_{b}^{+}$ regions, respectively; and
\item applying exterior IB forcing to exterior $\Omega_{b}^{+}$ region only. 
\end{enumerate}
In these tests, Peskin's four-point delta function is chosen as the underlying delta-kernel for generating MLS weights. We first consider same weights $\V{\Psi}^{m}_{\text{NCVS}}$ in the velocity interpolation and force spreading operators, so that $\cS = \cJ^{*}$.  In Figs.~\ref{fig:TG:error1} and~\ref{fig:TG:error2},  the rate of convergence of $L^1$, $L^2$, and $L^{\infty}$ error norms with respect to the grid cell size for both methods of imposing IB forces are presented. The convergence rates are found to be approximately 1 for velocity and pressure for all error norms. This suggests that the two methods impose the same (velocity) boundary conditions on the two sides of the interface, and hence the order of accuracy of the solution is insensitive to the IB force location.  The first-order convergence rates that occur when using modified weights $\V{\Psi}^{m}_{\text{NCVS}}$ in $\cS$ and $\cJ^{*}$ is also expected from the results of Sec.~  \ref{subsec_positivity}. For completeness, we also report in Fig.~\ref{fig:TG:error1CVS}, the convergence rates obtained by using $\V{\Psi}^{m}_{\text{CVS}}$ to construct $\cJ$ and $\cS$. IB forcing on both interior and exterior regions is considered (method 1). First-order accuracy is observed for this case as well.

Next, we compare the order of accuracy by employing non-shifted $\V{\Psi}$ weights in $\cJ$ and shifted $\V{\Psi}^{m}_{\text{NCVS}}$ weights in $\cS$. In this case $\cS \ne \cJ^{*}$. IB forcing on both interior and exterior regions is considered (method 1).  Fig.~\ref{fig:TG:errorIpsiSncvs} shows super-linear accuracy in all velocity error norms and linear accuracy in the pressure error norm. This confirms that one-sided MLS weights satisfy the linear polynomial reproduction constraint, which improves the overall order of accuracy of the scheme. We note that the non-adjointness of $\cS$ and  $\cJ$ operators did not affect the stability of this scheme. However, using $\V{\Psi}$ in $\cS$ resulted in flow instabilities for this case as well (data not shown). 

Fig.~\ref{fig:TGvort} shows the vorticity plots for all of the MLS cases considered in this section. All cases show similar trends, except for the CVS strategy employed in Fig.~\ref{fig:TG:vort_cvs_both}, in which \REVIEW{large} spurious oscillations in the vorticity field are observed near the interface. This again highlights the fact that the proposed CVS approach to mollifying the MLS weights leads to unphysical flow fields. \REVIEW{Mild oscillations are also observed in Figs.~\ref{fig:TG:vort_ncvs_both} and \ref{fig:TG:vort_ncvs_one}, which can be attributed to the first-order accuracy of $\V{\Psi}^{m}_{\text{NCVS}}$ weights used in the interpolation operator $\cJ$.}

In what follows, we report results by employing $\V{\Psi}^{m}_{\text{NCVS}}$ in both velocity interpolation and force spreading operators so that $\cS = \cJ^{*}$ in these cases. Although this coupling scheme reduces the order of accuracy as indicated by the results of this section, it will help bring-forth the worst-case scenario for one-sided IB/MLS simulations, if there is any. 

\REVIEW{
\subsection{Stokes' first problem}
To validate that the current approach can capture time-dependent viscous boundary layers accurately,  we consider the Stokes' first problem. The problem setup involves an infinitely long plate that is impulsively started with a velocity $ U_{p} $ parallel to its axis oriented along the horizontal $x$-direction. Exact analytical solutions for the time evolution of the velocity profile and viscous drag coefficient are available for this problem. The transient velocity profile and drag coefficient are given by 
\begin{align}
	&u(y,t) =U_p \left[ \textrm{erfc}\left( \frac{y}{2\sqrt{\nu t}}\right) \right], \\
	&C_D =\frac{2}{\sqrt{\pi t Re}}.
\end{align} 
The Reynolds number is defined as $Re  =U_p L/\nu $, in which $L$ is the characteristic length and $ \nu = \mu/\rho $ is the kinematic viscosity of the fluid. 

\begin{figure}
	\centering
	\subfigure[]{
		\includegraphics[width=0.47\textwidth]{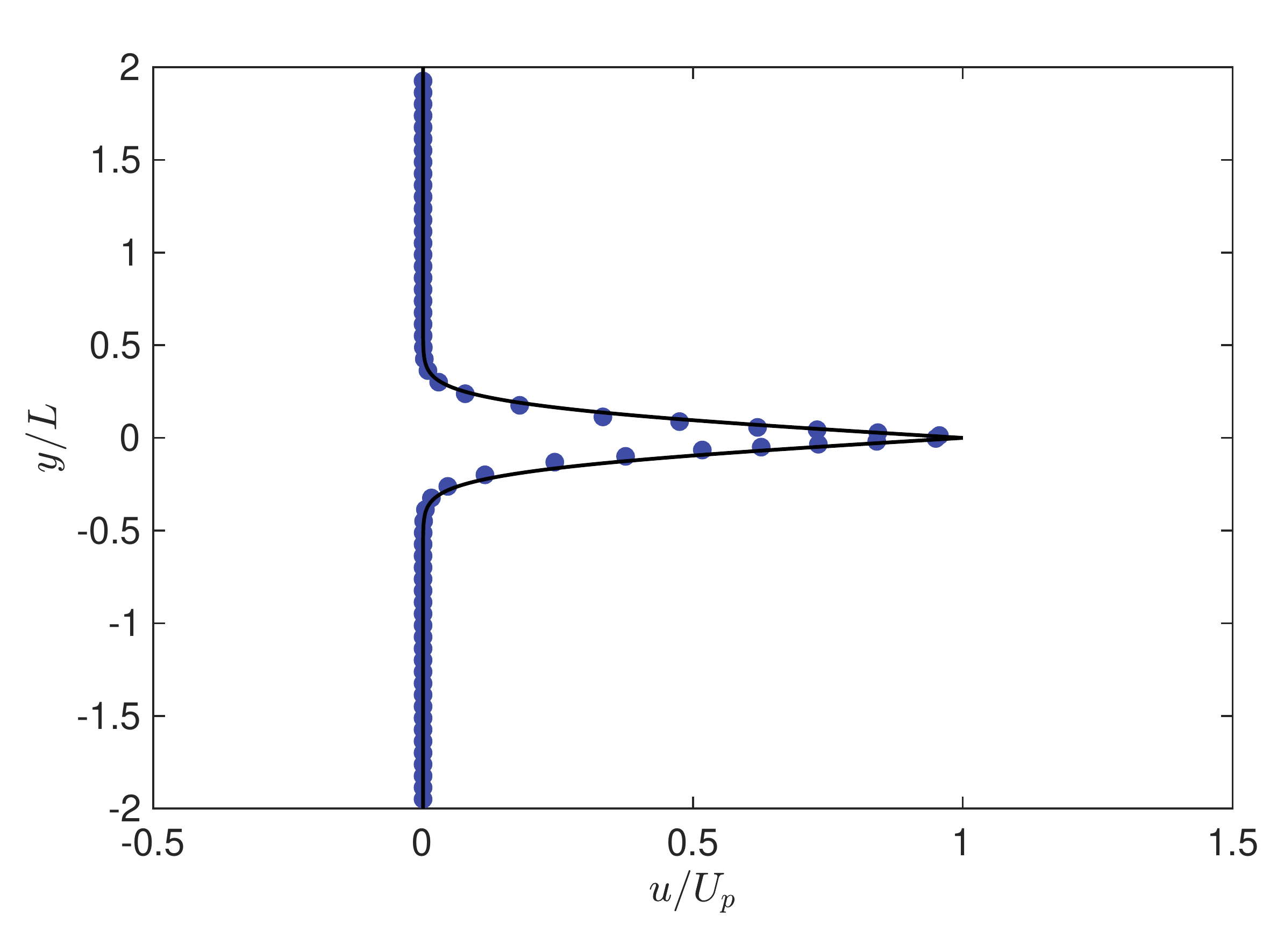} 
		\label{fig:Stokes:profile}
	}
	\subfigure[]{
		\includegraphics[width=0.47\textwidth]{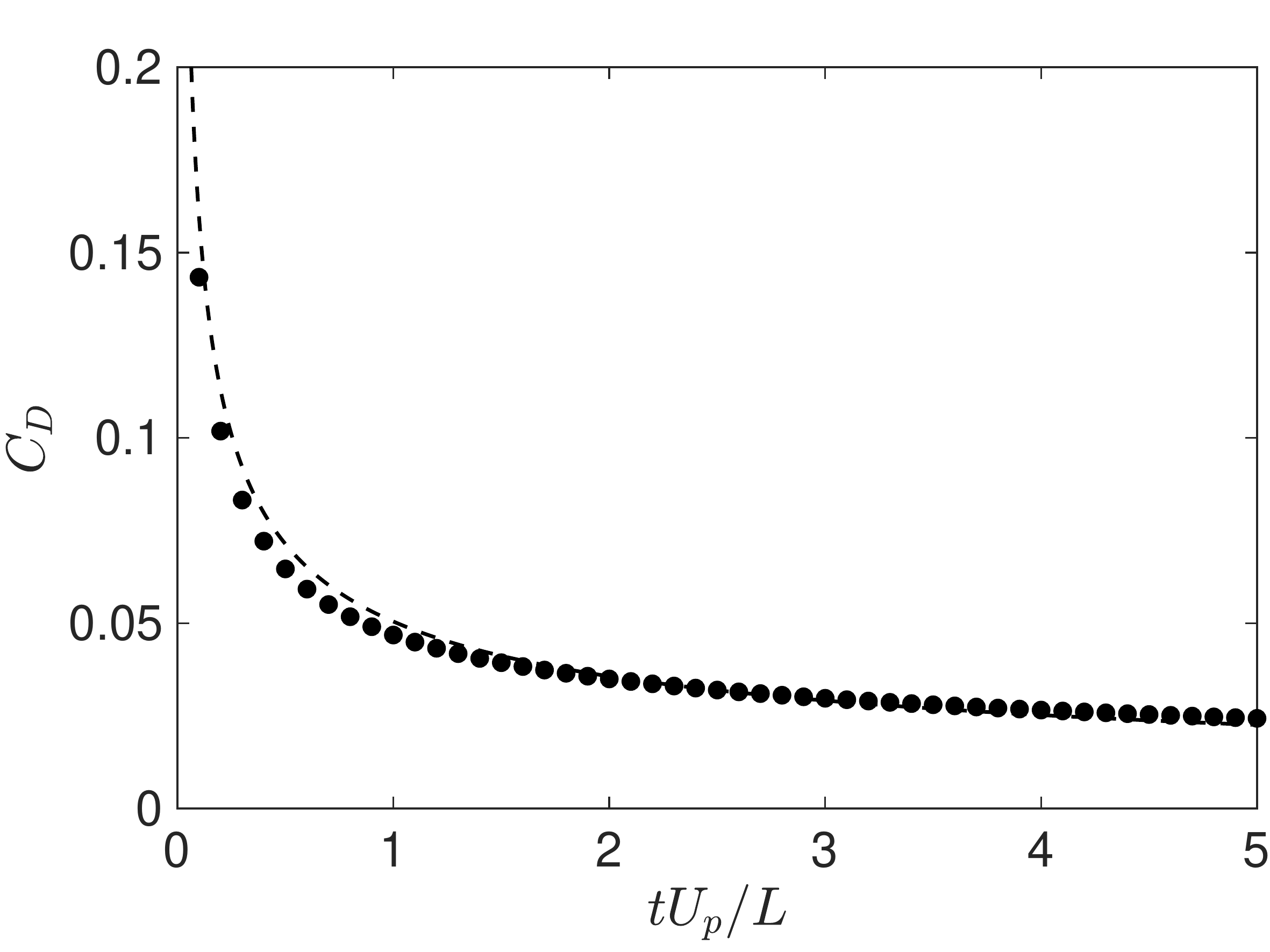} 
		\label{fig:Stokes:Cd}
	}
	\caption{\REVIEW{Comparison of analytical (solid line) and numerical (solid circle) solutions of the Stokes' first problem at $Re=500$.~\subref{fig:Stokes:profile} Vertical velocity profile $u(y,t)$ at non-dimensional time $ T=tU_p/L =5$.~\subref{fig:Stokes:Cd} Time evolution of viscous drag coefficient.}}
	\label{fig:Stokes}
\end{figure}

We carry out the simulation at $ Re=500 $ by setting $ U_p = 1 $ and $ L=1 $. The size of the computation domain is taken to be $ 4L \times 4L $ with periodic boundary conditions along both $ x $ and $ y $ directions. A plate stretching across the domain along the $ x $-axis is placed in the middle of the periodic box. The domain discretization is chosen such that the mesh spacing around the plate is $ \Delta x = \Delta y = 0.002L $.  Note that this problem setup corresponds to the open geometry configuration as shown in Fig.~\ref{subfig:open}, and the top and bottom sides of the plate interact separately with the corresponding fluid domains.   The velocity profile over the plate (on both sides) at non-dimensional time $ T=tU_p/L =5$ is plotted in Fig.~\ref{fig:Stokes:profile}; the numerical solution agrees quite well with the analytical solution. Next, we compare the evolution of the drag coefficient $ C_D = F_x/\half \rho U^2_p A$ ($F_x$ is the viscous force in the $x$-direction, and $ A $ is the area of the plate assuming a unit depth), with the analytical solution which is presented in Fig.~\ref{fig:Stokes:Cd}. The agreement between the analytical and simulated drag coefficient profiles is also excellent. We also simulated this case using the standard Peskin kernel. The simulation results were identical to the ones shown in Fig.~\ref{fig:Stokes}, and are not shown here for brevity. 


\begin{figure}
	\centering
	\includegraphics[width=0.55\textwidth]{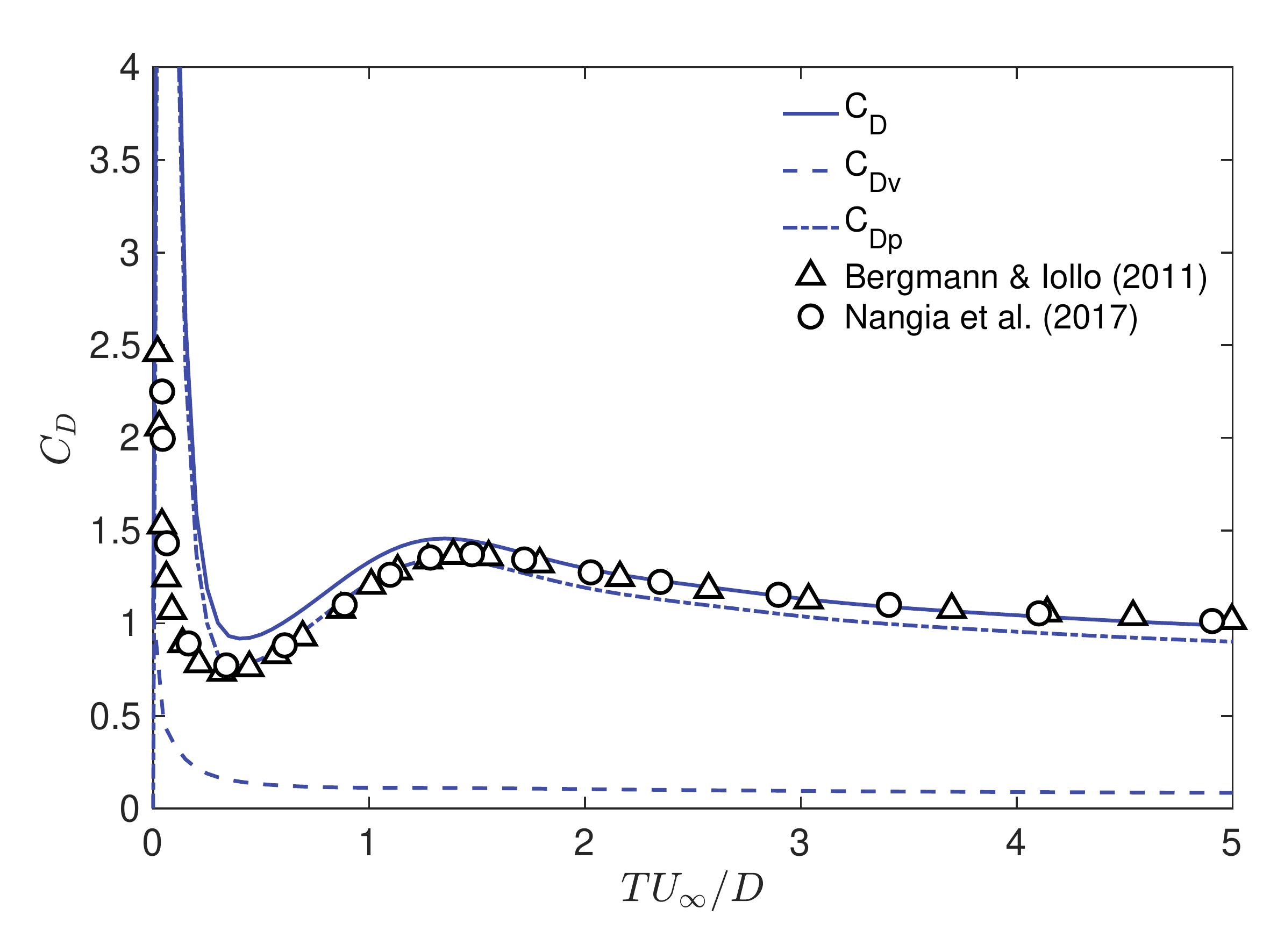} 
	\caption{\REVIEW{Comparison of time variation of drag coefficient $ C_D $ of cylinder in an impulsively started flow with Bergmann and Iollo~\cite{Bergmann11} and Nangia et al.~\cite{Nangia17}. Also plotted are the time evolution of skin friction and pressure coefficients, $ C_{Dv} $ and  $C_{Dp}$, respectively.  }}
	\label{fig:IC}
\end{figure}

\subsection{Impulsive flow over a cylinder}
Next, we consider an impulsively started flow over a 2D circular cylinder of diameter $D$, which also serves to quantify the transient boundary layer development. The simulation is performed in the computation domain $\Omega = [-4D,12D] \times [-16D, 16D]$. The diameter of the cylinder is $ D=1 $, and it is placed with its center at $(x,y) = (0,0)$. The Reynolds number is $Re=D U_{\infty} /\nu$, and is set to 500, the flow velocity in the horizontal $x$-direction is $ U_{\infty}= 1 $, and the kinematic viscosity of the fluid is $ \nu=1/Re = 0.002$. The mesh resolution on the cylinder surface is approximately $ 0.01D $. We validate the results of our simulation by comparing the evolution of the drag coefficient with the numerical results of Bergmann and Iollo~\cite{Bergmann11} and Nangia et al.~\cite{Nangia17}. The drag coefficient along the flow direction is defined as $ C_D = F_x/ \half\rho D U_{\infty}^2$. The comparison of $ C_D $ with data from the literature plotted in Fig.~\ref{fig:IC} shows a good agreement. In the same plot, we also present evolution of the viscous ($ C_{Dv} $) and pressure contribution ($C_{Dp}$) to the total drag. These coefficients are defined analogous to $C_D$ by using horizontal component of viscous and pressure forces, respectively.  Note that both prior studies~\cite{Bergmann11,Nangia17} compute the net hydrodynamic force on the body in an extrinsic manner, i.e.  indirectly through momentum-conservation principle, whereas we compute it directly in an intrinsic manner by integrating the hydrodynamic stress tensor on the surface of the body. The difference in the initial transient drag profile in Fig.~\ref{fig:IC} is attributed to the manner in which hydrodynamic forces are evaluated in the current and the prior works.   

}

\subsection{Impulsively started plate}
\label{sec_impulsiveplate}
We consider the flow over an impulsively started, infinitesimally-thin, two-dimensional plate to validate the present numerical method for moving and finite-sized non-closed geometry. The flow resulting from the plate motion is characterized by two counter-rotating vortices in the wake which are fed by the shear layer at the tip of the two ends of the plate. The size of the wake vortices gradually increases before eventual saturation. Our numerical method is validated by capturing the dynamics of the time evolution of the wake vortices.  We perform two-dimensional simulation at Reynolds numbers of 126 and 1000, in which the Reynolds number is based on the height $H_p$ and speed $U_p$ of the plate. The $Re$ is chosen to enable comparison with experimental and numerical results reported in the literature~\cite{tane71,koum96}. A computational domain of size $10H_p\times 10H_p$ with no-slip boundary condition on the domain boundaries in the lateral direction and periodic boundary condition along the axial (corresponding to the direction of plate motion) is used. A uniform mesh resolution of $0.01H_p$ is chosen for both $Re$ cases. This resolution is consistent with the that of Mittal and coworkers~\cite{mitt08}. At $t=0$ the plate is placed at the center of the computational domain and impulsively started with a velocity of $U_p$ in the axial direction. 

\begin{figure}
	\centering
	\subfigure[$Re = 126$]{
	   \includegraphics[scale=0.24]{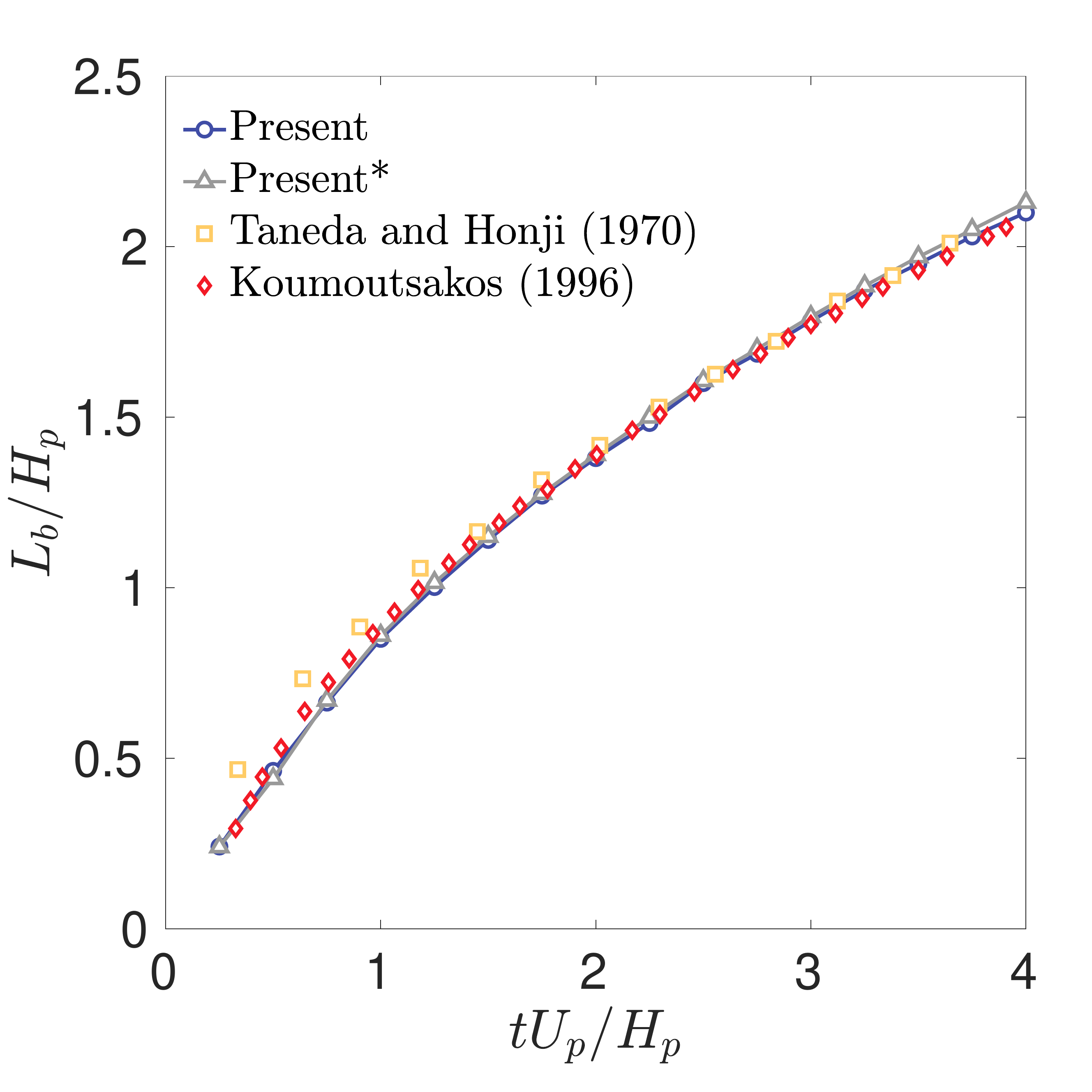} 
	   \label{fig:plate:bsize_evol-Re126}
	}
	\subfigure[$Re = 1000$]{
	    \includegraphics[scale=0.24]{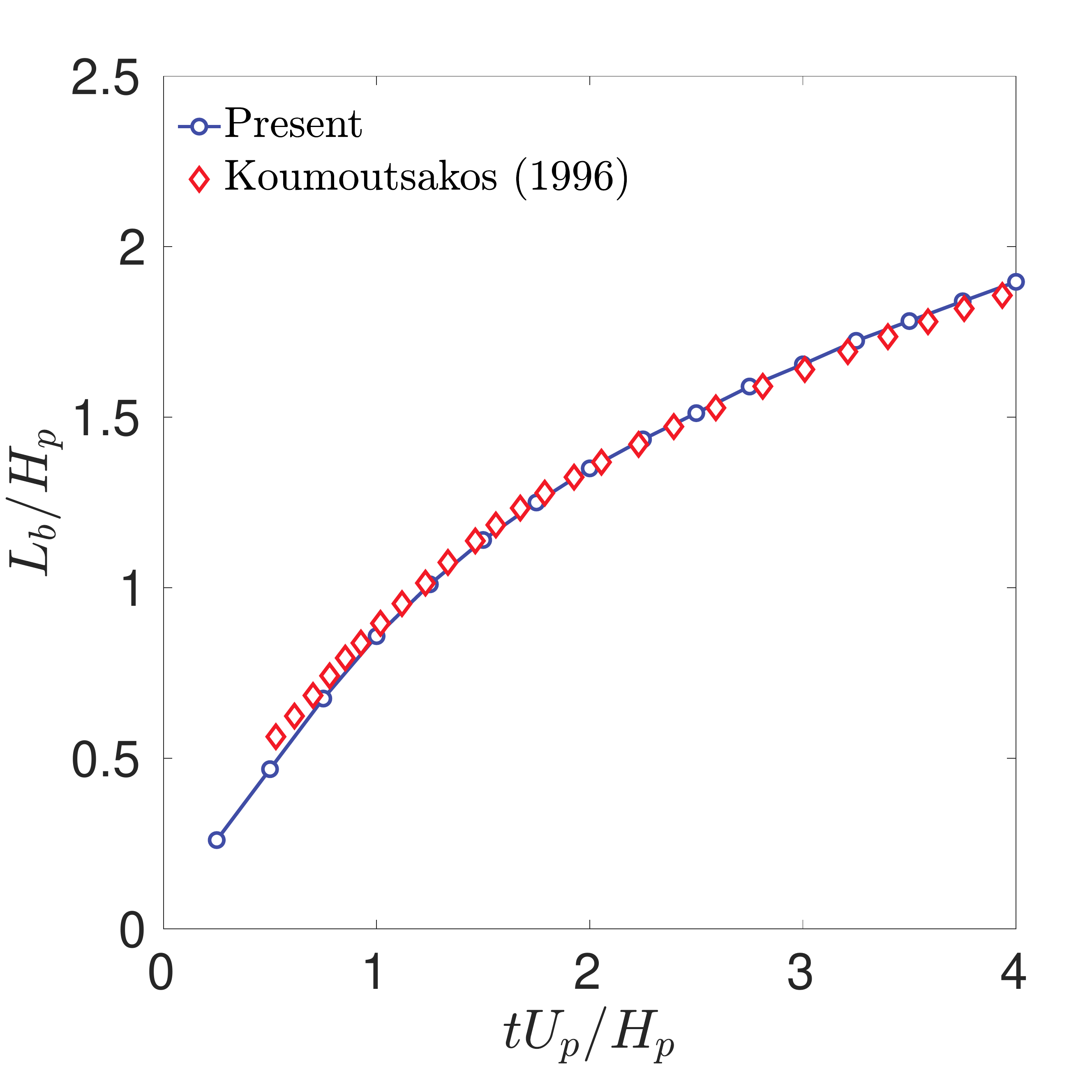} 
  	    \label{fig:plate:bsize_evol-Re1000}
	}
\caption{ The wake bubble size plotted against  time in dimensionless form at~\subref{fig:plate:bsize_evol-Re126} $Re=126$, and~\subref{fig:plate:bsize_evol-Re1000} $Re=1000$. Results of the present numerical method are compared against reported data in the literature~\cite{tane71,koum96}. In~\subref{fig:plate:bsize_evol-Re126},  the  curve (---o---) corresponds to using $\V{\Psi}^{m}_{\text{NCVS}}$ weights in both $\cJ$ and $\cS$, whereas the curve (---$\triangle$---) corresponds to using $\V{\Psi}$ weights in $\cJ$ and $\V{\Psi}^{m}_{\text{NCVS}}$ weights in $\cS$.}
\label{fig:plate:bsize_evol}
\end{figure}

\begin{figure}
	\centering
	\subfigure[$Re = 126$]{
	   \includegraphics[width=0.6\textwidth]{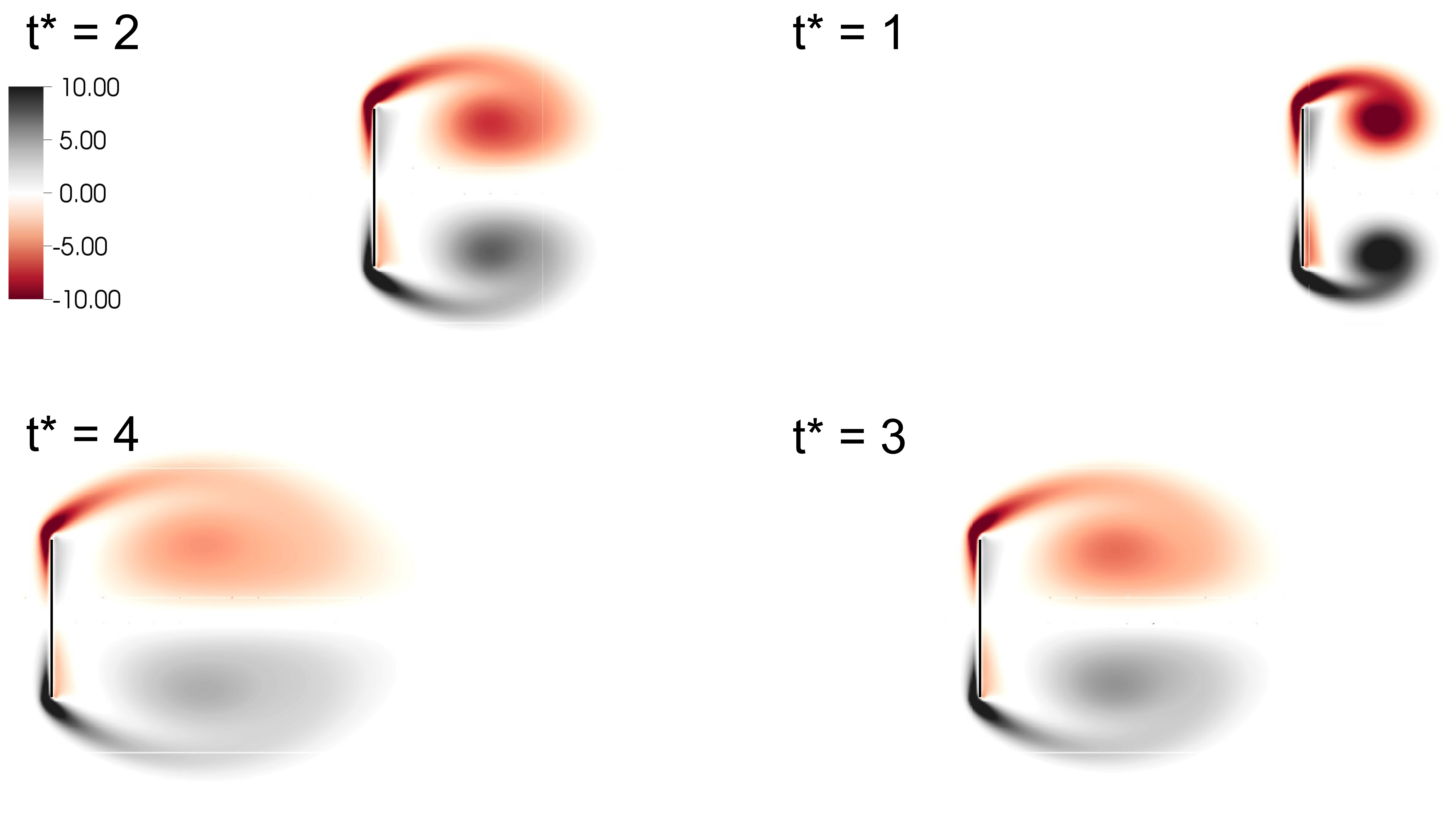} 
	   \label{fig:plate:vort-Re126}
	}
	\subfigure[$Re = 1000$]{
	    \includegraphics[width=0.6\textwidth]{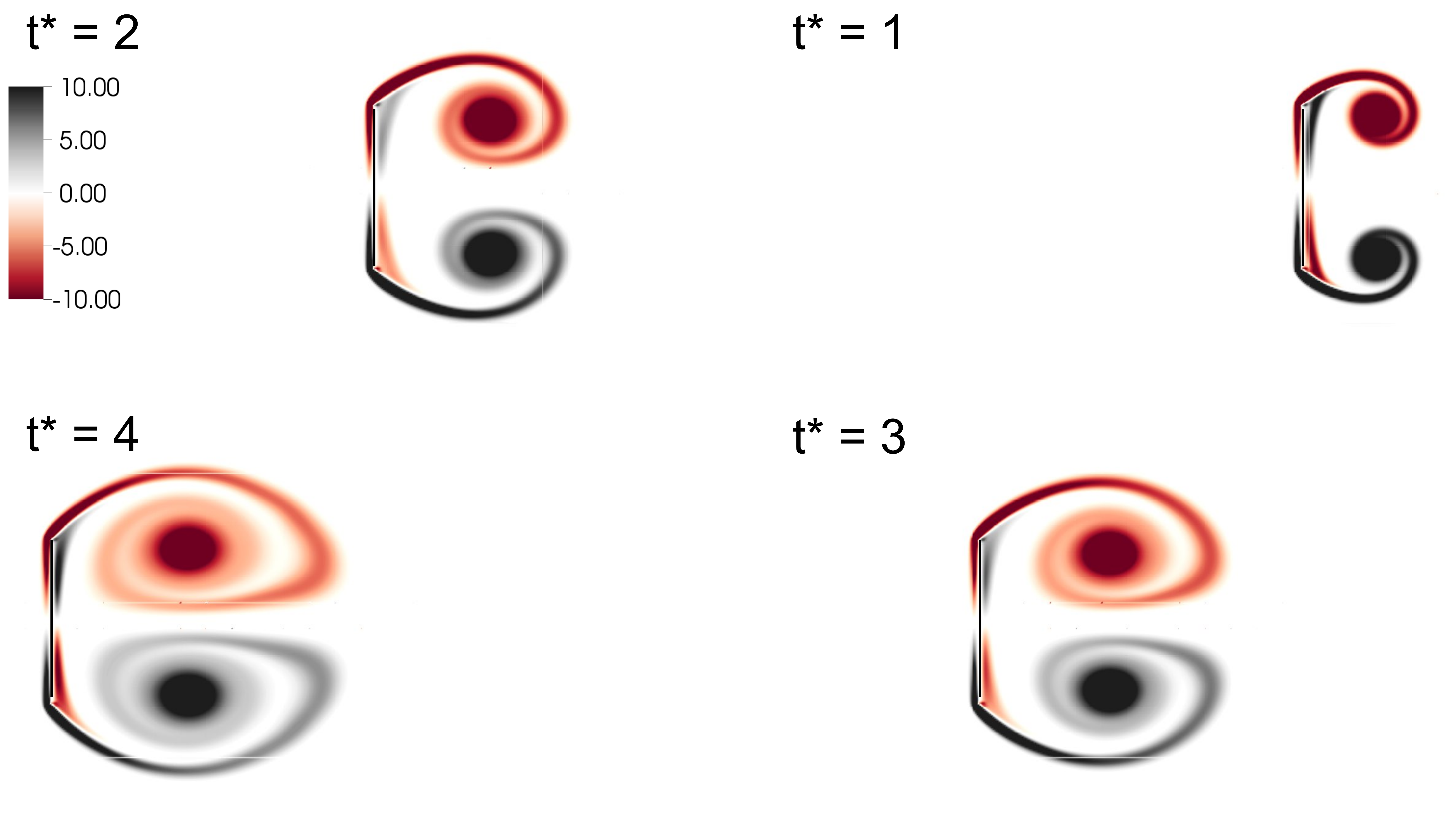} 
  	    \label{fig:plate:vort-Re1000}
	}
\caption{ Evolution of the separation bubble in the wake of an impulsively started plate at different instants of time for~\subref{fig:plate:vort-Re126} $Re=126$; and~\subref{fig:plate:vort-Re1000} $Re=1000$.}
\label{fig:plate:vort}
\end{figure}

The height and  velocity of the plate  are used as  characteristic length and velocity scale, respectively, to define the nondimensional wake bubble size $L_{b}^{*} = L_{b}/H_p$ and the nondimensional time $ t^{*} = tU_p/H_p$. We measure the evolution of the size of counter-rotating vortices in the wake bubble. The length of the wake bubble $L_b$ is the normal distance from the plate to the end of the region of reverse flow behind the plate.  Fig.~\ref{fig:plate:bsize_evol} shows the temporal evolution of the nondimensional wake bubble size, and we compare our results against the reported data in the literature. We find that our results are in good agreement with the experimental results of Taneda and Honji~\cite{tane71} at  $Re=126$, and with numerical results of Koumoutsakous and Sheils~\cite{koum96} at $Re=126$ and $Re=1000$. In Fig.~\ref{fig:plate:bsize_evol-Re126}, the results obtained by using non-shifted $\V{\Psi}$ weights in $\cJ$ and shifted $\V{\Psi}^{m}_{\text{NCVS}}$ weights in $\cS$ are compared against those obtained by using shifted $\V{\Psi}^{m}_{\text{NCVS}}$ weights in both $\cS$ and $\cJ$ operators. The results are found to be virtually indistinguishable. The progression of vorticity development behind the plate is visualized in Fig.~\ref{fig:plate:vort} for the two $Re$ cases.

\REVIEW{
\subsection{Flow around an oscillating cylinder}
In order to further investigate the robustness of the current numerical method in handling moving geometries, we consider the flow around an oscillating cylinder, which involves a time-varying motion of the body. The problem setup involves a circular cylinder oscillating about a mean position in a quiescent fluid. Assuming that the cylinder is oscillating along the $ x$-axis, the velocity of oscillation of the cylinder is 
\begin{equation}
	U = U_{m}\cos\left(\frac{2\pi}{T}t\right),
	\label{eqn:ocyl}
\end{equation} 
in which $ U_m $ and $ T $  are the peak velocity and time period of oscillation, respectively. The flow generated by the oscillating cylinder is characterized by two dimensionless parameters,  namely, the Keulegan-Carpenter number $ KC=U_m T/D $ and the Reynolds number $ Re=U_m D/\nu $. Here, $ D $ is the diameter of the cylinder. The case has been well studied in the literature, for example~\cite{Shen09,Bhalla13}. We compare our computations against the results of Shen et al.~\cite{Shen09} and Bhalla et al.~\cite{Bhalla13} by setting $ D=1 $ and $ U_m=1 $, $KC=5 $, and $ Re=100 $. The computational setup consists of a cylinder placed at the center of a domain at $ t=0 $, and a square domain that spans a length of $ 32D $ along both $ x $ and $ y $ directions. The mesh spacing on the cylinder surface and along its path of oscillation is chosen to be $ 0.008D $, which also matches the spacing employed in~\cite{Bhalla13}. 

\begin{figure}
	\centering
	\subfigure[]{
		\includegraphics[width=0.45\textwidth]{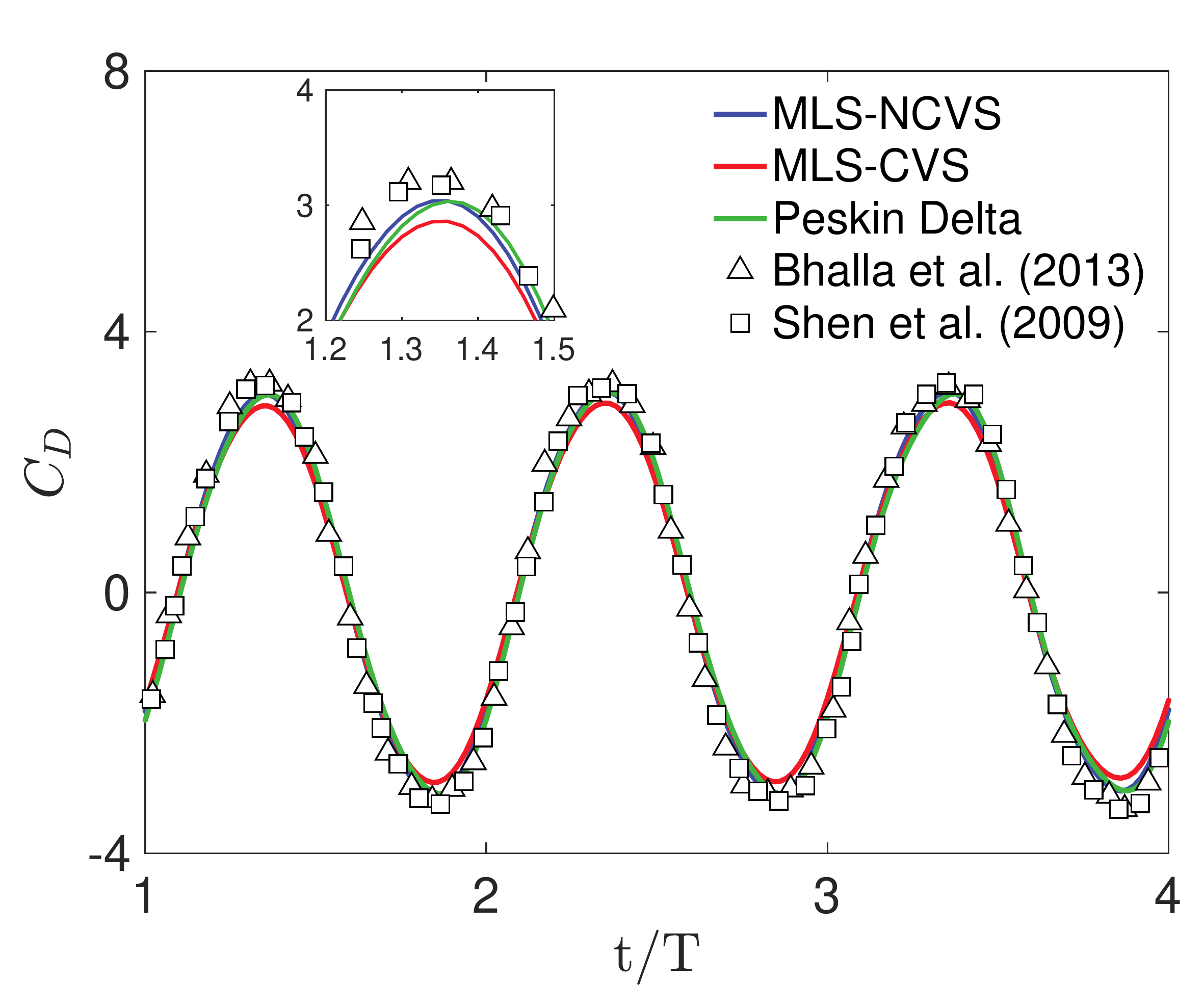} 
		\label{fig:OC:Cd}
	}
	\subfigure[]{
		\includegraphics[width=0.45\textwidth]{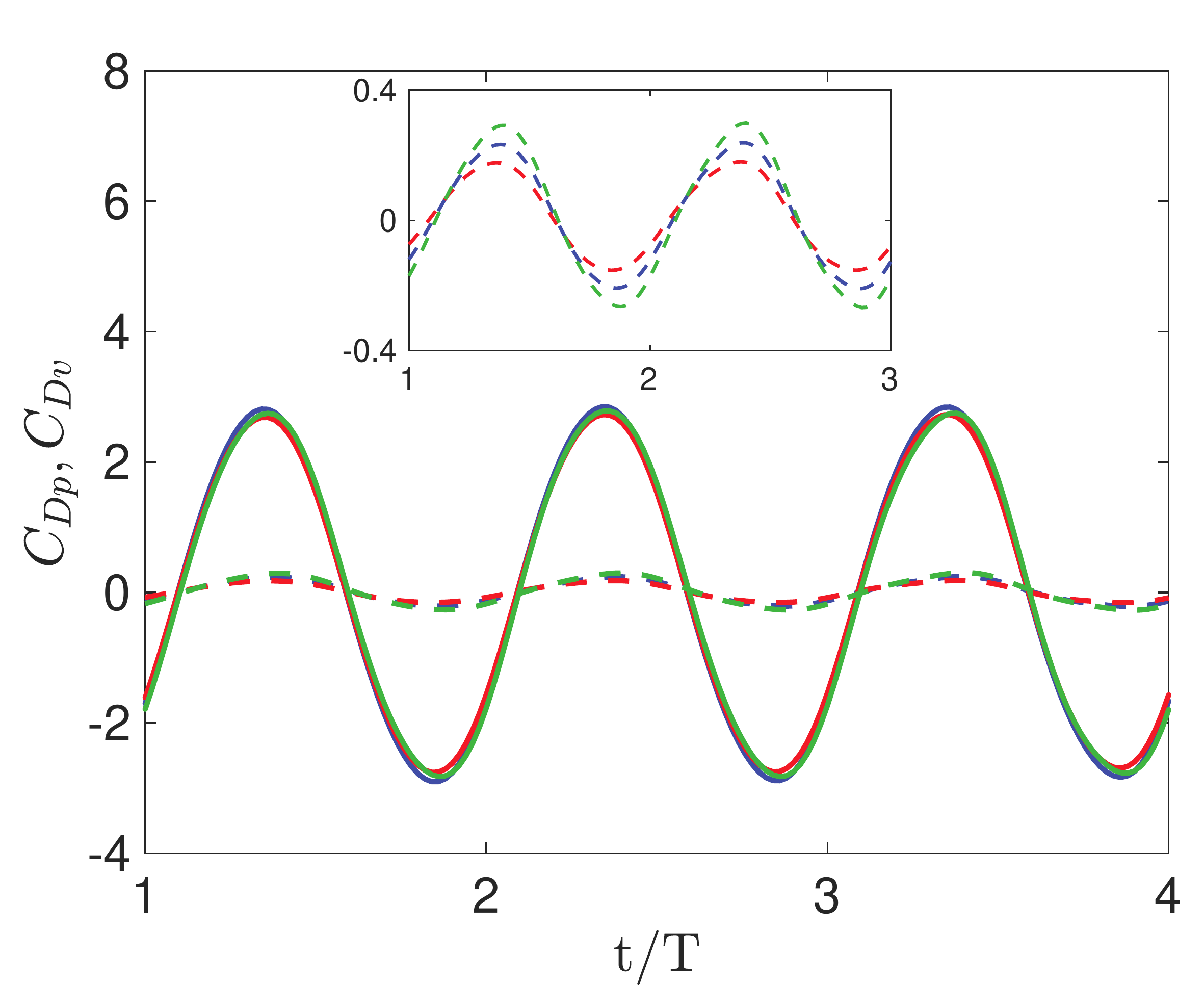} 
		\label{fig:OC:Cd-split}
	}
	\caption{\REVIEW{Flow around an oscillating cylinder at $Re = 100$ and $KC = 5$.~\subref{fig:OC:Cd} Comparison of the drag coefficient $C_D$ as a function of time with the literature.~\subref{fig:OC:Cd-split} Comparison of pressure $C_{Dp}$(solid lines) and skin friction $C_{Dv}$ (dashed lines) coefficients using one-sided CVS and NCVS mollification strategies with those obtained using the two-sided four-point Peskin kernel. The color coding of the curves in~\subref{fig:OC:Cd-split} is same as~\subref{fig:OC:Cd}.}}
	\label{fig:OC}
\end{figure}

To elucidate the differences between the CVS and NCVS mollification strategies, the current validation case is performed using both approaches and compared to results obtained using the standard four-point Peskin kernel. The drag coefficient, $ C_D=F_x/ \half \rho U^{2}_{m} D$, computed from the simulation is presented in Fig.~\ref{fig:OC:Cd}, and compared against the literature. There is good agreement between the present results and those reported previously~\cite{Shen09,Bhalla13}. However, closer inspection reveals that near the peaks of the $ C_D $ curve, the two-sided Peskin kernel and the one-sided NCVS kernel show better agreement with the literature data. Differences in the approaches are further investigated by comparing the total skin friction coefficient $C_{Dv}$ and pressure coefficient $C_{Dp}$, against those obtained using the two-sided delta function; see Fig~\ref{fig:OC:Cd-split}. Whereas the differences in the pressure coefficient $C_{Dp}$ are negligible, we find that the CVS approach underestimates the skin friction coefficient to a greater extent than the NCVS approach; see the inset of Fig.~\ref{fig:OC:Cd-split}.  We attribute this difference to the nature of the magnitude and (lesser) decay of the CVS weights away from the evaluation point (cf. Fig.~\ref{fig:MLSProc}), which leads to smaller velocity gradients. 

}

\subsection{Flow past a sphere} \label{sec_spherecase}
For validating three-dimensional flows with one-sided kernels, we consider the canonical case of flow past a sphere.  The availability of detailed local, global and topological flow characteristics via several experimental and direct numerical simulation based studies makes flow over a stationary sphere an excellent candidate for 
validation~\cite{john99,mitt08,tomb00,clif78}. The flow around a sphere exhibits steady axisymmetric, steady non-axisymmetric and unsteady non-axisymmetric behavior depending on the Reynolds number of the flow. For the present study, we choose Reynolds numbers of 100, 300, and 1000, which cover both steady and unsteady regimes.  The computational domain is  $\Omega = [-40D,-40D,-40D]\times[120D, 40D,40D]$ for $ Re=100$  and $Re=300$ cases, and $\Omega = [-25D,-60D,-60D]\times[100D, 65D,65D]$ for the $ Re=1000 $ case. In all of the cases, a sphere of diameter $D$ is placed with its center at $(x,y,z) = (0,0,0)$.  Local mesh refinement is used to enhance the mesh resolution around the sphere. The near-surface mesh resolution in  $ Re=100$ and $Re=300$ is $0.01D$, and for  $Re=1000$ is $0.008D$. A uniform inflow boundary condition and a homogeneous Neumann boundary condition are imposed at the inflow
and outflow boundaries, respectively, and  free-slip boundary conditions are imposed on the vertical and lateral walls of the computational domain.  

As discussed in Sec.~\ref{sec_mls_theory}, the number of interpolation points in an MLS procedure depends upon the choice of the weighting function. Kernels with very narrow width can result in a near-singular Gram matrix when used for one-sided interpolation.  Using flow around a sphere at $Re=100$, we investigate the stability and accuracy of various kernels.  Through this study, we can understand the smallest stable and accurate weighting kernel suitable for the MLS methodology. 

The steady axisymmetric flow at $Re=100$ is characterized by a recirculation bubble in the immediate wake of the sphere. We identify the key geometric characteristics of the recirculation bubble, namely its center and length and compare them against the reported data in the literature. The center of the recirculation $(x_c, y_c)$ is measured with respect to the center of the sphere, and the bubble length $L_b$ is measured from the leeward end of the sphere. In Table~\ref{tab:spherewake} we report the center and length of the wake bubble nondimensionalized by $D$ for all of the kernels, except the smoothed three-point IB kernel. Our tests found that the masking procedure using the smoothed three-point IB kernel results in a near-singular Gram matrix, and leading to numerical divergence. This test suggests that kernels whose width is three grid point or less are not suitable for the one-sided IB/MLS method presented in this work. The wider IB kernels are stable and pose no stability issues. 

\begin{figure}
\centering
\subfigure[Using $\delta_{4}$ weights in MLS]{
	\includegraphics[width=0.45\textwidth]{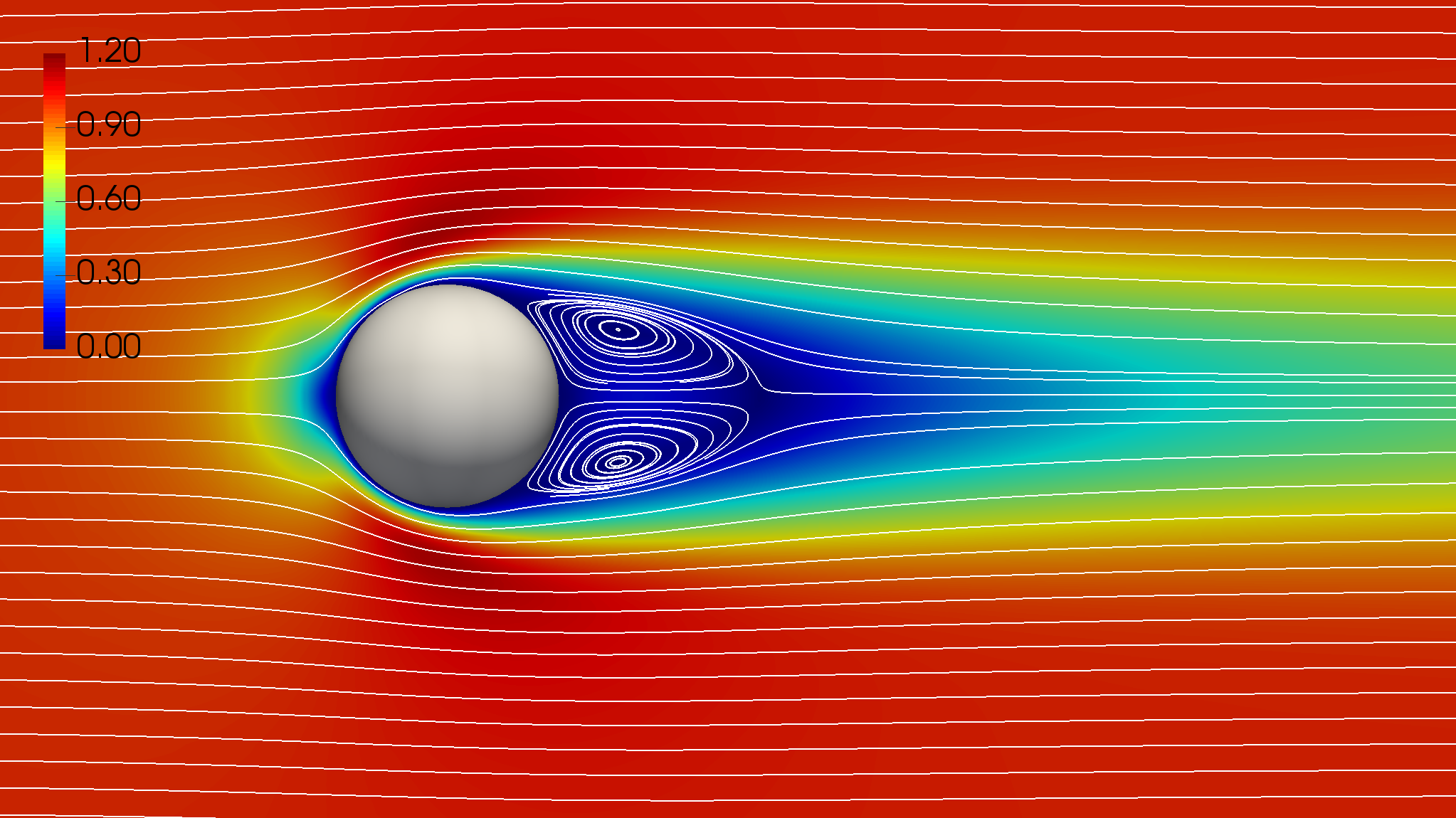} 
	\label{fig:SphereRe100:4pt}
}
\subfigure[Using $\delta_6^{\text{new}}$ weights in MLS]{
	\includegraphics[width=0.45\textwidth]{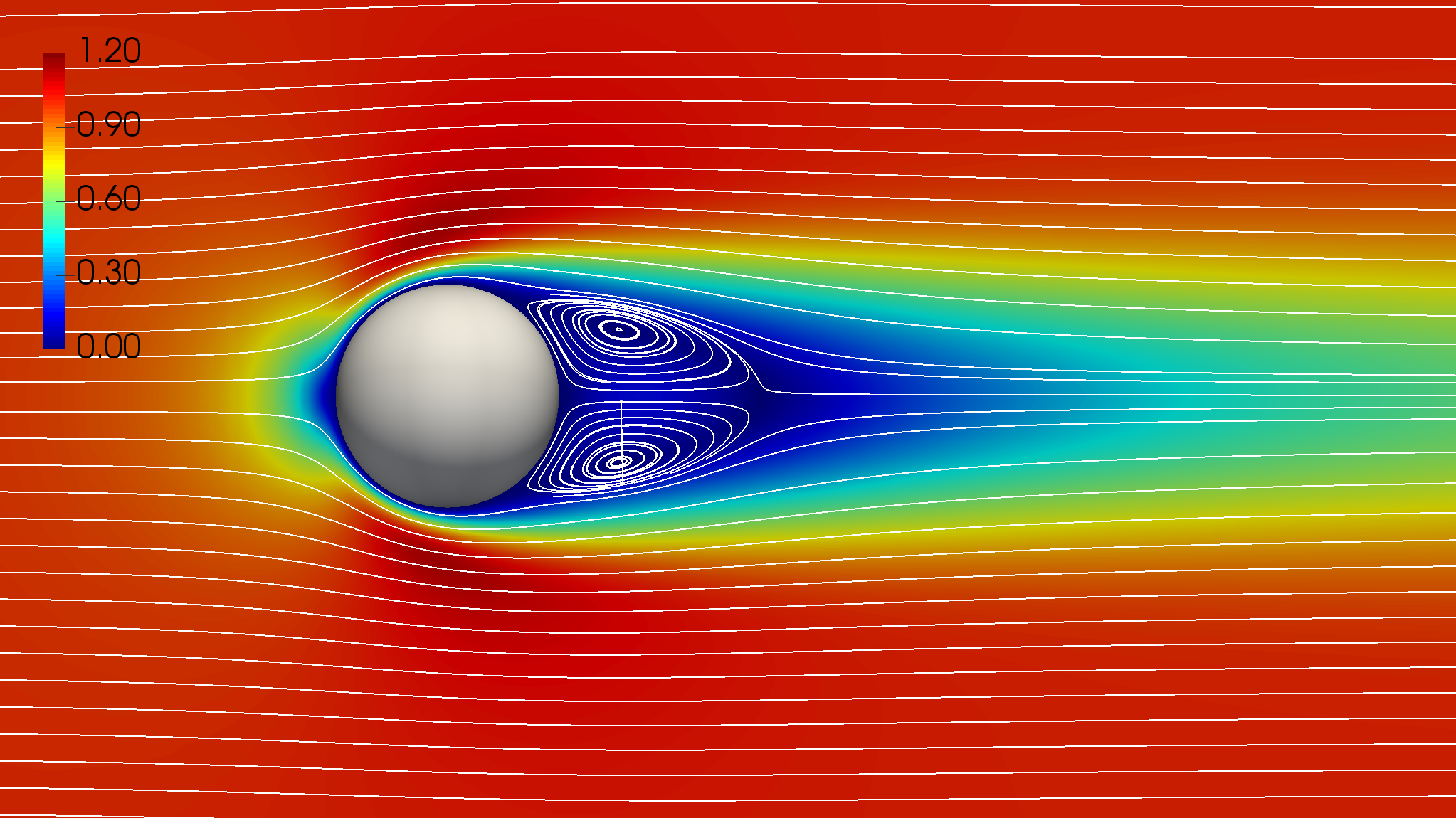} 
	\label{fig:SphereRe100:6pt}
}
\caption{Steady state streamlines and velocity magnitude for flow past a sphere at $Re=100$.}
\label{fig:sphere:Re100}
\end{figure}

{
\renewcommand{\arraystretch}{1.4}
\begin{table}[tb]
	\begin{center}
		\caption{\label{tab:spherewake} Comparison of steady state wake characteristics for flow over a sphere at $Re=100$ with data from the literature.}
		\vspace{0mm}  
		\begin{tabular}{llll}\hline
			& $y_c/D$ & $x_c/D$ & $L_b/D$  \\ \hline
			$\delta_{4}$        			                                     &  0.29   & 0.766 & 0.89  \\ 
			$\varphi_{5}$         						   &  0.29   & 0.76  & 0.9 \\
			$\varphi_{6}$         		 				   &  0.29   & 0.764  & 0.9 \\
			$\delta_5^{\text{new}}$          			   &  0.294   & 0.764  & 0.9 \\
			$\delta_6^{\text{new}}$         			   &  0.29   & 0.762  & 0.905 \\
			Mittal et al.~\cite{mitt08}    & 0.278  & 0.742 &  0.84 \\
			Johnson and Patel~\cite{john99}     & 0.29    & 0.75   & 0.88 \\
			Taneda~\cite{tane56}                      & 0.28    & 0.74   & 0.8 \\ 
			Tomboulides and Orszag~\cite{tomb00} &- &-& 0.88  \\ \hline
		\end{tabular} \\
	\end{center}
\end{table}
}

The comparison of geometric wake characteristics in Table~\ref{tab:spherewake} shows that all of the kernels reported yield virtually identical results. Based on this, it can be concluded that the IB kernel with the smallest width that is stable and accurate is the four-point IB kernel. In what follows, the four-point IB kernel will be used unless stated otherwise. With regards to the wake characteristics, the results presented in Table~\ref{tab:spherewake} show that our simulations agree well with the reported data in the literature. The streamlines and velocity magnitude of the flow in the vicinity of the sphere are shown in Fig.~\ref{fig:sphere:Re100}. 
Results obtained from $\delta_4$ and $\delta_6^{\text{new}}$ IB kernels are essentially indistinguishable. 

\begin{figure}
\centering
\subfigure[$C_d$]{
	\includegraphics[scale=0.24]{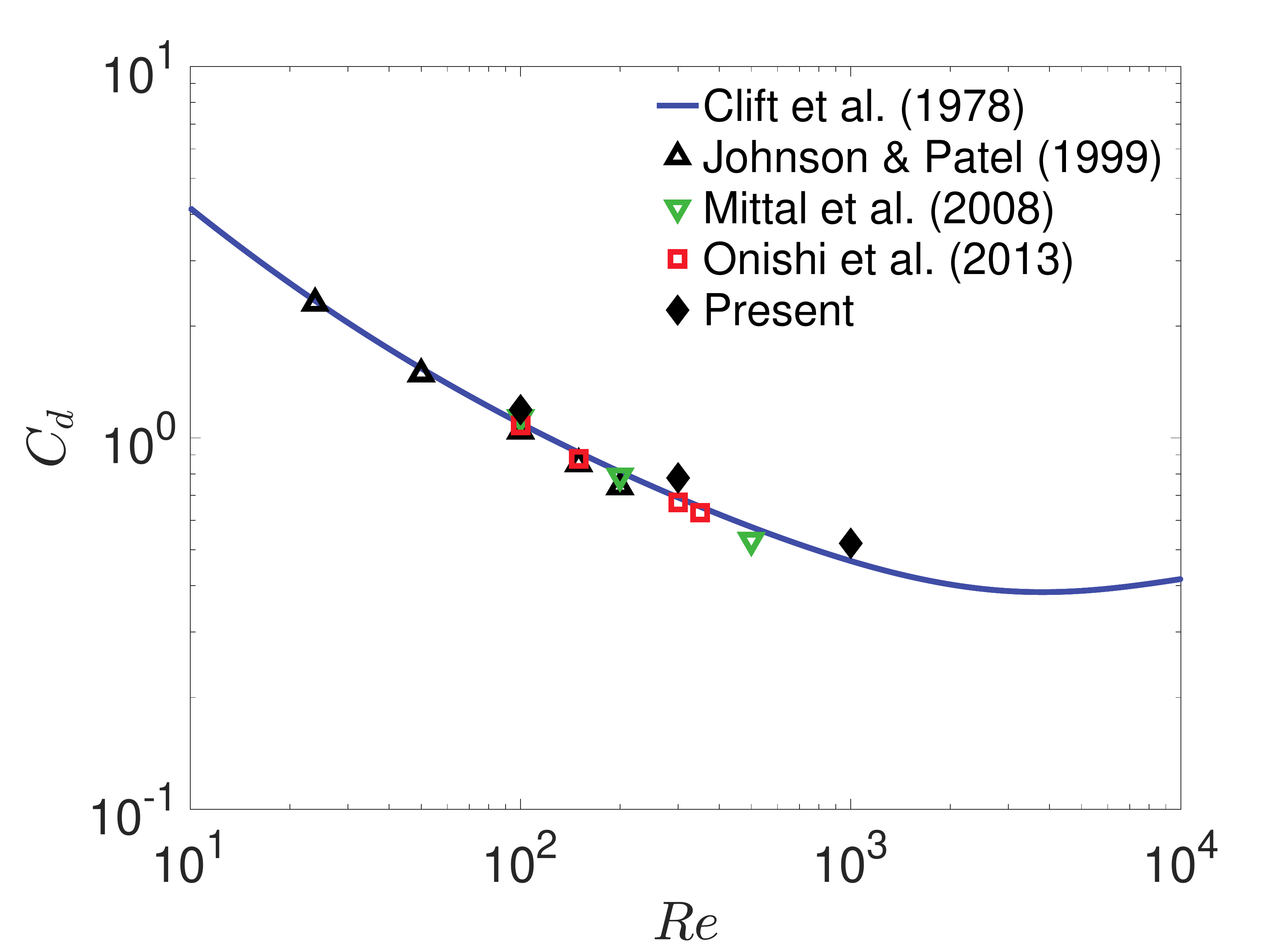} 
	\label{fig:sphere:Cd-Re}
}
\subfigure[$C_p$]{
	\includegraphics[scale=0.24]{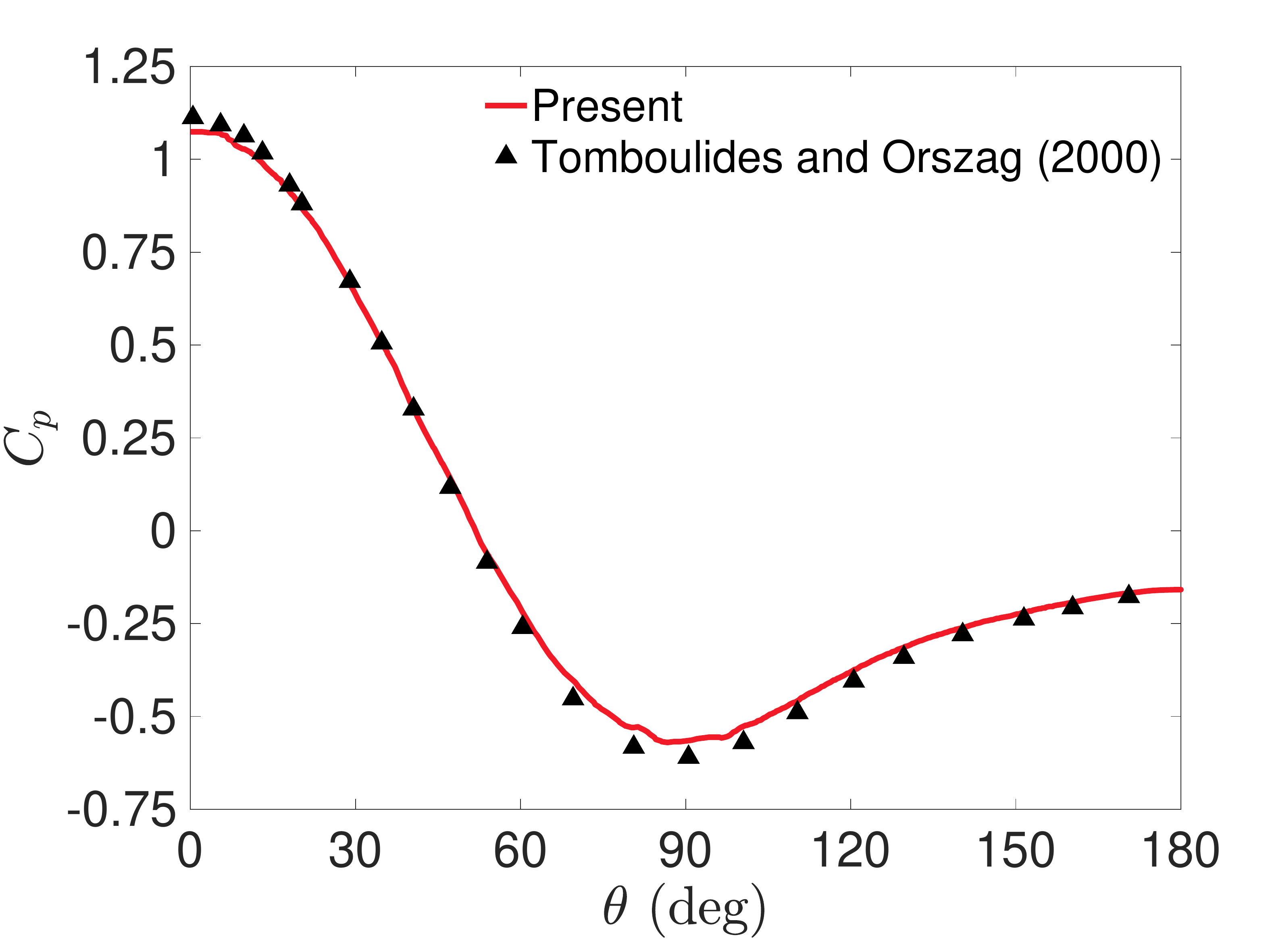} 
	\label{fig:sphere:Cp-Re100}
}
\caption{\subref{fig:sphere:Cd-Re} Time-averaged drag coefficient plotted against Reynolds number and compared with literature \cite{clif78,mitt08,john99,onis13}.~\subref{fig:sphere:Cp-Re100} Comparison of steady state surface pressure coefficient at $Re = 100$, plotted against polar angle and compared with the reported results of Tomboulides and Orszag~\cite{tomb00}. }
\label{fig:sphere:Cd-Cp}
\end{figure}

\begin{figure}
\centering
\subfigure[$\delta_4$]{
	\includegraphics[width=0.45\textwidth]{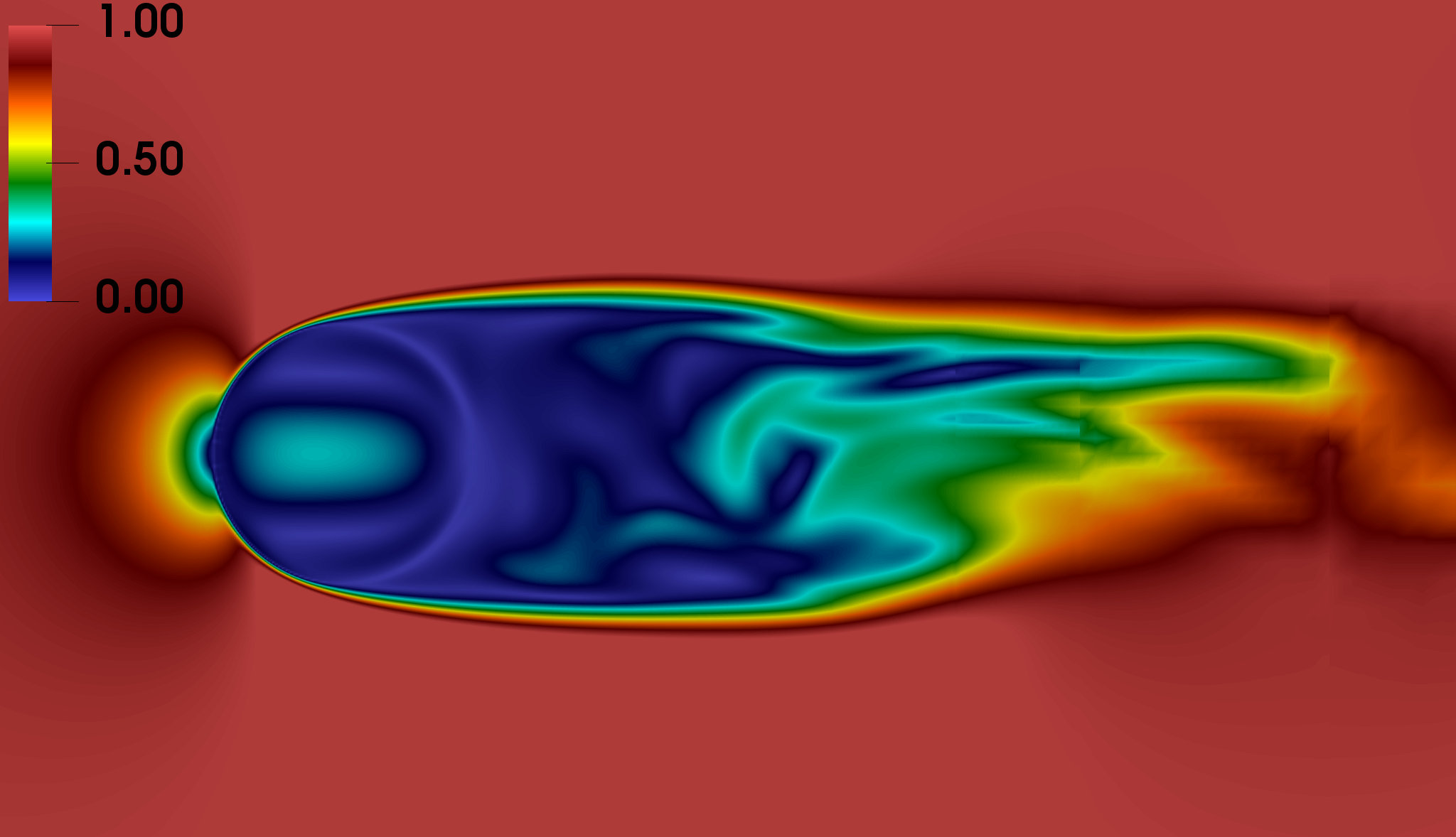} 
	\label{fig:sphere:Re1000-cIB}
}
\subfigure[$\V{\Psi}^m_{\text{NCVS}}$]{
	\includegraphics[width=0.45\textwidth]{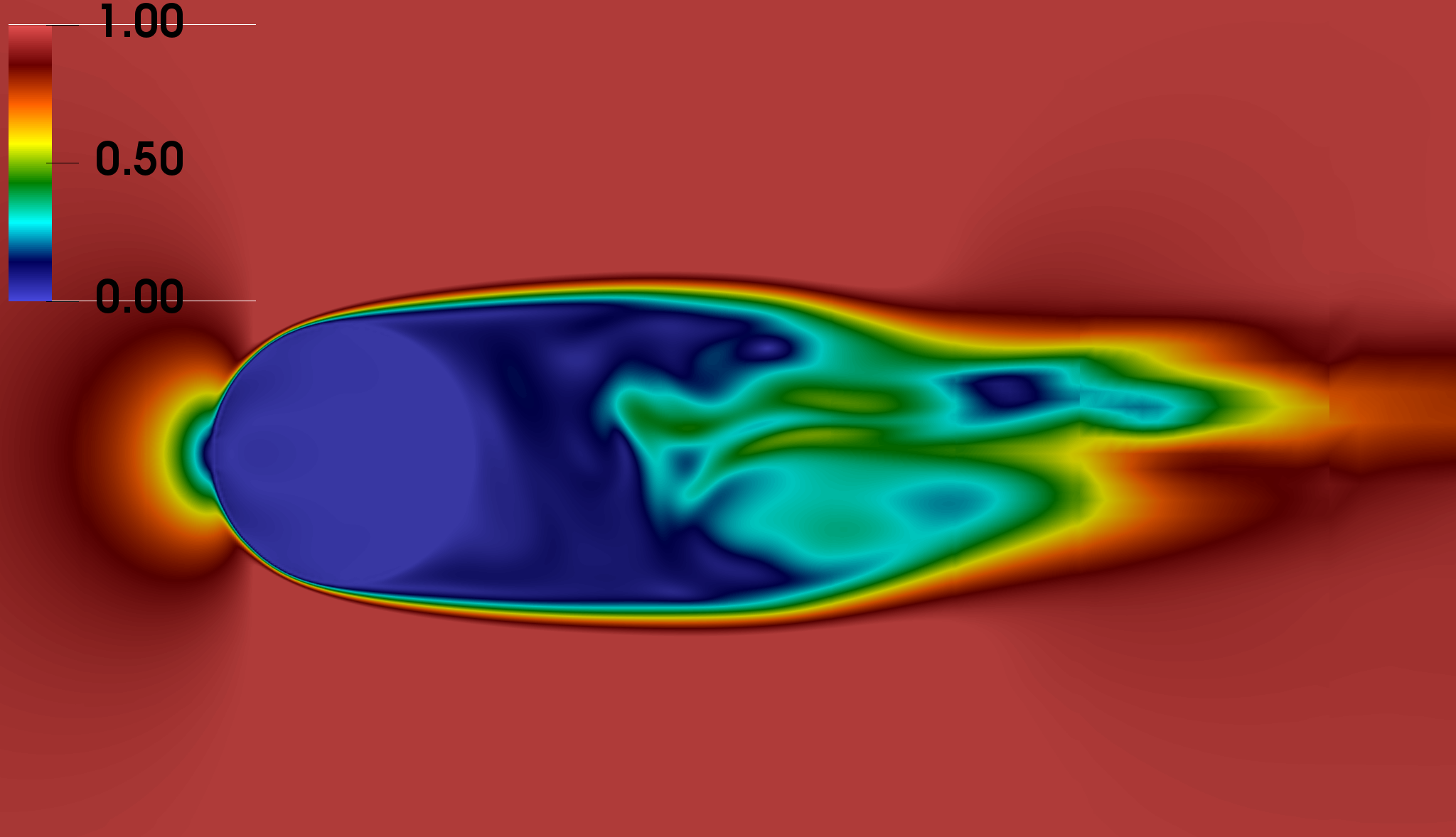} 
	\label{fig:sphere:Re1000-MLS}
}
\caption{Comparison of internal flow \REVIEW{(shown through velocity magnitude colorplot)} inside the sphere at $Re=1000$ using:~\subref{fig:sphere:Re1000-cIB} two-sided $\delta_4$ IB kernel; and ~\subref{fig:sphere:Re1000-MLS} one-sided $\V{\Psi}^m_{\text{NCVS}}$ MLS kernel for the direct forcing IB method at $ t^* = tU_\infty/D = 50 $, in which $U_\infty$ is the uniform inlet velocity.}
\label{fig:sphere:Re1000}
\end{figure}

{
\renewcommand{\arraystretch}{1.4}
\begin{table}[tb]
	\begin{center}
		\caption{\label{tab:sphereSt} Strouhal number of flow over sphere at $ Re = 300 $ and $ Re =1000 $.}
		\vspace{0mm}  
		\begin{tabular}{lll}\hline
			& $ Re = 300$ & $ Re =1000$ \\ \hline
			Present         										     &0.127 & 0.195 \\ 
			Tomboulides and Orszag~\cite{tomb00} & 0.136& 0.195 \\
			Mittal et al.~\cite{mitt08}   &0.137& -- \\
			Johnson and Patel~\cite{john99}    &0.135& -- \\  \hline
		\end{tabular} \\
	\end{center}
\end{table}
}

The time-averaged drag coefficients for the three $Re$ cases considered are plotted in Fig.\ref{fig:sphere:Cd-Re} and compared with experimental and numerical data from literature. A good agreement is obtained. The surface pressure coefficient for the $ Re=100 $ case, as shown in Fig.\ref{fig:sphere:Cp-Re100}, agrees very well with that reported in Tomboulides and Orszag~\cite{tomb00}, who used a body-fitted finite volume solver in their study.  The dimensionless dominant frequency 
of the unsteady wake, i.e. the Strouhal number of the flow for $Re=300$ and $Re=1000 $ cases are compared with the literature in Table~\ref{tab:sphereSt}. For both cases the agreement is good. We conclude this section with a comparison of velocity magnitude of the flow inside the sphere geometry at $ Re=1000 $ produced by using one-sided and regular IB kernels. As shown in Fig.~\ref{fig:sphere:Re1000} the two-sided four-point IB kernel leads to sloshing flow inside the sphere geometry, which can interact and influence the external flow through the discrete diffusion and advection operators. In contrast, the internal flow using the current one-sided IB kernel is insignificant. Hence, it results in a better separation of internal and external flow domains.

\subsection{Ahmed Body}

The study of the aerodynamics of road vehicles is challenging both numerically and experimentally due to the complexity of road vehicle geometries and the resulting high-speed flows. To enable characterization of the complex three-dimensional flow around road vehicles a simplified geometry, known as the Ahmed vehicle model or the Ahmed body has been used~\cite{ahme84}. The Ahmed vehicle model retains the main characteristic of the flow features around real vehicles such as the three-dimensional regions of separation and unsteady coherent vortex structures. Owing to its ability to mimic real vehicle flows, the Ahmed vehicle model has been the de-facto standard to study the road vehicle aerodynamics. 

\begin{figure}
	\begin{center}
		\includegraphics[width=0.8\textwidth]{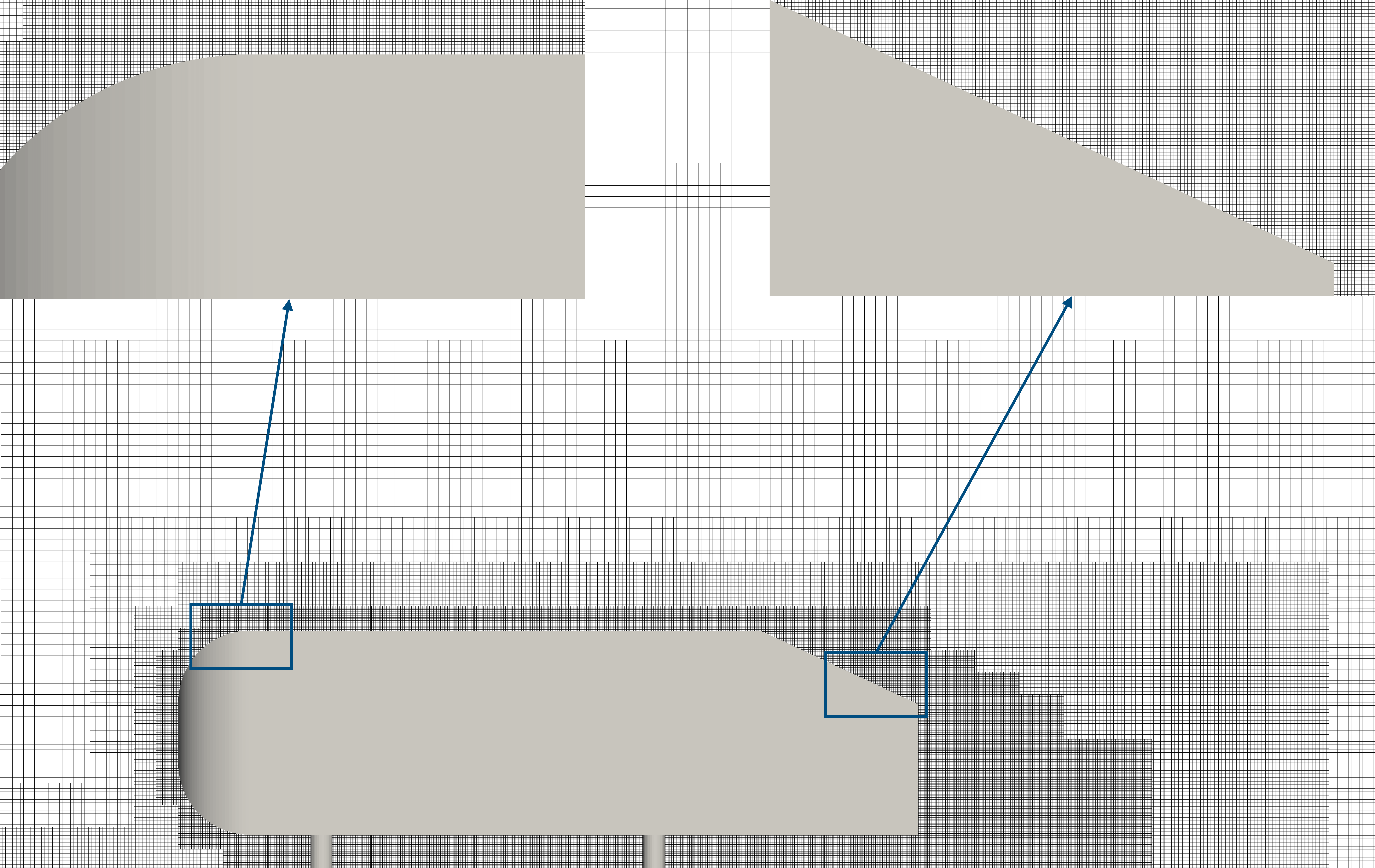} 
	\end{center}
	\caption{\label{fig:Abody-mesh} Ahmed vehicle model and the locally refined mesh used in the numerical simulation. The insets show a magnified view of the mesh near the geometry surface. }
\end{figure}

\begin{figure}
	\centering
	\subfigure[]{
	   \includegraphics[width=0.36\textwidth]{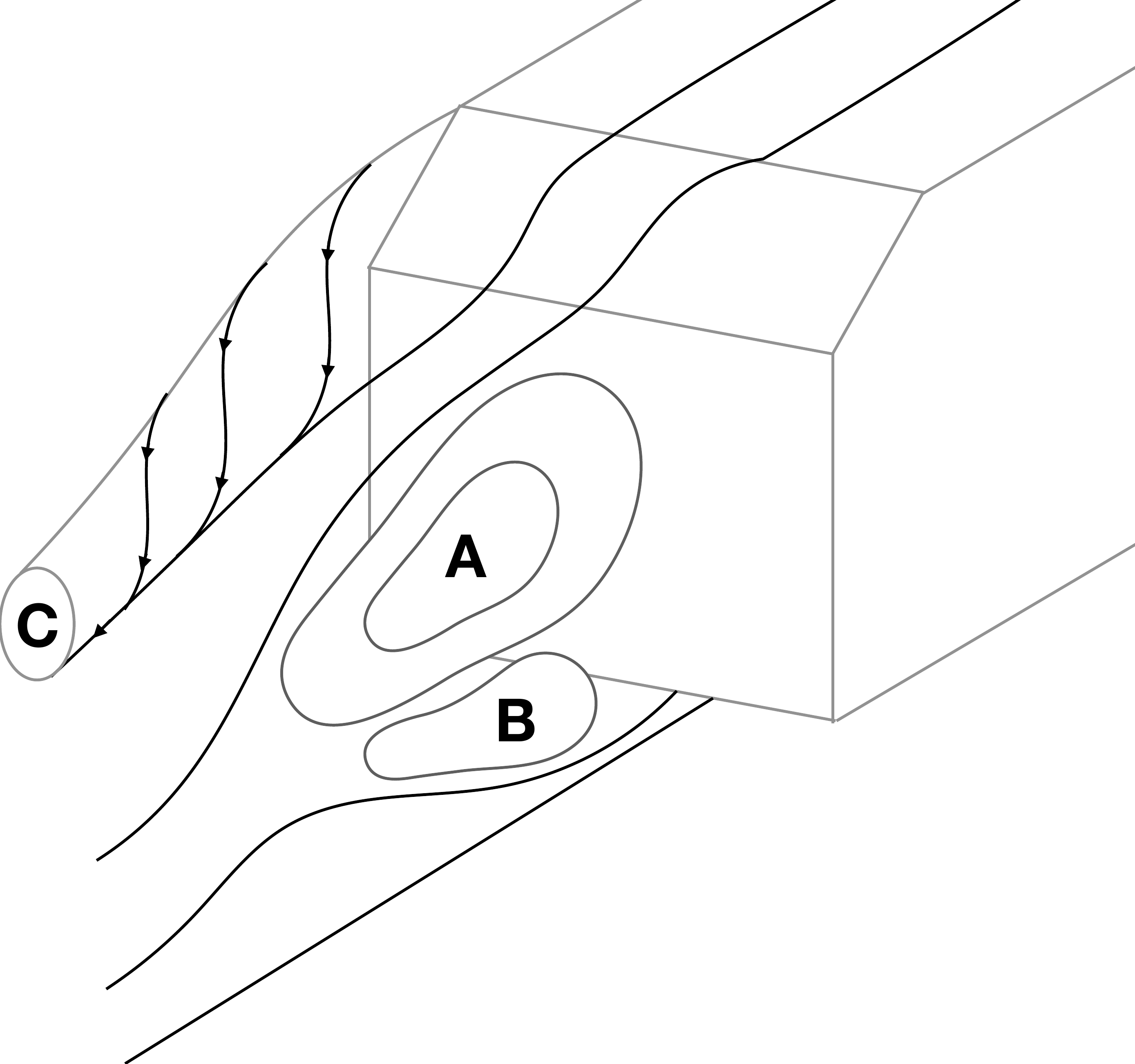} 
	   \label{fig:Ahmed:wake-schematic}
	}
	\subfigure[]{
	    \includegraphics[width=0.62\textwidth]{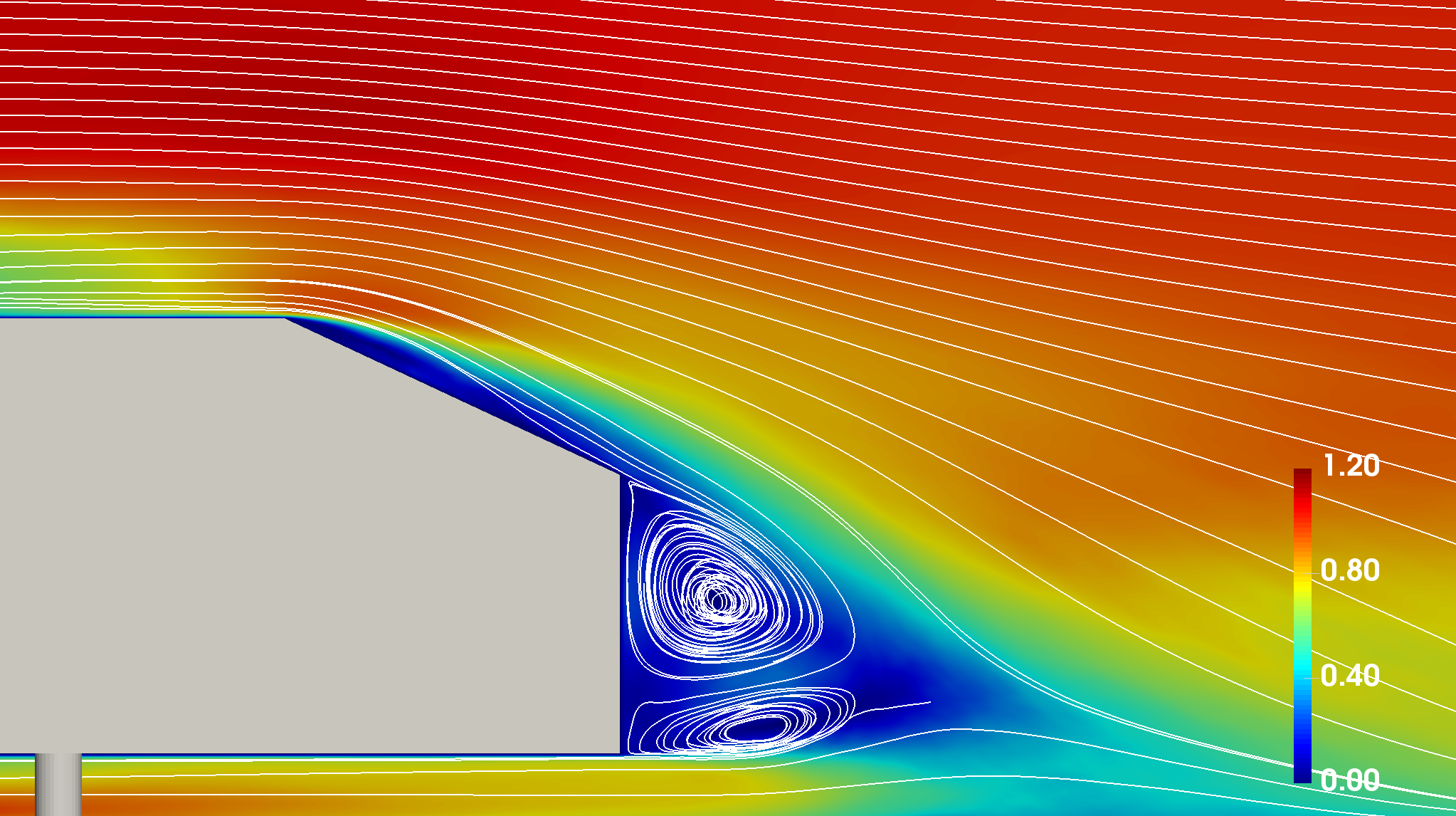} 
  	    \label{fig:Ahmed:wake-streamlines}
	}
	\subfigure[]{
	    \includegraphics[width=0.6\textwidth]{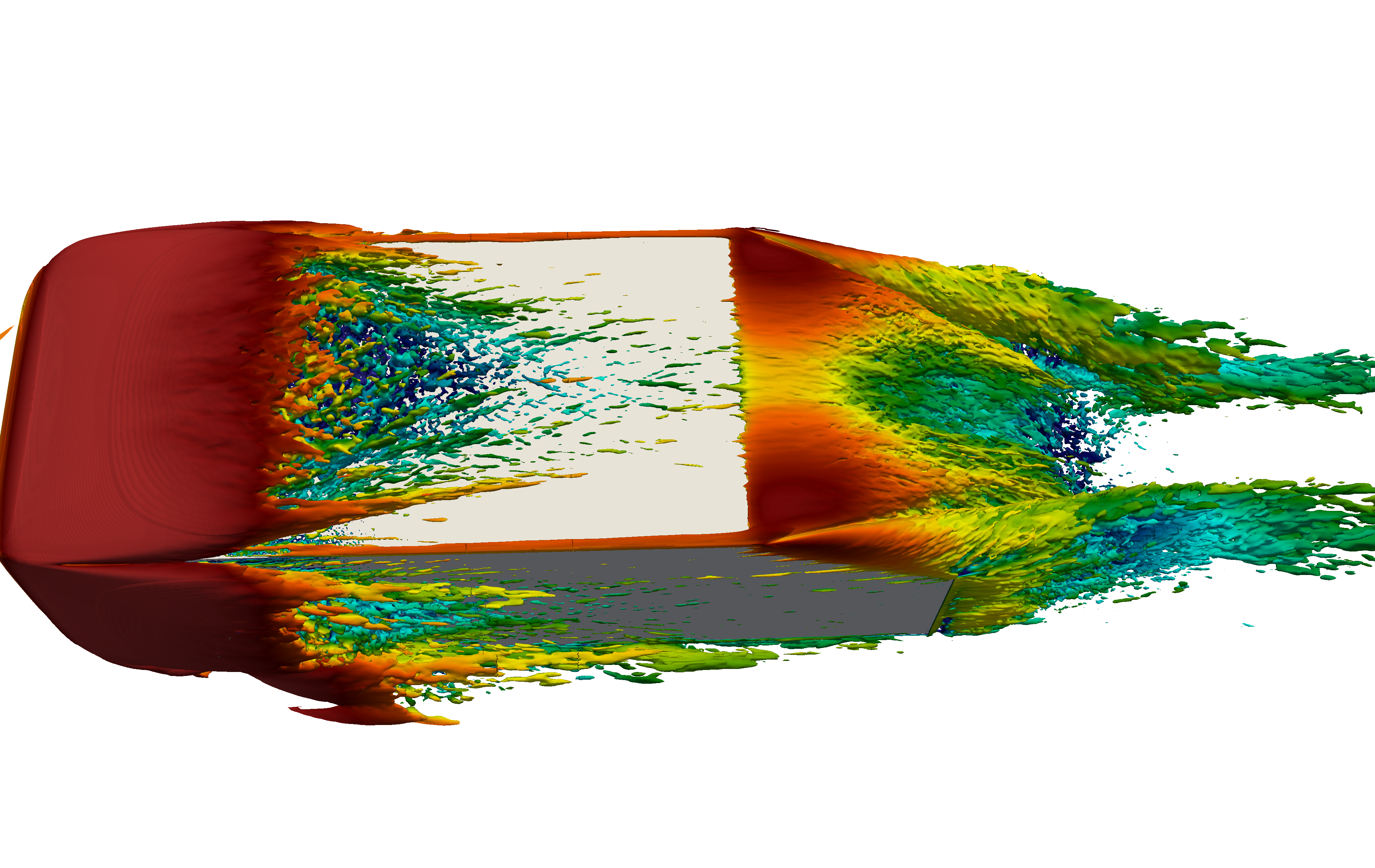} 
  	    \label{fig:Ahmed:Qcriteria}
	}
\caption{Flow past an Ahmed body at $Re_L  = 9.31\times 10^{4}$.~\subref{fig:Ahmed:wake-schematic} Schematic of time-averaged wake structures behind an Ahmed body when the slant angle is between $12.5^{\circ}$ and $ 30^{\circ}$\cite{ahme84,vino05};~\subref{fig:Ahmed:wake-streamlines} Visualization of the characteristic wake structures and the separation bubble over the slant surface of the Ahmed vehicle model through streamlines and velocity magnitude on a vertical plane;~\subref{fig:Ahmed:Qcriteria} Iso-surfaces of Q-criterion with dimensionless $Q^{*}=QL^2/U^2=10$ of the time-averaged flow.}
\label{fig:Abody-wakestructures}
\end{figure}

\begin{figure}
	\centering
	\subfigure[Two-sided Peskin $\delta_4$ kernel]{
	\includegraphics[width=0.6\textwidth]{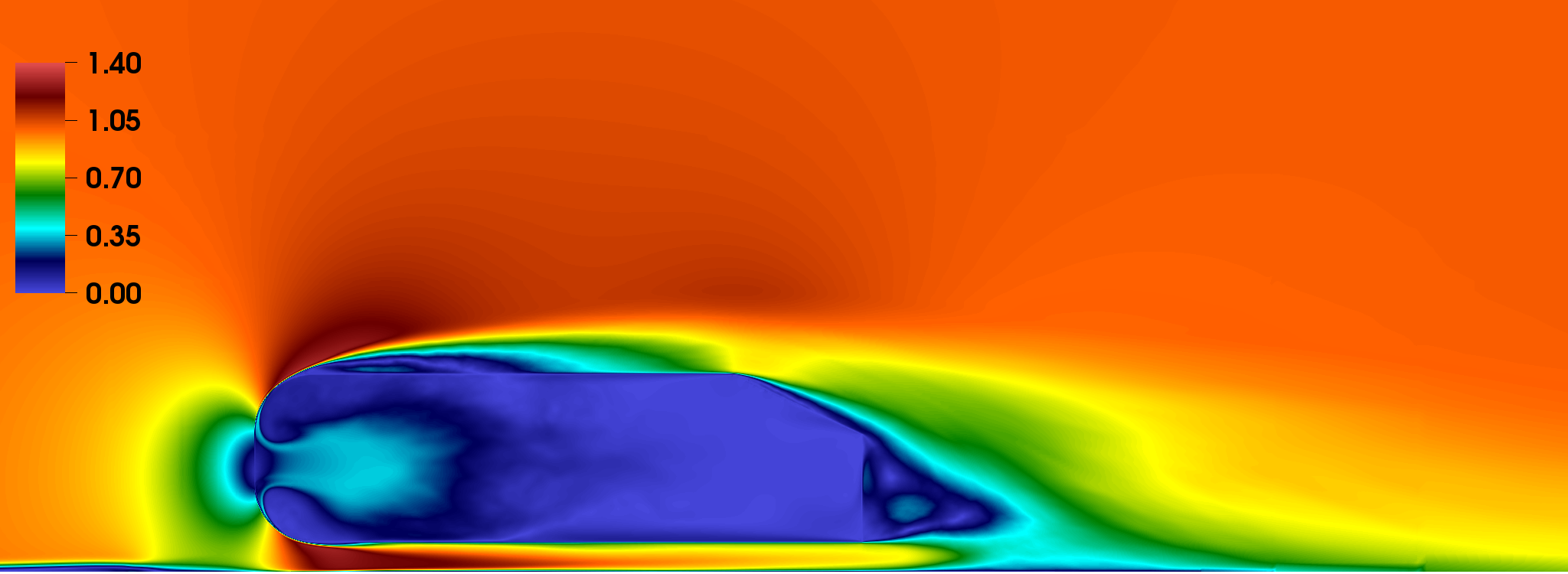}
	\label{fig:Abody-VMag:4pt} 
    } 
	\subfigure[One-sided $\V{\Psi}^m_{\rm NCVS}$ IB/MLS kernel]{
	\includegraphics[width=0.6\textwidth]{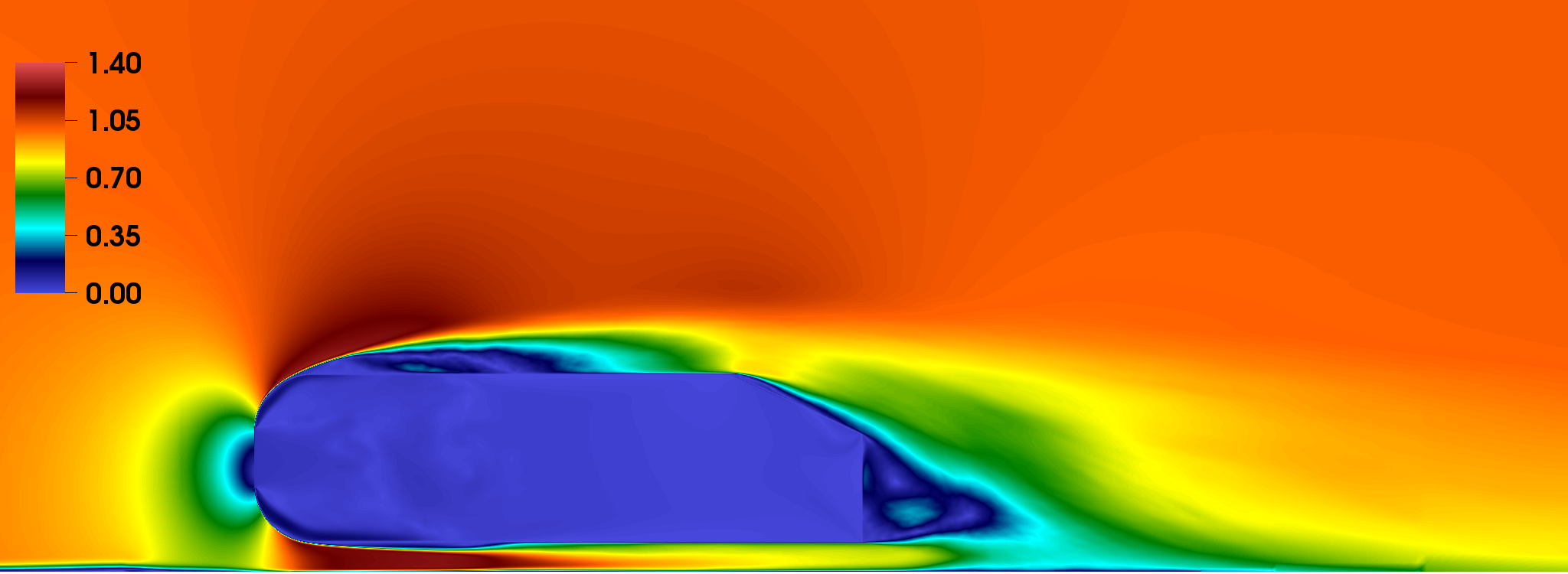}
	\label{fig:Abody-VMag:MLS} 
} 
	\caption{\REVIEW{Time-averaged velocity magnitude of flow over the Ahmed vehicle model using~\subref{fig:Abody-VMag:4pt} the standard two-sided  and~\subref{fig:Abody-VMag:MLS} the present one-sided IB method.}}
	\label{fig:Abody-VMag} 
\end{figure}

The applicability of the present method to moderately high $Re$ flows and complex geometries is demonstrated by using  the flow around an Ahmed body as the final validation case. The Ahmed vehicle model considered for present case has a rear slant surface that makes an angle of $ 25^{\circ}$ with the top horizontal surface. Our simulations are validated against the experimental data of Moghimi and Rafee~\cite{mogh18}. Moghimi and Rafee report drag and lift coefficient at $Re_L = 9.31\times 10^{4}$, in which the Reynolds number is based on the length $L$ of the Ahmed body.  The computational domain and the boundary conditions are set to match the experimental setup of~\cite{mogh18}. A uniform inflow boundary condition and a homogeneous Neumann boundary condition are imposed at the inflow and outflow boundaries, respectively, and a no-slip condition is imposed on the boundaries along the lateral and vertical direction.  To ensure that the boundary layer encountered over the Ahmed body is resolved by our computational mesh, a near geometry mesh resolution of $0.001L$ is employed. Local mesh refinement is used to restrict the fine mesh close to the geometry surface to reduce the computational cost of uniform grids. The locally refined mesh used for the present case is shown in Fig.~\ref{fig:Abody-mesh}. Considering the Blasius solution over a flat plate, the laminar boundary layer thickness is expected to scale as $\delta_{\text{lam}} \sim 5 L Re_{L}^{-\half}$,  which gives a boundary layer thickness of $\delta_{\text{lam}} \sim 0.016L$. The chosen near wall mesh resolution results in at least 15 cells inside the boundary layer.   The drag and lift coefficients from our results are compared with the wind tunnel data of~\cite{mogh18} in Table~\ref{tab:Abody}. The results of our simulation are in good agreement with the reported data of Moghimi and Rafee within the error range of the experimental data. 

{
	\renewcommand{\arraystretch}{1.4}
	\begin{table}[tb]
		\begin{center}
			\caption{\label{tab:Abody} Drag and lift coefficient of the Ahmed body.}
			\vspace{0mm}  
			\begin{tabular}{lll}\hline
				& $ C_d $& $C_L$ \\ \hline
				Present         									& 0.4 & 0.29 \\ 
				Moghimi and Rafee~\cite{mogh18}   & 0.44$ \pm $0.04& 0.25$ \pm $0.03 \\
				\hline
			\end{tabular} \\
		\end{center}
	\end{table}
}

Although the wake of the flow over a bluff-body is typically unsteady, the time-averaged flow over an Ahmed body shows persistent large-scale wake structures. Several experimental studies have characterized these flow structures based on the rear slant angle~\cite{ahme84, vino05, spoh02}. These studies have identified critical slant angles of $ 12.5^{\circ} $ and $ 30^{\circ}$ which appear to govern the topological structures of the time-averaged flow. When the slant angle is between these critical angles, three predominant three-dimensional vortical structures occur in the wake accompanied by a separation bubble on the slant surface. A schematic of these structures is shown in Fig.~\ref{fig:Ahmed:wake-schematic}. Right behind the rear end of the Ahmed body, two counter-rotating vortices A and  B, one larger than the other can be found. The third vortical structure C originates at the intersection of the slant surface and the side wall, which is stretched and elongated in the wake while it feeds the two primary vortices A and B. The characteristic wake signatures of the Ahmed body are identified in the results of our numerical simulation for qualitative validation. The near-wake time-averaged flow and streamlines on a vertical plane are presented in Fig.~\ref{fig:Ahmed:wake-streamlines} in which the primary vortices A and B can be identified and the separation bubble over the rear slant surface can also be seen. The elongated helical flow structure C can be identified through the visualization of the time-averaged Q-criterion. In Fig.~\ref{fig:Ahmed:Qcriteria}, iso-surfaces of Q-criterion with dimensionless $Q^{*} = 10$ are plotted. The helical wake structure C can be seen on the two sides of the rear end of the Ahmed body. \REVIEW{Fig.~\ref{fig:Abody-VMag} shows the side view of the velocity magnitude over the entire Ahmed car model. It is clearly observed in the figure that the two-sided IB method leads to substantial flows within the car interior, which are largely eliminated using the current one-sided approach.}


\section{Summary and Conclusions}
This work develops a one-sided IB kernel approach for diffuse-interface IB methods. For IB models using thin structure or interface representations, the two sides of the flow region often require segregation and separate treatments. This is necessary to avoid spurious flows inside closed geometries and to avoid interaction of flow on the two sides of an IB surface resulting from velocity interpolation and force spreading operations. To realize this, this paper introduces an approach that uses the moving least squares (MLS) method to dynamically generate one-sided IB kernels from  standard two-sided IB kernels. We find that the stability of direct forcing IB methods benefit from kernels that monotonically decreasing and positive (or have negligible negative tails) to avoid oscillatory feedback through force spreading. A straightforward application of the one-sided MLS construction can generate larger weights for nearby Eulerian grid nodes and negative weights for far-away grid nodes to satisfy the linear conditions. Two weight-shifting approaches, NCVS and CVS, are proposed to alleviate this issue. It was shown that the CVS approach can generate numerically thick boundary layers and spurious flow oscillations near the  interface, whereas the NCVS approach produced physically correct flow features around the structure. Therefore, based on these results, the NCVS mollification approach is recommended.  

The order of accuracy of one-sided IB kernels is tested through the Taylor-Green vortex flow problem. It was shown that by employing mollified weights in both velocity interpolation and force spreading operators the order of accuracy of the solution reduces to one. However, the accuracy of the scheme is improved (becomes super-linear) by adopting the original non-shifted MLS weights in the velocity interpolation operator. For more general problems, which typically include stress discontinuities along the fluid-structure interface, the present method will be only first-order accurate even with non-shifted MLS weights in the interpolation operator. This also holds true for regular diffuse-interface IB methods.  

The dependency of accuracy and stability of an MLS kernel on the underlying weighting function is also analyzed. The one-sided MLS kernel was found to be stable for all weighting functions tested in this work except for the narrower three-point kernel. Our observation is that a basic regularized delta function kernel appears to require at least a support of four grid cells in each coordinate direction. A wider kernel avoids (near) singular Gram matrix in the one-sided moving least squares problem. When considering the problem of flow past a sphere at $Re=100$, all weighting kernels were found to produce essentially equally accurate one-sided kernels.  The proposed IB/MLS method was further validated through the cases of flow past an oscillating cylinder, Stokes' first problem, flow over a sphere, impulsively started plate/cylinder, and the Ahmed vehicle model. In all cases, our results are in excellent agreement with the reported data in the literature. We also demonstrate that the present method can effectively eliminate the spurious internal flow that is typically generated by simulating flow past bluff bodies with regular diffuse-interface IB methods.


\section*{Acknowledgements}
A.P.S.B~acknowledges support from NSF award OAC 1931368. R.B.~acknowledges the support for this work by JSPS KAKENHI Grant Number 20K19503.  B.E.G~acknowledges support from NSF awards OAC 1450327, OAC 1652541, OAC 1931516, and DMS 1664645.

\appendix
\renewcommand\thesection{\Alph{section}}
\section{Kernel functions} \label{sec_kernel_appendix}

Here we write the one-dimensional form of some of the kernels used in this work. In multiple dimensions, a tensor product of the one-dimensional kernels is formed. The functional form of the kernel is expressed in terms of $r = (x - X)/h$.

\begin{align}
&\text{Smoothed three-point IB kernel~\cite{Yang09}}: \delta_{3}(r) = \begin{cases}
  \frac{3}{4} - r^2, & 0 \le |r| < 0.5, \\
    \frac{1}{2}(\frac{9}{4} -{3|r|} + {r^2}), & 0.5 \le |r| < 1.5,  \\
       0, & 1.5 \le |r|.
\end{cases} \\
&\text{Peskin's four-point IB kernel~\cite{Peskin02}}: \delta_{4}(r) =  \begin{cases}
\frac{1}{8} \left( 3 - 2|r| + \sqrt{1 + 4|r| - 4r^2}\right), & 0 \le |r| < 1, \\
\frac{1}{8} \left( 5 - 2|r| - \sqrt{-7 + 12|r| - 4r^2}\right), & 1 \le |r| < 2,  \\
0, & 2 \le |r|.
\end{cases} \\
&\text{Radial basis function}: \text{RBF}(r) =  \begin{cases}
e^{-2r^2}, & |r| < 2,  \\
0, & 2 \le |r|.
\end{cases} \\
&\text{Two-point cubic spline kernel~\cite{Vanella09}}: \varphi_{2}(\bar{r})  =  \begin{cases}
\frac{2}{3} -4 \bar{r}^2 + 4 \bar{r}^3 & 0 \le |\bar{r}| < 0.5, \\
\frac{4}{3} -4\bar{r} + 4 \bar{r}^2 - \frac{4}{3}\bar{r}^3 & 0.5 \le |\bar{r}| < 1,  \\
0, & 1 \le |\bar{r}|.
\end{cases} \\
&\text{Five-point spline kernel}: \varphi_{5}(r)  =  \begin{cases}
\frac{1}{24}( 6 \kappa^4 -60 \kappa^3 + 210\kappa^2 - 300\kappa+155) & 0 \le |r| < 0.5, \\
\frac{1}{24}( -4\kappa^4 +60 \kappa^3 - 330\kappa^2 + 780\kappa + 655) & 0.5 \le |r| < 1.5,  \\
\frac{1}{24}( \kappa^4 -20 \kappa^3 + 150\kappa^2 - 500\kappa +625 ) & 1.5 \le |r| < 2.5,  \\
0, & 2.5 \le |r|,
\end{cases} \\
&\text{Six-point spline kernel}: \varphi_{6}(r)  =  \begin{cases}
\frac{1}{60}( -5\kappa^5 +90\kappa^4 -630 \kappa^3 + 2130\kappa^2 - 3645\kappa+2193) & 0 \le |r| < 1, \\
\frac{1}{120}( 5\kappa^5 -120\kappa^4 +1140 \kappa^3 - 5340\kappa^2 + 12270\kappa - 10974) & 1 \le |r| < 2,  \\
\frac{1}{120}( -\kappa^5 +30\kappa^4 -360 \kappa^3 + 2160\kappa^2 - 6480\kappa+7776) & 2 \le |r| < 3,  \\
0, & 3 \le |r|.
\end{cases} 
\end{align}
Here $\bar{r} = r/1.2$ for the two-point (more specifically 2.4 grid cells wide) cubic spline function used in Vanella and Balaras~\cite{Vanella09},  and $ \kappa = |r|+2.5 $ and $ \kappa = |r| +3 $ for the five-point and six-point spline function, respectively. 
In addition to the aforementioned kernels, we also consider the new five- and six-point IB kernel, $\delta_5^{\text{new}}$ and $\delta_6^{\text{new}}$, respectively. 
These new kernels remove the negative tail of the standard five- and six-point IB kernels by imposing a weaker second moment condition. We refer readers to Bao et al.~\cite{Bao2016} for their functional form.

\section*{Bibliography}
\begin{flushleft}
 \bibliography{References}
\end{flushleft}

\end{document}